\documentclass[12pt]{amsart}
\usepackage{amscd,amsmath,amssymb,amsfonts}
\usepackage[cmtip, all]{xy}
\begin{document}
\newtheorem{thm}{Theorem}[subsection]
\newtheorem{prop}[thm]{Proposition}
\newtheorem{lem}[thm]{Lemma}
\newtheorem{cor}[thm]{Corollary}
\newtheorem{conj}[thm]{Conjecture}
\theoremstyle{definition}
\newtheorem{Def}[thm]{Definition}
\theoremstyle{remark}
\newtheorem{rem}[thm]{Remark}
\newtheorem{rems}[thm]{Remarks}
\newtheorem{ex}[thm]{Example}
\newtheorem{exs}[thm]{Examples}
\numberwithin{equation}{subsection}

\newcommand{\sA}{{\mathcal A}}
\newcommand{\sB}{{\mathcal B}}
\newcommand{\sC}{{\mathcal C}}
\newcommand{\sD}{{\mathcal D}}
\newcommand{\sE}{{\mathcal E}}
\newcommand{\sF}{{\mathcal F}}
\newcommand{\sG}{{\mathcal G}}
\newcommand{\sH}{{\mathcal H}}
\newcommand{\sI}{{\mathcal I}}
\newcommand{\sJ}{{\mathcal J}}
\newcommand{\sK}{{\mathcal K}}
\newcommand{\sL}{{\mathcal L}}
\newcommand{\sM}{{\mathcal M}}
\newcommand{\sN}{{\mathcal N}}
\newcommand{\sO}{{\mathcal O}}
\newcommand{\sP}{{\mathcal P}}
\newcommand{\sQ}{{\mathcal Q}}
\newcommand{\sR}{{\mathcal R}}
\newcommand{\sS}{{\mathcal S}}
\newcommand{\sT}{{\mathcal T}}
\newcommand{\sU}{{\mathcal U}}
\newcommand{\sV}{{\mathcal V}}
\newcommand{\sW}{{\mathcal W}}
\newcommand{\sX}{{\mathcal X}}
\newcommand{\sY}{{\mathcal Y}}
\newcommand{\sZ}{{\mathcal Z}}
\newcommand{\A}{{\mathbb A}}
\newcommand{\B}{{\mathbb B}}
\newcommand{\C}{{\mathbb C}}
\newcommand{\D}{{\mathbb D}}
\newcommand{\E}{{\mathbb E}}
\newcommand{\F}{{\mathbb F}}
\newcommand{\G}{{\mathbb G}}
\renewcommand{\H}{{\mathbb H}}
\newcommand{\I}{{\mathbb I}}
\newcommand{\J}{{\mathbb J}}
\newcommand{\M}{{\mathbb M}}
\newcommand{\N}{{\mathbb N}}
\renewcommand{\P}{{\mathbb P}}
\newcommand{\Q}{{\mathbb Q}}
\newcommand{\R}{{\mathbb R}}
\newcommand{\mS}{{\mathbb S}}
\newcommand{\T}{{\mathbb T}}
\newcommand{\U}{{\mathbb U}}
\newcommand{\V}{{\mathbb V}}
\newcommand{\W}{{\mathbb W}}
\newcommand{\X}{{\mathbb X}}
\newcommand{\Y}{{\mathbb Y}}
\newcommand{\Z}{{\mathbb Z}}

\renewcommand{\1}{{\mathbf{1}}}

\newcommand{\an}{{\rm an}}
\newcommand{\alg}{{\rm alg}}
\newcommand{\cl}{{\rm cl}}
\newcommand{\Alb}{{\rm Alb}}
\newcommand{\CH}{{\rm CH}}
\newcommand{\mc}{\mathcal}
\newcommand{\mb}{\mathbb}
\newcommand{\surj}{\twoheadrightarrow}
\newcommand{\red}{{\rm red}}
\newcommand{\codim}{{\rm codim}}
\newcommand{\rank}{{\rm rank}}
\newcommand{\Pic}{{\rm Pic}}
\newcommand{\Div}{{\rm Div}}
\newcommand{\Hom}{{\rm Hom}}
\newcommand{\im}{{\rm im}}
\newcommand{\Spec}{{\rm Spec \,}}
\newcommand{\sing}{{\rm sing}}
\newcommand{\Char}{{\rm char}}
\newcommand{\Tr}{{\rm Tr}}
\newcommand{\Gal}{{\rm Gal}}
\newcommand{\Min}{{\rm Min \ }}
\newcommand{\Max}{{\rm Max \ }}
\newcommand{\supp}{{\rm supp}\,}
\newcommand{\0}{\emptyset}
\newcommand{\sHom}{{\mathcal{H}{om}}}

\newcommand{\Nm}{{\operatorname{Nm}}}
\newcommand{\NS}{{\operatorname{NS}}}
\newcommand{\id}{{\operatorname{id}}}
\newcommand{\Zar}{{\text{\rm Zar}}} 
\newcommand{\Ord}{{\mathbf{Ord}}}
\newcommand{\FSimp}{{{\mathcal FS}}}
\newcommand{\inj}{{\text{\rm inj}}}  
\newcommand{\Sch}{{\operatorname{\mathbf{Sch}}}} 
\newcommand{\cosk}{{\operatorname{\rm cosk}}} 
\newcommand{\sk}{{\operatorname{\rm sk}}} 
\newcommand{\subv}{{\operatorname{\rm sub}}}
\newcommand{\bary}{{\operatorname{\rm bary}}}
\newcommand{\Comp}{{\mathbf{SC}}}
\newcommand{\IComp}{{\mathbf{sSC}}}
\newcommand{\Top}{{\mathbf{Top}}}
\newcommand{\holim}{\mathop{{\rm holim}}}
\newcommand{\op}{{\text{\rm op}}}
\newcommand{\<}{\mathopen<}
\renewcommand{\>}{\mathclose>}
\newcommand{\Sets}{{\mathbf{Sets}}}
\newcommand{\del}{\partial}
\newcommand{\fib}{{\operatorname{\rm fib}}}
\renewcommand{\max}{{\operatorname{\rm max}}}
\newcommand{\bad}{{\operatorname{\rm bad}}}
\newcommand{\Spt}{{\mathbf{Spt}}}
\newcommand{\Spc}{{\mathbf{Spc}}}
\newcommand{\Sm}{{\mathbf{Sm}}}
\newcommand{\cofib}{{\operatorname{\rm cofib}}}
\newcommand{\hocolim}{\mathop{{\rm hocolim}}}
\newcommand{\colim}{\mathop{{\rm colim}}}
\newcommand{\Glu}{{\mathbf{Glu}}}
\newcommand{\can}{{\operatorname{\rm can}}}
\newcommand{\Ho}{{\mathbf{Ho}}}
\newcommand{\GL}{{\operatorname{\rm GL}}}
\newcommand{\sq}{\square}
 \newcommand{\Ab}{{\mathbf{Ab}}}
\newcommand{\Tot}{{\operatorname{\rm Tot}}}
\newcommand{\loc}{{\operatorname{\rm s.l.}}}
\newcommand{\HZ}{{\operatorname{\sH \Z}}}
\newcommand{\Cyc}{{\operatorname{\rm Cyc}}}
\newcommand{\cyc}{{\operatorname{\rm cyc}}}
 \newcommand{\RCyc}{{\operatorname{\widetilde{\rm Cyc}}}}
 \newcommand{\Rcyc}{{\operatorname{\widetilde{\rm cyc}}}}
\newcommand{\Sym}{{\operatorname{\rm Sym}}}
\newcommand{\fin}{{\operatorname{\rm fin}}}
\newcommand{\SH}{{\operatorname{\sS\sH}}}
\newcommand{\Anna}{{$\mathcal{ANNA\ SARAH\ LEVINE}$}}
\newcommand{\Wedge}{{\Lambda}}
\newcommand{\eff}{{\operatorname{\rm eff}}}
\newcommand{\rcyc}{{\operatorname{\rm rev}}}
\newcommand{\DM}{{\operatorname{\mathcal{DM}}}}
\newcommand{\GW}{{\operatorname{\rm{GW}}}}
\newcommand{\sSets}{{\mathbf{sSets}}}
\newcommand{\Nis}{{\operatorname{Nis}}}
\newcommand{\Cat}{{\mathbf{Cat}}}
\newcommand{\ds}{{/\kern-3pt/}}
\newcommand{\bD}{{\mathbf{D}}}
\newcommand{\res}{{\operatorname{res}}}

\newcommand{\qf}{{\operatorname{q.fin.}}}

\newcommand{\fil}{\phi}
\newcommand{\barfil}{\sigma}

\newcommand{\Fac}{{\mathop{\rm{Fac}}}}
\newcommand{\Fun}{{\mathbf{Func}}}
\newcommand{\Th}{{\mathop{\rm{Th}}}}
\newcommand{\hofib}{{\mathop{\rm{hofib}}}}
\newcommand{\Aff}{{\mathop{\rm{Aff}}}}

\addtocounter{section}{-1}
\title{The homotopy coniveau tower}

\author{Marc Levine}

\address{
Department of Mathematics\\
Northeastern University\\
Boston, MA 02115\\
USA}

\email{marc@neu.edu}


\thanks{The author gratefully acknowledges the support of the Humboldt Foundation through the Wolfgang Paul Program, and support of the NSF via grants DMS 0140445 and DMS-0457195}

\renewcommand{\abstractname}{Abstract}
\begin{abstract} We examine the ``homotopy coniveau tower" for a general cohomology theory on smooth $k$-schemes and give a new proof that the layers of this tower for $K$-theory agree with motivic cohomology.  In addition, the homotopy coniveau tower agrees with Voevodsky's slice tower for $S^1$-spectra, giving a proof of a connectedness conjecture of Voevodsky.

The homotopy coniveau tower construction extends to a tower of functors on the Morel-Voevodsky stable homotopy category, and we identify this $\P^1$-stable homotopy coniveau tower with Voevodsky's slice filtration for $\P^1$-spectra. We also show that the 0th layer for the motivic sphere spectrum is the motivic cohomology spectrum, which gives the layers for a general $\P^1$-spectrum the structure of a module over motivic cohomology. This recovers and extends recent results of Voevodsky on the 0th layer of the slice filtration, and yields a spectral sequence that is reminiscent of the classical Atiyah-Hirzebruch spectral sequence.\\
\end{abstract}

\maketitle
\tableofcontents
\section{Introduction} The original purpose of this paper was to give an alternative argument
for the technical underpinnings of the papers \cite{BL, FriedSus}, in which the constructruction
of a spectral sequence from motivic cohomology to 	$K$-theory is given.  As in the method used by Suslin  \cite{Suslin} to analyze the Grayson spectral sequence, we rely on localization properties of the relevant spectra.

Having done this, it becomes clear that the method applies more generally to a
presheaf of spectra on  smooth schemes over a perfect base field $k$, satisfying certain
conditions.  We therefore give a general discussion for a presheaf of spectra $E:\Sm/k^\op\to\Spt$ which is homotopy invariant and satisfies Nisnevic excision.

For such a functor, and an $X$ in $\Sm/k$, we construct the {\em homotopy coniveau tower}
\begin{equation}\label{eqn:HCTower0}
\ldots\to E^{(p+1)}(X,-)\to  E^{(p)}(X,-)\to\ldots\to E^{(0)}(X,-)\sim E(X)
\end{equation}
where the $E^{(p)}(X,-)$ are simplicial spectra with $n$-simplices the limit of the spectra with support $E^W(X\times\Delta^n)$, where $W$ is a closed subset of codimension $\ge p$ in ``good position". This is just the evident extension of the tower used by Friedlander-Suslin in \cite{FriedSus}. One can consider this tower as the algebraic analog of the one in topology formed by applying a cohomology theory to the skeletal filtration of a CW complex. The main objects of our study are the layers $E^{(p/p+1)}(X,-)$ in this tower. 

We have discussed the functoriality of the $E^{(p)}(X,-)$ and the tower \eqref{eqn:HCTower0} in  \cite{LevineChowMov}. The method used is a variation on the classical Chow's moving lemma, aided by Gabber's approach to moving lemmas. The results of \cite{LevineChowMov} replace the total spectra $|E^{(p)}(X,-)|$ with functors
\[
E^{(p)}:\Sm/k^\op\to \Spt,
\]
and similarly for the tower \eqref{eqn:HCTower0} and layers $E^{(p/p+1)}$; these results require only that $E$ is homotopy invariant and satisfies Nisnevic excision. The homotopy invariance of the  $E^{(p)}(X,-)$ is also verified in  \cite{LevineChowMov}. The methods of  \cite{Loc} allows us to prove a localization property for the homotopy coniveau tower.

We iterate the functors $E^{(p)}$, and thereby find a simple description of the layers. In this, a crucial role is played by the 0th layers $(\Omega_{\P^1}^pE)^{(0/1)}$ of the $p$-fold $\P^1$-loop space $\Omega_{\P^1}^pE$, where
\[
(\Omega_{\P^1}E)(X):=\fib[E(X\times\P^1)\to E(X\times\infty)].
\]
Indeed, for $W$ smooth, the restriction map
\[
E^{(0/1)}(W)\to E^{(0/1)}(k(W))
\]
is a weak equivalence, which enables us to extend the functor $E^{(0/1)}$ to all $k$-schemes as a locally constant sheaf for the Zariski topology. One can then identify $E^{(p/p+1)}(X)$ with the simplicial spectrum $E_\loc^{(p/p+1)}(X,-)$ having $n$-simplices
\[
E_\loc^{(p/p+1)}(X,n)=\coprod_{x\in X^{(p)}(n)}(\Omega_{\P^1}^pE)^{(0/1)}(k(x)).
\]
Here $X^{(p)}(n)$ is the set of codimension $p$ points of $X\times\Delta^n$, with closure in good position.

Once one has this description of the layers, it is easy to compute the layers for $K$-theory, as one can easily show that $\Omega_{\P^1}^pK=K$ and that $K^{(0/1)}(F)$ is canonically a $K(K_0(F),0)=K(\Z,0)$ for $F$ a field. This gives a direct identification of $K_\loc^{(p/p+1)}(X,n)$ with Bloch's higher cycle group $z^p(X,n)$, and thus the weak equivalence $K^{(p/p+1)}(X,-)$ with $z^p(X,-)$.

The natural setting for the homotopy coniveau tower is in the Morel-Voevodsky homotopy category of $S^1$ spectra, $\SH_{S^1}(k)$. This is essentially the category of presheaves of spectra on $\Sm/k$, localized with respect to the Nisnevic topology and $\A^1$ weak equivalence. Voevodsky \cite{VoevSlice1} constructs the {\em slice tower}
\[
\ldots\to f_{n+1}E\to f_n E\to\ldots\to f_0E=E
\]
for $E\in \SH_{S^1}(k)$, where $f_nE\to E$ is universal for maps of the form $\Sigma_{\P^1}^d F\to E$, in analogy with the classical Postnikov tower for spectra. Our main result in this direction is the identification of this $S^1$-slice tower with our homotopy coniveau tower (for $k$ an infinite perfect field). The localization properties of the homotopy coniveau tower enable us to show how the truncation functors $f_n$ commute with the $\P^1$-loops functor: 
\[
f_{n+1}\circ\Omega_{\P^1}\cong \Omega_{\P^1}\circ f_n.
\]
This  in turn proves a connectedness conjecture of Voevodsky:
\begin{conj} Let $\Sigma^d_{\P^1}\SH_{S^1}(k)$ be the localizing subcategory of $\SH_{S^1}(k)$ generated by the presheaves  $\Sigma^d_{\P^1} F$, $F\in \SH_{S^1}(k)$. If $E$ is in $\Sigma^d_{\P^1}\SH_{S^1}(k)$, then $\Omega_{\P^1}\Sigma_{\P^1}E$ is also in $\Sigma^d_{\P^1}\SH_{S^1}(k)$.
\end{conj}

After this, we turn to the $\P^1$-stable theory. Let $\SH(k)$ be the Morel-Voevodsky stable homotopy category of $\P^1$-spectra over $k$. There is an infinite $\P^1$-suspension functor
\[
\Sigma_{\P^1}^\infty:\SH_{S^1}(k)\to \SH(k)
\]
with essential image denoted $\SH^\eff(k)$. Voevodsky defines the slice tower
\[
\ldots\to f_{p+1}\sE\to f_p\sE\to\ldots\to f_0\sE\to f_{-1}\sE\to\ldots\to\sE.
\]
in $\SH(k)$,
where the map $f_d\sE\to \sE$ is universal for maps $\sF\to\sE$ with $\sF$ in $\Sigma^d_{\P^1}\SH^\eff(k)$.

The localization property for the $E^{(p)}$ allows one to define a $\P^1$-spectrum $\fil_p\sE$ for a $\P^1$-spectrum $\sE:=(E_0, E_1,\ldots)$, by the formula
\[
(\fil_p\sE)_n:=  E_n^{(n+p)}.
\]
This gives the stable homotopy coniveau tower
\[
\ldots\to \fil_{p+1}\sE\to\fil_p\sE\to\ldots\to\fil_0\sE\to\fil_{-1}\sE\to\ldots\to\sE.
\]

 We show that $\fil_d\sE$ is in  $\Sigma^d_{\P^1}\SH^\eff(k)$, and that the canonical map
\[
 \fil_{d}\sE\to f_d\sE
\]
is an isomorphism, thus identifying the slice tower with the stable homotopy coniveau tower.

Finally, we compute of the 0th layer $\barfil_0$ in the homotopy coniveau tower for the motivic sphere spectrum $\1$, assuming that the base-field $k$ is perfect. The idea here is that the cycle-like description of $E_\loc^{(p/p+1)}(X,-)$ enables one to define a ``reverse cycle map"
\[
\rcyc:\HZ\to \barfil_0\1.
\]
It is then rather easy to show that $\rcyc$ induces a weak equivalence after applying the 0th layer functor $ \barfil_0$ again. However, since motivic cohomology is already the 0th layer of $K$-theory, applying $ \barfil_0$ leaves $\HZ$ unchanged, and similarly for $\barfil_0\1$, giving the desired weak equivalence $\HZ\sim \barfil_0\1$. The analogous statement for the slice filtration  in characteristic zero has been proven by Voevodsky  \cite{VoevSlice2} by a different method. 

In any case, for a $\P^1$-spectrum $\sE$,  each layer $\barfil_p\sE$ is a $\HZ$-module. Assuming that Voevodsky's slice tower has a lifting to a natural tower in Jardine's category of motivic symmetric spectra \cite{Jardine2}, work of Roendigs-Ostvar implies that the layers $s_n\sE$ in the slice tower for $\sE$ are naturally the Eilenberg-Maclane spectra of objects of the category of motives over $k$, $\DM(k)$. 

We thus have the objects $\pi^\mu_p\sE$ of $\DM(k)$ whose Eilenberg-Maclane $\P^1$-spectrum satisfies
\[
\sH((\pi^\mu_p\sE)[p])\cong \barfil_p\sE.
\]
For $E:\Sm/k\to\Spt$ the $0$th $S^1$-spectrum of $\sE$ and $X\in\Sm/k$, the spectral sequence associated to the homotopy coniveau tower can be expressed as
\[
E_2^{p,q}:=\H^p(X,\pi^\mu_{-q}\sE)\Longrightarrow \hat{E}_{-p-q}(X),
\]
wheere $\hat{E}$ is the completion of $E$ with respect to the homotopy coniveau tower. Under certain connectivity properties of $E$, one has $E=\hat{E}$.

Using the bi-graded homotopy groups of a $\P^1$-spectrum, this gives the weight-shifted spectral sequence
\[
E_2^{p,q}=\H^p(X,(\pi^\mu_{-q}\sE)\otimes\Z(b))\Longrightarrow
\hat{\sE}^{p+q,b}(X).
\]

For the $K$-theory $\P^1$-spectrum
\[
\sK:=(K,K,\ldots),
\]
our computation of the layers $K^{(p/p+1)}$ gives
\[
\pi^\mu_p\sK=\Z(p)[p];\quad p=0,\pm1,\pm2,\ldots,
\]
and we recover the Bloch-Lichtenbaum, Friedlander-Suslin spectral sequence
\[
E_2^{p,q}:=\H^{p-q}(X,\Z(-q))\Longrightarrow K_{-p-q}(X).
\]

We have had a great deal of help in developing the techniques that went in to this paper. Fabien Morel played a crucial role in numerous discussions on the $\A^1$-stable homotopy category and related topics. Conversations with Bruno Kahn, Jens Hornbostel and Marco Schlichting were very helpful, as were the lectures of Bjorn Dundas and Vladamir Voevodsky at the Sophus Lie Summer Workshop in $\A^1$-stable homotopy theory. I would also like to thank Paul Arne Ostvar for his detailed comments on an earlier version of this manuscript. Finally, I am very grateful to the Humboldt Foundation and the Universitit\"at Essen for their support of this reseach, especially my colleagues H\'el\`ene Esnault and Eckart Viehweg.
  
\hbox{}
\newpage

\section{Spaces, spectra and homotopy categories}

\subsection{Presheaves of spaces} Let $\Spc$ denote the category of {\em spaces}, i.e., simplicial sets, and $\Spc_*$ the category of pointed simplicial sets. For a category $\sC$, we have the category $\Spc(\sC)$ of  presheaves of spaces on $\sC$, and $\Spc_*(\sC)$ of  presheaves of pointed spaces. 

We give $\Spc$ and $\Spc_*$ the standard model structures: cofibrations are  monomorphisms, weak equivalences are weak equivalences on the geometric realization, and fibrations are detemined by the RLP with respect to trivial cofibrations; the fibrations are then exactly the Kan fibrations. We let  $|A|$ denote the geometric realization, and $[A,B]$  the homotopy classes of (pointed) maps $|A|\to |B|$.

We give  $\Spc(\sC)$  and $\Spc_*(\sC)$ the model structure of functor categories described by Bousfield-Kan \cite{BousfieldKan}. That is, the cofibrations and weak equivalences are the pointwise ones, and the fibrations are determined by the RLP with respect to trivial cofibrations. We let $\sH\Spc(\sC)$ and $\sH \Spc_*(\sC)$ denote the associated homotopy cateogies. 

Note that $\Spc(\sC)$ and $\Spc_*(\sC)$ inherit operations from $\Spc$ and $\Spc_*$, for instance limits and colimits. In particular, in $\Spc_*(\sC)$ we have wedge product $A\wedge B:=A\times B/A\times*\cup
*\times B$. We also have the inclusions $\Spc\to \Spc(\sC)$, $\Spc_*\to \Spc(\sC)$ as constant presheaves, giving us the suspension functor on $\Spc_*(\sC)$, $\Sigma A:=S^1\wedge A$, and the inclusion $+:\Spc(\sC)\to \Spc(\sC)$ by adding a disjoint base-point. These operations pass to the homotopy category.

\subsection{Spectra} Let $\Spt$ denote the category of
spectra. To fix ideas, a spectrum will be a sequence of pointed spaces
$E_0, E_1,\ldots$ together with maps of pointed simplicial sets $\epsilon_n:S^1\wedge E_n\to E_{n+1}$.  Maps of spectra are maps of the underlying spaces which are compatible with the attaching maps $\epsilon_n$. The stable homotopy groups $\pi_n^s(E)$ are defined by
\[
\pi_n^s(E):=\lim_{m\to\infty} [S^{m+n},E_m].
\]

The category $\Spt$ has the following model structure: Cofibrations are maps $f:E\to F$ such that $E_0\to F_0$ is a cofibration, and for each $n\ge0$, the map
\[
E_{n+1}\amalg_{S^1\wedge E_n}S^1\wedge F_n\to F_{n+1}
\]
is a cofibration.  Weak equivalences are the stable weak equivalences, i.e., maps $f:E\to F$ which induce an isomorphism on $\pi_n^s$ for all $n$. Fibrations are characterized by having the RLP with respect to trivial cofibrations.

Let $\Spt(\sC)$ be the category of presheaves of spectra on $\sC$. We use the following model structure on $\Spt(\sC)$ (see \cite{Jardine1}): Cofibrations and weak equivalences are given pointwise, and fibrations are characterized by having the RLP with respect to trivial cofibrations. We denote   the associated homotopy category by $\sH\Spt(\sC)$.   We write $\SH$ for the homotopy category of $\Spt$.

\subsection{Notation}
For a scheme $B$, $\Sm/B$ is the category of smooth separated $B$-schemes of finite type.
For a morphism $f:Y\to X$ in $\Sm/B$, and $E\in\Spt(B)$, we write $E(X/Y)$ for the homotopy fiber of $f^*:E(X)\to E(Y)$. Similarly, if $(Z,z:B\to Z)$ is a pointed $B$-scheme, $E(Z\wedge(X/Y))$ is the homotopy fiber of $z^*:E(Z\times X/Z\times Y)\to E(X/Y)$. In case $j:U\to X$ is an open immersion with closed complement $W\subset X$, we write $E^W(X)$ for $E(X/U)$. We write $X_+$ for $X\amalg\Spec k/\Spec k$.

Given an $E\in\Spt(B)$ and a $Y\to B\in \Sm/B$, we let $E^{(Y)}$ denote the presheaf $E^{(Y)}(Z):=E(Z\times_BY)$.

\subsection{Nisnevic model structure} Fix a noetherian base-scheme $B$  and  let $\sC$ be a subcategory of $\Sm/B$ with the same objects as $\Sm/B$ and containing all the smooth $B$-morphisms. In particular, the Nisnevic topology is defined on $\sC$.  

For a point $x\in X$, with $X\in\sC$, and $E=(E_0, E_1,\ldots)$ a presheaf of spectra on $\sC$, the  stalk of $E$ at $x$, $E_x$, is the spectrum $(E_{0x}, E_{1x},\ldots)$, where $E_{nx}$ is the stalk (in the Nisnevic topology) of the presheaf of spaces $E_n$ at $x$.

Let $\Spt_\Nis(\sC)$ denote the model category with the same underlying category and cofibrations as $\Spt(\sC)$, where  a map $f:E\to F$ is a weak equivalence if  $\tilde{f}_x:\tilde{E}_x\to \tilde{F}_x$ is a weak equivalence in $\Spt$ for all $x\in X\in\sC$, and the fibrations are characterized by having the right lifting property (RLP) with respect to trival cofibrations. We let $\sH\Spt_\Nis(\sC)$ denote the associated homotopy category and write $\Spt_\Nis(B)=\Spt_\Nis(\Sm/B)$,  $\SH_s(B):=\sH\Spt_\Nis(\Sm/B)$. For details, we refer the reader to \cite{Jardine1}.

\subsection{$\A^1$-local model structure} One imposes the relation of $\A^1$-weak equivalence in 
$\SH_s(B)$ by means of Bousfield localization applied to the model category $\Spt_\Nis(B)$. 

An object $E\in \Spt(B)$ is called {\em $\A^1$-local}   if the map $E\to E^{(\A^1)}$ induced by the projections $Y\times\A^1\to Y$ is a weak equivalence in  $\Spt_\Nis(B)$. A map $f:F\to F'$ in $\Spt_\Nis(B)$ is an {\em $\A^1$-weak equivalance} if 
\[
f^*:\Hom_{\SH_s(B)}(F',E)\to \Hom_{\SH_s(X)}(F,E)
\]
is an isomorphism for all $\A^1$-local $E$. $\Spt_{S^1}(B)$ is the model category with the same underlying category and cofibrations as $\Spt_\Nis(B)$, the weak equivalences being the $\A^1$-weak equivalences, and the fibrations determined by the RLP with respect to trivial cofibrations. The fact that this is indeed a model category is discussed in \cite{Jardine2}.

We  refer to $\Spt_{S^1}(B)$ as the category of {\em $S^1$-spectra over $B$}. The homotopy category $\sH\Spt_{S^1}(B)$ is denoted $\SH_{S^1}(B)$.

\subsection{Simplicial spectra}

For a spectrum $E$, we have the Postnikov tower
\[
\xymatrix{
\ldots \ar[r]& \tau_{\ge N}E\ar[rr]\ar[dr] &&\ar[dl] \tau_{\ge N-1}E\ar[r]& \ldots\\
&&E
}
\]
with $\tau_{\ge N}E\to E$ the $N-1$-connected cover of $E$, i.e., $\tau_{\ge N}E\to E$ is an isomorphism on homotopy groups $\pi_n$ for $n\ge N$, and $\pi_n(\tau_{\ge N}E)=0$ for $n<N$. One can make this tower functorial in $E$, so we can apply the construction $\tau_{\ge N}$ to functors $E:\sC\to\Spt$.

We have the category $\Ord$ with objects the finite ordered sets $[n]:=\{0<\ldots<n\}$, $n=0, 1, \ldots$, and maps order-preserving maps of sets. Let $\Ord_{\le N}$ be the full subcategory with objects $[n]$, $0\le n\le N$.

Let $E:\Ord^\op\to \Spt$ be a simplicial spectrum. We have the $N$-truncated simplicial spectrum $E_{\le N}:\Ord^\op_{\le N}\to \Spt$,  the associated total spectrum $|E_{\le N}|$, and the tower  of spectra
\begin{equation}\label{eqn:SimplTower}
|E_{\le 0}|\to \ldots\to |E_{\le N}|\to\ldots\to |E|
\end{equation}
Since taking the total spectrum commutes with filtered colimits, we have the natural weak equivalences
\[
\xymatrix{
\hocolim_M|E_{\le M}| \ar[r]^-\sim& |E|\\
 \hocolim_{N, M}\tau_{\ge -N}|E_{\le M}|\ar[r]_-\sim\ar[u]^\sim& \hocolim_N\tau_{\ge -N}|E|\ar[u]_\sim .}
\]

When the context makes the meaning clear, we will often omit the separate notation for the total spectrum, and freely pass between a simplicial spectrum and its associated total spectrum.

 \section{The homotopy coniveau tower}

\subsection{The construction} We fix a  noetherian base scheme $S$, separated and   of finite Krull dimension.  

We have the cosimplicial scheme $\Delta^*$, with
\[
\Delta^r=\Spec(\Z[t_0,\ldots, t_r]/\sum_jt_j-1).
\]
The {\em vertices} of $\Delta^r$ are the closed subschemes $v^r_i$ defined by
$t_i=1$, $t_j=0$ for $j\neq i$. A {\em face} of $\Delta^r$ is a closed subscheme
defined by equations of the form $t_{i_1}=\ldots = t_{i_s}=0$.

Let $E\in\Spt(S)$ be a presheaf of spectra. For  $X$ in $\Sm/S$ with closed subscheme $W$ and open complement $j:X\setminus W\to X$, we have the 
 homotopy fiber $E^W(X)$ of $j^*:E(X)\to E(X\setminus W)$. If we have a chain of closed subsets
$W'\subset W\subset X$, we have a natural map $i_{W',W*}:E^{W'}(X)\to E^W(X)$ and a natural weak equivalence
\begin{equation}\label{Excision1}
\cofib(i_{W',W*}:E^{W'}(X)\to E^W(X))\sim
E^{W\setminus W'}(X\setminus W').
\end{equation}
Here ``cofib" means homotopy  cofiber in the category of spectra.

For $X$ in $\Sm/S$, we let $\sS_X^{(p)}(r)$ denote the set of closed subsets
$W$ of $X\times\Delta^r$ such that 
\[
\codim_{X\times F}(W\cap(X\times F))\ge p
\]
for all faces $F$ of $\Delta^r$. Clearly, sending $r$ to $\sS_X^{(p)}(r)$ defines a
simplicial set $\sS_X^{(p)}(-)$. We let $X^{(p)}(r)$ be the set of codimension $p$ points $x$ of $X\times\Delta^r$ with closure $\bar{x}\in\sS_X^{(p)}(r)$.

We let $E^{(p)}(X,r)$ denote the (filtered) limit
\[
E^{(p)}(X,r)=\hocolim_{W\in \sS_X^{(p)}(r)}E^W(X\times\Delta^r).
\]
Sending $r$ to $E^{(p)}(X,r)$ defines a simplicial spectrum $E^{(p)}(X,-)$. Since
$\sS_X^{(p+1)}(r)$ is a subset of $\sS_X^{(p)}(r)$, we have the tower of simplicial spectra
\begin{equation}\label{HConiveauTower}
\ldots\to E^{(p+1)}(X,-)\to E^{(p)}(X,-)\to\ldots\to E^{(0)}(X,-),
\end{equation}
which we call the {\em homotopy coniveau tower}.  We let $E^{(p/p+1)}(X,-)$ denote
the cofiber of the map $E^{(p+1)}(X,-)\to E^{(p)}(X,-)$.

Two properties of $E$ that we shall often require are:
\begin{enumerate}
\item[A1.] $E$ is {\em homotopy invariant}: For each $X$ in $\Sm/S$, the map
$p^*:E(X)\to E(X\times\A^1)$ is a weak equivalence.
\item[A2.] $E$ satisfies {\em Nisnevic excision}: Let $f:X'\to X$ be an \'etale morphism in $\Sm/S$, and $W\subset X$ a closed subset. Let $W'=f^{-1}(W)$, and suppose that $f$ restricts to an isomorphism $W'\to W$. Then $f^*:E^W(X)\to E^{W'}(X')$ is a weak equivalence.
\end{enumerate} 
Instead of A2, we will occasionally require the weaker condition:
\begin{enumerate}
\item[A2'.] $E$ satisfies {\em Zariski excision}: Let $j:U'\to X$ be an open immersion in $\Sm/S$, and $W\subset X$ a closed subset contained in $U$. Then $j^*:E^W(X)\to E^W(U)$ is a weak equivalence.
\end{enumerate}

We introduce one final axiom to handle the case of finite residue fields. Suppose we have a finite Galois extension  $k\to k'$ with group $G$. Given $E\in\Spt(k)$, define $\pi_*^G\pi^*E$ by
\[
\pi_*^G\pi^*E(X):= E(X_{k'})^G,
\]
where $(-)^G$ denotes a functorial model for the $G$ homotopy fixed point spectrum. We have as well the natural transformation
\[
\pi^*:E\to \pi_*^G\pi^*E.
\] 
\begin{enumerate}
\item[A3.]   Suppose that $k$  is a finite field. Let  $k\to k'$ be a finite Galois extension of $k$ with group $G$. Then after inverting $|G|$, the natural transformation $\pi^*:E\to \pi^G_*\pi^* E$ is a weak equivalence.
\end{enumerate}
\begin{rem}\label{rem:A3}It is shown in \cite[Corollary 9.4.2]{LevineChowMov} that $E$ satisfies A3 if $E$ is the 0-spectrum of some $\P^1$-$\Omega$-spectrum $\sE\in\Spt_{\P^1}(k)$  (see \S~\ref{sec:P1Spectra}).
\end{rem}

\begin{Def} Let $X$ be in $\Sm/S$. The {\em weight-completed} spectrum $\hat{E}(X)$ is
\[
\hat{E}(X)=\holim_pE^{(0/p)}(X,-).
\]
\end{Def}

\begin{prop}\label{SpecSeqProp2} Take $E\in \Spt(S)$, where $S$ is a noetherian scheme of finite Krull dimension.
\begin{enumerate}
\item  There is a weakly convergent spectral sequence
\[
E_1^{p,q}=\pi_{-p-q}(E^{(p/p+1)}(X,-))\Longrightarrow \hat{E}_{-p-q}(X).
\] 
\item If $E=\tau_{\ge N}E$ for some $N$, then the above spectral sequence is strongly convergent and the canonical map $E^{(0)}(X,-)\to \hat{E}(X)$ is a weak equivalence.
\item If $E$ is homotopy invariant, the canonical map $E(X)\to E^{(0)}(X,-)$ is a weak equivalence.
\end{enumerate}
\end{prop}

\begin{proof}  
(1) The spectral sequence is  
constructed by the standard process of linking the
long exact sequences of homotopy groups arising from the homotopy cofiber sequences
\[
E^{(p+1)}(X,-)\to E^{(p)}(X,-)\to E^{(p/p+1)}(X,-). 
\]
The first assertion then follows from the general theory of homotopy limits (see \cite{BousfieldKan}).

For (2), suppose $E=\tau_{\ge N}E$. We first show that the sequence is strongly convergent.

By \eqref{Excision1} and a limit argument, we have
\[
\pi_m(E^{(p/p+1)}(X,r))=\lim_{\substack{\to\\W'\subset W}}\pi_m\big( E^{W\setminus
W'}(X\times\Delta^r\setminus W')\big),
\]
where the limit is over $W'\in \sS_X^{(p+1)}(r)$, $W\in \sS_X^{(p)}(r)$.  It follows that 
$\pi_m(E^{(p/p+1)}(X,r))=0$ for $m<N-1$. 

From the tower \eqref{eqn:SimplTower}, we thus have the the strongly
convergent spectral sequence
\[
E^{a,b}_1=\pi_{-a}(E^{(p/p+1)}(X,-b))\Longrightarrow \pi_{-a-b}(E^{(p/p+1)}(X,-)).
\]
Since $ \sS_X^{(p)}(r)=\0$ for $p>\dim X+r$, this implies that
\[
\pi_{-p-q}E^{(p/p+1)}(X,-)=0\]
 for $p>-p-q+\dim X+N+1$, from which it follows that the
spectral sequence \eqref{ConSpSeq} is strongly convergent.  Similarly, it follows that the natural map $E^{(0)}(X,-)\to \hat{E}(X)$ is a weak equivalence. 

For (3) the simplicial spectrum $E^{(0)}(X,-)$ is just the simplicial
spectrum $E(X\times\Delta^*)$, i.e., $r\mapsto E(X\times\Delta^r)$. Since $E$ is homotopy invariant, the natural map
\[
E(X)\to E(X\times\Delta^*)
\]
is a weak equivalence, completing the proof.
\end{proof}

\begin{cor} Take $E\in \Spt(S)$, where $S$ is a noetherian scheme of finite Krull dimension.
If  $E$ is homotopy invariant and  $E=\tau_{\ge N}E$ for some $N$, then the homotopy coniveau tower \eqref{HConiveauTower} yields a strongly convergent spectral sequence
\begin{equation}\label{ConSpSeq}
E_1^{p,q}=\pi_{-p-q}(E^{(p/p+1)}(X,-))\Longrightarrow E_{-p-q}(X).
\end{equation}
which we call the homotopy coniveau spectral sequence.
\end{cor}

\subsection{First properties}\label{subsec:FirstProperties} We give a list of elementary properties of the spectra $E^{(p)}(X,-)$
\begin{enumerate}
\item Sending $X$ to $E^{(p)}(X,-)$ is functorial for equi-dimensional (e.g. flat) maps $Y\to X$ in $\Sm/S$.
\item The pull-back $p_1^*:E^{(p)}(X,-)\to E^{(p)}(X\times\A^1,-)$ is a weak equivalence. The proof is the same as that for Bloch's cycle complexes, given in \cite{AlgCyc}. For details, see \cite[Theorem 3.3.5]{LevineChowMov}.
\item Sending $E$ to $E^{(p)}(X,-)$ is functorial in $E$.
\item The functor $E\mapsto E^{(p)}(X,-)$ sends weak equivalences to weak equivalences, 
and send homotopy (co)fiber sequences to homotopy (co)fiber sequences.
\end{enumerate}

Exactly the same properties hold for the layers $E^{(p/p+r)}$.

\section{Localization} We now show that the simplicial spectra $E^{(p)}(X,-)$ behave well with respect to localization.

\subsection{Stable homology of spectra} For a simplicial set $S$, we have the
simplicial abelian group $\Z S$, with $n$-simplices $\Z S_n$ the free abelian group on
$S_n$. If $S,*$ is a pointed simplicial set, define $\Z(S,*)_n:=\Z S_n/\Z*$.

Let $E=\{E_n,
\phi_n:\Sigma E_n\to E_{n+1}\}$ be a spectrum; we take the $E_n$ to be pointed simplicial
sets, and the $\phi_n$ to be maps of pointed simplicial sets. Form the spectrum $\Z E$ by
taking $(\Z E)_n=\Z (E_n,*)$, where $\Z\phi_n:\Sigma(\Z E)_n\to(\Z E)_{n+1}$ is the map
induced by $\phi_n$, composed with the natural map $\Sigma(\Z E)_n\to\Z(\Sigma E_n)$.
The natural maps $E_n\to\Z E_n$ give a natural map $E\to \Z E$ of spectra; one shows
that this construction respects weak equivalence and taking homotopy cofibers. The {\em stable homology}
$H_n(E)$ is defined by $H_n(E)=\pi_n(\Z E)$. Using the Dold-Thom theorem, one has the
formula for $H_n(E)$ as
\[
H_n(E)=\lim_\to \tilde{H}_{n+m}(E_m),
\]
where $\tilde{H}$ is  reduced homology and the maps in the limit are  the composition
\[
\tilde{H}_{n+m}(E_n)\cong \tilde{H}_{n+m+1}(\Sigma E_n)\xrightarrow{\phi_{n*}} 
\tilde{H}_{n+m+1}(E_{n+1}).
\]

The Hurewicz theorem for simplicial sets gives the following analogous result for
spectra:

\begin{prop}\label{Hurewicz} Let $E$ be a spectrum which is $N$-connected for some
$N\in\Z$. Then $\pi_n(E)=0$ for all $n$ if and only if $H_n(E)=0$ for all $n$.
\end{prop}

\begin{proof} Since both $\pi_n$ and $H_n$ respect weak equivalence, and are
compatible with suspension of spectra, we may assume that $N\ge 1$, and that $E$ is an
$\Omega$-spectrum, i.e., the natural maps $E_n\to\Omega E_{n+1}$ are weak
equivalences. Then $\pi_n(E)=\pi_{n+m}(E_m)$ for all $m$. Suppose $H_n(E)=0$ for all
$n$; we prove by induction that $\pi_{n+m}(E_m)=0$ for all $n$ and $m$. 

By assumption $\pi_{N+m}(E_m)=0$ for all $m$, with $N\ge 1$. We may therefore proceed
by induction on $n$ to show that $\pi_{n+m}(E_m)=0$ for all $n$ and $m$. Supposing
that $\pi_{n+m-1}(E_m)=0$ for all $m$, the Hurewizc theorem implies that the Hurewicz
map  $\pi_{n+m}(E_m)\to \tilde{H}_{n+m}(E_m)$ is an isomorphism for all $m$, and one easily
checks that the Hurewicz map is compatible with the limits defining $H_n$ and $\pi_n$.
Thus, the maps $\tilde{H}_{n+m}(E_m)\to \tilde{H}_{n+m+1}(E_{m+1})$ are isomorphisms for all $m$;
since the limit is zero by assumption, this implies that $\tilde{H}_{n+m}(E_m)=0$ for all $m$,
whence $\pi_{n+m}(E_m)=0$ for all $m$.

The proof that $\pi_n(E)=0$ for all $n$ implies $H_n(E)=0$ for all $n$ is similar, but
easier, and is left to the reader.
\end{proof}

\subsection{The localization theorem}  

Let $X$ be smooth and essentially of finite type over $S$, and let $j:U\to X$ be an
open subscheme, with complement $i:Z\to X$. We let $\sS_{X,Z}^{(p)}(r)$ denote the subset of
$\sS_X^{(p)}(r)$ consisting of those $W$ contained in $Z\times\Delta^r$. Let $\sS_{U/X}^{(p)}(r)$
be the image of $\sS_X^{(p)}(r)$ in $\sS_{U}^{(p)}(r)$ under $(j\times\id)^{-1}$. Taking the colimit
of $E^W(X\times\Delta^r)$ over $W\in \sS_{X,Z}^{(p)}(r)$ and varying $r$ and $p$ gives us the
tower of  simplicial spectra
\[
\ldots\to E^{(p+1)}_Z(X,-)\to E^{(p)}_Z(X,-)\to \ldots\to
E^{(d)}_Z(X,-)=E^{(0)}_Z(X,-),
\]
where $d$ is any integer satisfying $d\le\codim_XZ_j$ for all irreducible components
$Z_j$ of $Z$. Similarly, taking the colimit of the $E^W(U\times\Delta^r)$ over $W\in \sS_{U/X}^{(p)}(r)$ for varying $p$ and
$r$ gives the 
tower of  simplicial spectra
\[
\ldots\to E^{(p+1)}(U_X,-)\to E^{(p)}(U_X,-)\to \ldots\to
E^{(0)}(U_X,-).
\]
We have as well the natural maps 
\begin{gather*}
i_*:E^{(p)}_Z(X,r)\to E^{(p)}(X,r),\
j^{*!}:E^{(p)}(X,r)\to E^{(p)}(U_X,r),\\
\iota:E^{(p)}(U_X,r)\to E^{(p)}(U,r),\ j^*:E^{(p)}(X,r)\to E^{(p)}(U,r), 
\end{gather*}
with $j^*=\iota\circ j^{*!}$. 
  
Let $E^{(p/p+s)}(-)$ denote the cofiber of the maps $E^{(p+s)}(-)\to
E^{(p)}(-)$. Supposing that $E$ satisfies Zariski excision, we have the homotopy fiber sequences
\begin{align*}
&E^{(p)}_Z(X,r)\xrightarrow{i_*} E^{(p)}(X,r)\xrightarrow{j^{*!}}
E^{(p)}(U_X,r)\\
&E^{(p/p+s)}_Z(X,r)\xrightarrow{i_*} E^{(p/p+s)}(X,r)\xrightarrow{j^{*!}}
E^{(p/p+s)}(U_X,r).
\end{align*}
These give the homotopy fiber sequences of simplicial spectra 
\begin{align*}
&E^{(p)}_Z(X,-)\xrightarrow{i_*} E^{(p)}(X,-)\xrightarrow{j^{*!}}
E^{(p)}(U_X,-)\\
&E^{(p/p+s)}_Z(X,-)\xrightarrow{i_*} E^{(p/p+s)}(X,-)\xrightarrow{j^{*!}}
E^{(p/p+s)}(U_X,-)\\
\end{align*}

The localization theorem is

\begin{thm}\label{LocThm} Let $E$ be in $\Spt(S)$. Suppose the base-scheme $S$ is a scheme essentially of
finite type over a semi-local DVR with infinite residue fields. Then the maps
\begin{align*}
&E^{(p)}(U_X,-)\to E^{(p)}(U,-)\\
&E^{(p/p+s)}(U_X,-)\to E^{(p/p+s)}(U,-)\\
\end{align*}
are weak equivalences
\end{thm}

\begin{proof} The second weak equivalence follows from the first by taking cofibers.

For the first map,  this result follows by exactly the same method as used in the proof of
\cite[Theorem 8.10]{Loc}. Indeed,  to show that the map
$E^{(p)}(U_X,-)\to E^{(p)}(U,-)$ is a weak equivalence, it suffices to prove the result with $E^{(p)}(-, n)$ replaced by $\tau_{\ge m}E^{(p)}(-,n)$ for all $m$. By the Hurewicz theorem
(Proposition~\ref{Hurewicz}), it suffices to show that $E^{(p)}(U_X,-)\to E^{(p)}(U,-)$
is a homology isomorphism. This follows by applying
\cite[Theorem 8.2]{Loc}, just as in the proof of Theorem 8.10 ({\it loc. cit.}). For the reader's convenience, we include a sketch of the argument.

Let $E=(E_0, E_1,\ldots)$ be a spectrum. Using the Dold-Kan correspondence, we can identify the stable homology spectrum $\Z E$ with the complex formed by taking the normalized complex of the simplicial abelian group $\Z(E_n,*)$ and then taking the limit over the bonding maps $\Z(E_n,*)[n]\to
\Z(E_{n+1},*)[n+1]$. Abusing notation, we denote this complex also by $\Z E$; for the remainder of the proof, $\Z E$ will mean the complex, not the spectrum.

For $W\subset U\times \Delta^r$, we have the complex
\[
\Z(\tau_{\ge m}E^W(U\times\Delta^r))
\]
computing the stable homology of $\tau_{\ge m}E^W(U\times\Delta^r)$. Taking the limit of $W\in
 \sS_U^{(p)}(r)$ or in $\sS_{U/X}^{(p)}(r)$ gives us  the complexes $\Z E^{(p)}_m(U,r)$ and 
$\Z E^{(p)}_m(U_X,r)$, which  compute the stable homology of $\tau_{\ge m}E^{(p)}(U,r)$ and $\tau_{\ge m}E^{(p)}(U_X,r)$.

For $W\subset \amalg_{r=0}^N U\times\Delta^r$, let $W_n\subset U\times\Delta^n$ be the union of $(\id\times g)^{-1}(W)$, as $g:\Delta^n\to\Delta^r$ runs over structure morphisms for the cosimplicial scheme $\Delta$. Using the usual alternating sum of the  pullback by coface maps $\id\times\delta^r_i:U\times\Delta^r\to U\times\Delta^{r+1}$, we form the double complex $n\mapsto
\Z \tau_{\ge m}E^{W_n}(U\times\Delta^n)$ and denote the associated total complex by $\Z E^W_m(U\times\Delta^*)$. Thus the limit of the complexes $\Z E^W_m(U\times\Delta^*)$ over $W\in 
 \sS_U^{(p)}(r)$ or in $\sS_{U/X}^{(p)}(r)$, $r=1, 2,\ldots$,  computes the stable  homology of the simplicial spectra $n\mapsto \tau_{\ge m}E^{(p)}(U,n)$ and $n\mapsto \tau_{\ge m}E^{(p)}(U_X,n)$. We denote the  limits of these complexes by 
$\Z E_m^{(p)}(U)^*$ and $\Z E_m^{(p)}(U_X)^*$, respectively.  It thus suffices to show that
\[
\iota_\Z: \Z E^{(p)}_m(U_X)^*\to  \Z E^{(p)}_m(U)^*
\]
is a quasi-isomorphism for all $m$.

For 
$W\in  \sS_U^{(p)}(r)$, $W'\in \sS_{U/X}^{(p)}(r)$, let
\[
\iota_W:\Z E^W_m(U\times\Delta^*)\to\Z E^{(p)}_m(U)^*
\]
 and 
 \[
\iota_{W'}^X:\Z E^{W'}_m(U\times\Delta^*)\to\Z E^{(p)}_m(U_X)^*
\]
be the canonical maps.

Next, we construct another pair of complexes which approximate $\Z E^{(p)}_m(U)^*$ and $\Z E^{(p)}_m(U_X)^*$. For this, fix an integer $N\ge0$. Let $\partial\Delta^N_i\subset \Delta^N$ be the subscheme defined by $t_i=0$; for $I\subset \{0,\ldots, N\}$ let $\partial\Delta^N_I$ be the face $\cap_{i\in I}\partial\Delta^N_i$. For $I\supset J$, let $i_{J,I}:\Delta^N_I\to\Delta^N_J$ be the inclusion.

Let $\Z\Sm/S$ be the additive category generated by $\Sm/S$, i.e., for connected $X$, $Y$, $\Hom_{\Z\Sm/S}(X,Y)$ is the free abelian group on the set of morphisms  $\Hom_{\Sm/S}(X,Y)$, and disjoint union becomes direct sum. We will construct objects in the category of complexes $C(\Z\Sm/S)$.

Form the complex $(\Delta^N,\partial\Delta^N)^*$ which is $\oplus_{I,\ |I|=n}\partial\Delta^N_I$ in degree $-n$, and with differential 
\[
d^{-n}:(\Delta^N,\partial\Delta^N)^{-n}\to (\Delta^N,\partial\Delta^N)^{-n+1}
\]
given by $d^{-n}:=\prod_{I,\ |I|=n}d^{-n}_I$, where
\[
d^{-n}_I:\partial\Delta^N_I\to \oplus_{J,\ |J|=n-1}\partial\Delta^N_J
\]
is the sum
\[
d^{-n}_I:=\sum_{j=1}^n i_{I\setminus\{i_j\},I},
\]
where $I=(i_1,\ldots, i_n)$, $i_1<\ldots<i_n$.

We also have the complex $\Z\Delta^*$, which is $\Delta^n$ in degree $n$, with differential the usual alternating sum of coboundary maps.

The identity map on $\Delta^N$ extends to a map of complexes
\[
\Phi^N:\Z\Delta^*\to (\Delta^N,\partial\Delta^N)[-N];
\]
the maps in degree $r<N$ are all $\pm\id_{\Delta^r}$. We can take the product of   this construction with $U$, giving us the complex 
$U\times(\Delta^N,\partial\Delta^N)$ and the map  of complexes
\[
 \Phi^N:U\times\Z\Delta^*\to U\times(\Delta^N,\partial\Delta^N)[-N]
\]

 For $W\in  \sS_U^{(p)}(N)$, form the complex $\Z  E^W_m(U\times(\Delta^N,\partial\Delta^N))$ by taking $\oplus_{I, |I|=n} \Z\tau_{\ge m}E^{W_{N-n}}(U\times\Delta^{N-n})$ in degree $-n$, using the differentials in $U\times(\Delta^N,\partial\Delta^N)$ to form a double complex and then taking the total complex. We thus have the map of complexes
\[
\Phi^{N*}_W: \Z E^W_m(U\times(\Delta^N,\partial\Delta^N))[-N]\to 
\Z E^W(U\times\Delta^*).
\]
One shows that $\Phi^{N*}_W$ induces a homology isomorphism in degrees $<N$.

Take $W\in \sS_U^{(p)}(N)$. The main construction of \cite{Loc} gives a map of complexes
\[
\Psi_W:U\times\Z\Delta^*\to U\times(\Delta^N,\partial\Delta^N)[-N]
\]
and a degree -1 map
\[
H_W:U\times\Z\Delta^*\to U\times(\Delta^N,\partial\Delta^N)[-N]
\]
with the following properties:
\begin{enumerate}
\item $dH_W=\Psi_W-\Phi^N$.
\item Write $\Psi_W$ as a sum 
\[
\Phi_W=\sum_{i=0}^N\sum_{I,j\ |I|=i} n_{ij}\psi_{ijI}
\]
with $\psi_{ijI}:\Delta^{N-i}\to \partial\Delta^N_I=\Delta^{N-i}$ maps in $\Sm/S$. Then $\psi_{ijI}^{-1}(W_{N-i})$ is in $\sS_{U/X}^{(p)}(N-i)$.
\item Write $H_W$ as a sum 
\[
H_W=\sum_{i=0}^N\sum_{I,j\ |I|=i} n_{ij}H_{ijI}
\]
with $H_{ijI}:\Delta^{N-i+1}\to \partial\Delta^N_I=\Delta^{N-i}$ maps in $\Sm/S$. Then $H_{ijI}^{-1}(W_{N-i})$ is in $\sS_{U}^{(p)}(N-i+1)$. If $W'\subset W_{N-i}$ is in $\sS_{U/X}^{(p)}(N-i)$, then 
$H_{ijI}^{-1}(W')$ is in $\sS_{U/X}^{(p)}(N-i+1)$.
\end{enumerate}

Thus $\Psi_W$ induces the map of complexes
\[
\Psi_W^*: \Z E^W_m(U\times(\Delta^N,\partial\Delta^N))[-N]\to
 \Z E^{(p)}_m(G,U_X)^*
 \]
 and $H_W$ gives a degree 1 map
 \[
 H_W^*: \Z E^W_m(U\times(\Delta^N,\partial\Delta^N))[-N]\to
 \Z E^{(p)}_m(U)^*
 \]
 with 
 \[
 dH_W^*=\iota_\Z\circ\Psi_W^*-\iota_W\circ \Phi^{N*}_W
 \]
 Furthermore, if $W'\subset W$ is in $\sS_{U/X}^{(p)}(N)$, then $H_W$ gives a
  degree 1 map
 \[
H_W^{X*}: \Z E^{W'}_m(U\times(\Delta^N,\partial\Delta^N))[-N]\to
 \Z E^{(p)}_m(U_X)^*
 \]
 with 
 \[
 dH_W^{X*} = \Psi_W^*-\iota_{W'}^X\circ \Phi^{N*}_{W'}.
 \]
 Since $\Phi^{N*}_W$ is a homology isomorphism in degrees $<N$ and $\Z E^{(p)}_m(U)^*$ and
 $\Z E^{(p)}_m(U_X)^*$ are the limits of $\Z E^W_m(U\times\Delta^*)$ and $\Z E^{W'}_m(U\times\Delta^*)$, respectively, this shows that $\iota_\Z$ is a quasi-isomorphism, completing the proof.
\end{proof}

\begin{cor}\label{cor:LocThm} Let $E$ be in $\Spt(k)$ satisfying Zariski excision; if $k$ is finite, we assume in addition that $E$ satisfies axiom A3. Let $j:U\to X$ be an open immersion in $\Sm/k$ with complement $i:Z\to X$. Then the sequences of  spectra 
\begin{align*}
&E^{(p)}_Z(X,-)\xrightarrow{i_*} E^{(p)}(X,-)\xrightarrow{j^*}
E^{(p)}(U,-)\\
&E^{(p/p+s)}_Z(X,-)\xrightarrow{i_*} E^{(p/p+s)}(X,-)\xrightarrow{j^*}
E^{(p/p+s)}(U,-)
\end{align*}
extend canonically to distinguished triangles in $\SH$. 
\end{cor}

\begin{proof} If $k$ is infinite, this follows directly from the weak homotopy fiber sequences 
\begin{align*}
&E^{(p)}_Z(X,-)\xrightarrow{i_*} E^{(p)}(X,-)\xrightarrow{j^{*!}}
E^{(p)}(U_X,-)\\
&E^{(p/p+s)}_Z(X,-)\xrightarrow{i_*} E^{(p/p+s)}(X,-)\xrightarrow{j^{*!}}
E^{(p/p+s)}(U_X,-)
\end{align*}
and Theorem~\ref{LocThm}. For $k$ finite, one uses A3 and the existence of infinite extensions of $k$ of relatively prime power degree to reduce to the case of an infinite field.
\end{proof}

\subsection{The de-looping theorem}

Let $\Sm\ds S$ be the subcategory of $\Sm/S$ with the same objects and with morphisms $\Hom_{\Sm\ds S}(Y,Y')$ the smooth $S$-morphisms $Y\to Y'$.

\begin{Def} (1) For $E\in\Spt(S)$, define the presheaf of spectra $\Omega_TE$ on $\Sm/S$ by 
\[
\Omega_TE(Y):=E^{Y\times0}(Y\times\P^1).
\]
The same formula defines $\Omega_TE$ in $\Spt(\Sm\ds S)$ for $E\in \Spt(\Sm\ds S)$.\\
\\
For $E\in \Spt(X)$, we define the presheaf of spectra $\Omega_{\P^1}E$ on $\Sm/X$ by 
\[
\Omega_TE(Y):=\fib[E(Y\times\P^1)\xrightarrow{i_\infty^*}E(Y\times\infty)]
\]
(2)  For $E\in\Spt(S)$,   define the functor $\Omega_{\P^1}E$ by 
\[
\Omega_{\P^1}E(X):= E(\P^1\wedge X_+)=\fib(E(X\times\P^1)\xrightarrow{\res}E(X\times\infty)).
\]
We use the same formula to define $\Omega_{\P^1}E\in\Spt(\Sm\ds S)$ for $E\in \Spt(\Sm\ds S)$.
\end{Def}

\begin{rems} (1) If $E$ is homotopy invariant and satisfies Nisnevic excision, the same holds for $\Omega_TE$ and $\Omega_{\P^1}E$.  

\medskip
\noindent
(2) The commutative diagram
\[
\xymatrix{
E(X\times\P^1)\ar[r]^{\res}\ar@{=}[d]&E(X\times(\P^1\setminus0))\ar[d]^{\res}\\
E(X\times\P^1)\ar[r]^{\res}&E(X\times\infty)
}
\]
gives us the homotopy fiber sequence
\[
\Omega_TE(X)\to \Omega_{\P^1}E(X)\to \fib(E(X\times(\P^1\setminus0))\xrightarrow{\res}E(X\times\infty)).
\]
If $E$is homotopy invariant, $\fib(E(X\times(\P^1\setminus0))\xrightarrow{\res}E(X\times\infty))$ is weakly contractible, hence the natural map $\Omega_TE\to \Omega_{\P^1}E$ is a weak equivalence.
\end{rems}

Besides the usual uses of localization (e.g., reducing problems to the case of fields) the localization theorem tells us how to  commute the operation $E\mapsto E^{(p)}(X,-)$ with the  $T$-loops functor $E\mapsto \Omega_TE$.

For $W\subset Y$ a closed subset of some $Y\in\Sm/S$, and for $E\in\Spt(\Sm\ds S)$, the spectrum with support $(\Omega_TE)^W(Y)$ is the iterated homotopy fiber over the diagram
\[
\xymatrix{
E(Y\times\P^1)\ar[r]\ar[d]&E((Y\setminus W)\times\P^1)\ar[d]\\
E(Y\times(\P^1\setminus0))\ar[r]&E((Y\setminus W)\times(\P^1\setminus0))
}
\]
Similarly, the spectrum with support $E^{W\times0}(Y\times\P^1)$ is  the iterated homotopy fiber over the diagram
\[
\xymatrix{
E(Y\times\P^1)\ar[r]\ar[d]&E(Y\times\P^1\setminus W\times0)\ar@{=}[d]\\
E(Y\times\P^1\setminus W\times0)\ar@{=}[r]&E(Y \times\P^1\setminus W\times0)
}
\]
The evident restriction maps yield a map of the second diagram to the first, and hence a canonical
map of spectra
\[
\theta_p^W(Y): E^{W\times0}(Y\times\P^1)\to (\Omega_TE)^W(Y);
\]
if $E$ satisfies  Zariski excision, then $\theta_p^W(Y)$ is a weak equivalence, since $(Y\setminus W)\times\P^1\cup Y\times(\P^1\setminus0)=
Y\times\P^1\setminus W\times0$.

\begin{Def} For $E\in \Spt(k)$, let $E^{(p)}\ds k$ be the presheaf on $\Sm\ds k$
\[
Y\mapsto E^{(p)}(Y,-),
\]
and let $E^{(\P^1)(p)}_{(0)}\ds k$ be the presheaf on $\Sm\ds k$
\[
Y\mapsto E^{(p)}_{Y\times 0}(Y\times\P^1,-).
\]
\end{Def}
The maps $\theta_p^W(Y\times\Delta^n)$ yield the map in $\Spt(\Sm\ds k)$
\[
\theta_E: E^{(\P^1)(p)}_{(0)}\ds k \to (\Omega_TE)^{(p-1)}\ds k.
\]

The sequence
\[
E^{(p)}_{Y\times 0}(Y\times\P^1,-)\to E^{(p)}(Y\times\P^1,-)\to
E^{(p)}(Y\times(\P^1\setminus0),-)
\]
gives rise to the map in $\Spt(\Sm\ds k)$
\[
\tau_E:E^{(\P^1)(p)}_{(0)}\ds k\to \Omega_T(E^{(p)}\ds k)
\]

\begin{thm}\label{thm:Delooping}  Suppose that $E\in\Spt(k)$ satisfies Zariski excision. Then the diagram
\[
\xymatrix{
&E^{(\P^1)(p)}_{(0)}\ds k\ar[rd]^{\tau_E}\ar[ld]_{\theta_E}\\
 (\Omega_TE)^{(p-1)}\ds k&&\Omega_T(E^{(p)}\ds k)
 }
 \]
 defines an isomorphism in $\sH\Spt(\Sm\ds k)$
 \[
 \xi_p:(\Omega_TE)^{(p-1)}\ds k\to \Omega_T(E^{(p)}\ds k)
 \]
 \end{thm}

\begin{proof} By the localization theorem, $\tau_E$ is a pointwise weak equivalence.
 Since the map $\theta_E$ is a pointwise weak equivalence if $E$ satisfies Zariski excision, the result follows.
 \end{proof}

 \section{Functoriality and Chow's moving lemma} Fix a field $k$. In this section, we discuss the extension of the presheaf $E^{(p)}\ds k$ on $\Sm\ds k$ to a presheaf on $\Sm/k$. 
\subsection{Functoriality}  Take  $E\in \Spt(k)$. Recall from the previous section the presheaf $E^{(p)}\ds k$, $X\mapsto E^{(p)}(X,-)$. Let $\rho:\Sm\ds k\to \Sm/k$ be the inclusion.

\begin{thm}\label{thm:Funct} Suppose that $E\in\Spt(k)$ is homotopy invariant and satisfies Nisnevic excision; if
$k$ is finite,  assume in addition
that $E$ satisfies the axiom A3. Then, 
\begin{enumerate}
\item For each $p\ge0$ there is a presheaf  $E^{(p)}\in \Spt(k)$, together with an isomorphism
\[
\phi_p:E^{(p)}\circ \rho \to E^{(p)}\ds k
\]
in  $\sH\Spt(\Sm\ds k)$. 

\item There are natural transformations $\xi_p:E^{(p)}\to E^{(p-1)}$, $p\ge0$, making the diagram
\[
\xymatrix{
E^{(p)}\circ\rho\ar[r]^{\phi_p}\ar[d]&  {E^{(p)}\ds k}\ar[d]\\
E^{(p-1)}\circ\rho\ar[r]_{\phi_{p-1}}&  {E^{(p-1)}\ds k}
}
\]
commute in  $\sH\Spt(\Sm\ds k)$. 
\item The presheaf $E^{(\P^1)(p)}_{(0)}\ds k$ and natural transformation 
\[
E^{(\P^1)(p)}_{(0)}\ds k\to E^{(\P^1)(p-1)}_{(0)}\ds k
\]
   extends as in (1) and (2) to a presheaf  $E^{(\P^1)(p)}_{(0)}\in \Spt(k)$ and natural transformation  $\xi_p:E^{(\P^1)(p)}_{(0)}\to E^{(\P^1)(p-1)}_{(0)}$, and the diagram of Theorem~\ref{thm:Delooping}  extends as in (2) to a diagram of weak equivalences
\[
\xymatrix{
&E^{(\P^1)(p)}_{(0)}\ar[rd]^{\tau_E}\ar[ld]_{\theta_E}\\
 (\Omega_TE)^{(p-1)}&&\Omega_T(E^{(p)})
 }
 \]
 intertwining the transformations $\xi_{p-1}$,  $\xi^{(0)}_{(\P^1)(p)}$ and $\Omega_T(\xi_{p})$.
Setting $\psi_p:=\tau_E\circ\theta_E^{-1}$, we have
isomorphisms 
\[
\psi_p:(\Omega_TE)^{(p-1)}\to\Omega_T(E^{(p)}), \ p\ge0,
\]
intertwining the transformations $\xi_{p-1}$ and $\Omega_T(\xi_{p})$ (here we set $(\Omega_TE)^{(-1)}:=(\Omega_TE)^{(0)}$, $\xi_{-1}=\xi_0$).  

The same result holds after  replacing $\Omega_T$ with $\Omega_{\P^1}$.
\end{enumerate}
Additionally,  $E^{(p)}$ is a bifibrant object in $\Spt_\Nis(\Sm/k)$ and the operation $E\mapsto (E^{(p)}, \phi_p, \xi_p,\psi_p)$ is natural in $E$ and preserves weak homotopy fiber sequences in $\Spt_\Nis(k)$.
\end{thm}

\begin{proof} The theorem follows  essentially from the main result of \cite{LevineChowMov}, with some modifications and extensions. 

Let $f:Y\to X$ be a morphism in $\Sm/k$. We have defined in \cite[\S7.4]{LevineChowMov} a homotopy coniveau tower on $X$ adapted to $f$:
\[
\ldots\to E^{(p)}(X,-)_f\to E^{(p-1)}(X,-)_f\to\ldots\to E^{(0)}(X,-)_f
\]
The simplicial spectrum $E^{(p)}(X,-)_f$ is defined using the support conditions $\sS^{(p)}_X(n)_f$ adapted to $f$:
\[
\sS^{(p)}_X(n)_f:=\{W\subset X\times\Delta^n\ |\ W\in \sS_X^{(p)}(n)\text{ and }(f\times\id)^{-1}(W)\in
\sS_Y^{(p)}(n)\},
\]
with
\[
E^{(p)}(X,n)_f:=\hocolim_{W\in \sS^{(p)}_X(n)_f}E^W(X\times\Delta^n).
\]

We have also defined a category $\sL(\Sm/k)$ with objects morphisms $f:Y\to X$ in $\Sm/k$; the operation $(f:Y\to X)\mapsto X$ defines a faithful functor $\sL(\Sm/k)\to\Sm/k$, making 
$\sL(\Sm/k)^\op$ a lax fibered category over $\Sm/k^\op$. In addition, sending $f:Y\to X$ to $E^{(p)}(X,-)_f$ defines a functor $E^{(p)}(-)_?$ on $\sL(\Sm/k)^\op$. Sending $X$ to $\holim_{\pi^{-1}(X)}E^{(p)}(-)_?$
gives a lax functor from $\Sm/k^\op$ to spectra, which is then regularized to an honest presheaf  by applying a type of homotopy colimit construction adapted from work of Dwyer-Kan (see \cite[\S7.3]{LevineChowMov}). The bifibrant replacement of this presheaf (for r the Nisnevic-local model structure on $\Spt(k)$ gives the desired presheaf $E^{(p)}$.

To make the whole homotopy coniveau tower functorial, replace the presheaf category $\Spt(k)$ with $\Spt(\Sm/k\times\N)$, where $\N$ is the sequence category $0\to 1\to\ldots\to n\to\ldots$. Taking $\N$ to be discrete, the Nisnevic topology on $\Sm/k$ induces a topology on $\Sm/k\times\N$. We proceed exactly as above, constructing a presheaf of spectra $\hat{E}^{(*)}$ on $\Sm/k\times\N$, and then take $E^{(*)}$ to be the functorial fibrant model for the Nisnevic-local model structure. This defines the natural transformations $\xi_p$.

The same approach, applied to the diagram in Theorem~\ref{thm:Delooping}, gives the natural tranformations $\psi_p$ such that the whole package has the desired compatibilities.
\end{proof}

\begin{rem} If $E\in\Spt_{S^1}(k)$ is fibrant, then $E$ is homotopy invariant and satisfies Nisnevic excision. It follows from the naturality in $E$ and the fact that $E\mapsto E^{(p)}$ preserves homotopy cofiber sequences that the operations  $E\mapsto (E^{(p)}, \phi_p, \xi_p,\psi_p)$ descends to exact functors, resp. natural transformations, on $\SH_{S^1}(k)$, at least if $k$ is an infinite field.

If $k$ is a finite field, we can consider the full subcategory $\SH_{S^1}(k)_\fin$ with objects those $E$ which satisfy axiom A3. It is obvious that $\SH_{S^1}(k)_\fin$ is a triangulated subcategory. In the case of  a finite base-field, we have the functors, resp. natural transformations, as above on  $\SH_{S^1}(k)_\fin$.
\end{rem}

A useful consequence of Theorem~\ref{thm:Funct}(3) is 

\begin{cor}\label{cor:Delooping} Take $E\in\Spt(k)$ satisfying the hypotheses of Theorem~\ref{thm:Funct}. Then the canonical map
\[
E^{(p)}\to E
\]
induces an isomorphism on taking $p$th loop spaces
\[
\Omega_T^pE^{(p)}\to \Omega^p_TE.
\]
\end{cor}

\begin{proof} The composition $\psi_1\circ\ldots\circ\psi_p$ gives the isomorphism $\Omega^p_TE\to
\Omega^p_TE^{(p)}$; it follows from the explicit construction of $\psi_p$ on $\Omega_T(E^{(p-1)}\ds k)$ that the composition 
\[
\Omega^p_TE\to
\Omega^p_TE^{(p)}\to \Omega^p_TE
\]
is the identity.
\end{proof}

\begin{rem}
To state the next result, we need to describe how one extends a presheaf $E\in\Spt(k)$ to Zariski localizations of $X\in\Sm/k$. Let $S=\{x_1,\ldots, x_n\}$ be a finite set of points in $X$, and let $\sO$ be the semi-local ring $\sO_{X,S}$. We set
\[
E(\sO):=\colim_{S\subset U\subset X}E(U)
\]
where $U$ runs over open subschemes of $X$ containing $S$. This defines  $E(F)$ for $F$ a finitely and separately generated field extension of $k$ by choosing a smooth model $X$ with $F\cong k(X)$.
\end{rem}

\begin{cor}\label{cor:LayerFunct} Under the hypotheses of Theorem~\ref{thm:Funct}, for integers $p,r\ge0$, there is a presheaf  $E^{(p/p+r)}\in\Spt(k)$ whose restriction to
$\Sm\ds k$ is isomorphic to $E^{(p/p+r)}(?,-):\Sm\ds k^\op\to\Spt$ in $\sH\Spt(\Sm\ds k)$.
In addition
\begin{enumerate} 
\item The functor $E^{(0/1)}$ is {\em birational}: The restriction map 
\[
E^{(0/1)}(X)\to E^{(0/1)}(k(X))
\]
 is a weak equivalence.
\item The functor $E^{(0/1)}$ is {\em rationally invariant}: If $F\to F(t)$ is a pure transcendental extension of fields (finitely and separably generated over $k$), then $E^{(0/1)}(F)\to E^{(0/1)}(F(t))$ is a weak equivalence.
\end{enumerate}
\end{cor}

\begin{proof}  The main statement  follows from Theorem~\ref{thm:Funct}. For (1), fix an irreducible $X\in\Sm/k$, and let $Z\to X$ be a proper  closed subset. We have the localization fiber sequence
\[
E^{(0/1)}_Z(X,-)\to E^{(0/1)}(X,-)\to E^{(0/1)}(X\setminus Z,-).
\]
with $E^{(0/1)}_Z(X,-)$ the cofiber of $E^{(1)}_Z(X,-)\to E^{(0)}_Z(X,-)$. Since each closed subset $W\subset Z\times\Delta^n$ has codimension at least one on $X\times\Delta^n$, the map
$E^{(1)}_Z(X,n)\to E^{(0)}_Z(X,n)$ is an isomorphism for each $n$. Thus $E^{(0/1)}_Z(X,-)=0$ in $\SH$ and $E^{(0/1)}(X,-)\to E^{(0/1)}(X\setminus Z,-)$ is a weak equivalence. (2) follows by taking limits.

For (2), the homotopy property implies that
\[
E^{(0/1)}(F,-)\to E^{(0/1)}(\A^1_F,-)
\]
is a weak equivalence. Since $E^{(0/1)}(\A^1_F,-)\to
E^{(0/1)}(F(t),-)$ is a weak equivalence by (1), the result is proved.
\end{proof}

\subsection{The purity theorem}  
Using the functoriality of the $E^{(p)}$, we can extend the de-looping theorem to a version of the Thom isomorphism.

Fix a scheme $X$ in $\Sm/k$ and an $E\in\Spt(k)$. We may restrict $E$ to $\Sm/X$, giving the presheaf $E_X:\Sm/X^\op\to\Spt$. If we have a closed subset $Z$ of $X$, we have the functor  
\[
(f:U\to X)\mapsto E^{f^{-1}(Z)}(U),
\]
which we denote by $E_X^Z$. If $f:Y\to X$ is a morphism in $\Sm/k$, we have the pushforward $f_*:\Spt(Y)\to \Spt(X)$, defined by
\[
f_*F(U\to X):=F(U\times_XY).
\]
Clearly, $f_*$ preserves weak equivalences, hence descends to 
\[
f_*:\sH\Spt(Y)\to \sH\Spt(X).
 \]

\begin{lem}\label{ThomIso} Let $i:Z\to X$ be a codimension $d$ closed embedding, with $X$ and $Z$ in $\Sm/k$. Suppose that $E:\Sm/k^\op\to \Spt$ is homotopy invariant  and satisfies Nisnevic excision, and that the normal bundle $N_{Z/X}$ is trivial. Then a choice of isomorphism $\phi:N_{Z/X}\cong Z\times\A^d$ determines a natural isomorphism in $\sH\Spt(X)$,
\[
\omega_\psi:E^Z_X \to i_*(\Omega^d_TE_Z),
\]
natural in $(Z,X,\phi)$.
\end{lem}

\begin{proof} Using Nisnevic excision, the inclusion $\A^d\to(\P^1)^d$ induces a natural isomorphism 
\[
\Omega^d_TE_Z(Y)\to (E_Z)^{Y\times0}(Y\times\A^d).
\]
for $Y\to Z$ in $\Sm/Z$. Letting $E_Z^{-\times0}(-\times\A^d)$ denote the presheaf
\[
Y\mapsto  (E_Z)^{Y\times0}(Y\times\A^d),
\]
we thus have the isomorphism
\[
\Omega^d_TE_Z \to (E_Z)^{-\times0}(-\times\A^d).
\]
in $\Spt(Z)$.

Let $s:Z\to N_{Z/X}$ be the zero-section, $p:N_{Z/X}\to Z$ the projection and denote the presheaf
on $\Sm/Z$. Taking a deformation to the normal bundle as in \cite{MorelVoev}  gives a natural isomorphism of $E^Z_X$ with  $i_*E^{s(Z)}(N_{Z/X})$ in    $\sH\Spt$. The chosen isomorphism $\phi:N_{Z/X}\cong Z\times\A^d$ sends $s(Z)$ over to $Z\times0$. As the deformation diagram is preserved by pullback with respect to a smooth $U\to X$, we actually have an isomorphism in $\sH\Spt(X)$, proving the result.
\end{proof}

This immediately yields

\begin{prop}\label{prop:Loops} Let $E:\Sm/k^\op\to \Spt$ be a homotopy invariant presheaf satisfying Nisnevic excision. Let $i:Z\to X$ be a codimension $d$ closed embedding in $\Sm/k$ such that the normal bundle $N_{Z/X}$ is trivial. Then for all $p\ge0$ we have isomorphisms in $\SH$:
\begin{align*}
&E^{(p)}_Z(X,-)\cong (\Omega_T^dE)^{(p-d)}(Z,-)\\
&E^{(p/p+1)}_Z(X,-)\cong (\Omega_T^dE)^{(p-d/p-d+1)}(Z,-),
\end{align*}
where, for $n<0$,  we set $(\Omega_T^dE)^{(n)}=(\Omega_T^dE)^{(0)}$ and $E^{(n/n+1)}=*$.
The isomorphisms may depend on the choice of trivialization of $N_{Z/X}$, but are natural in the category of closed embeddings $i$ with trivialization of $N_i$.
\end{prop}

We also have
\begin{cor} \label{cor:GerstenSS} Suppose $k$ is perfect. Let $X$ be in $\Sm/k$, and let $E$ be as in Proposition~\ref{prop:Loops}. For each $N\ge0$, there is a spectral sequence
\begin{multline*}
E^1_{p,q}(E):=\oplus_{x\in X^{(p)}}\pi_{p+q}\Omega_T^pE^{(N-p/N-p+s)}(k(x),-)\\\Longrightarrow
\pi_{p+q}E^{(N/N+s)}(X,-).
 \end{multline*}
\end{cor}

\begin{proof} This follows from the localization property Corollary~\ref{cor:LocThm} and Proposition~\ref{prop:Loops} by the usual limit process.
\end{proof}
 
 Finally, the birationality and rational invariance of $E^{(0/1)}$ enable us to prove an extended form of the purity isomorphism for the layers of the homotopy coniveau tower. Let $X$ be in $\Sm/k$ and let $W\subset X$ be a closed subset with $\codim_XW\ge d$. We let $W^0\subset W$ be the smooth locus of $W$ and $W^0(d)\subset W^0$ the union of those components of $W^0$ having codimension exactly $d$ on $X$.

\begin{cor}\label{cor:Gysin} Suppose $k$ is perfect.  Let$W\subset X$ be a closed subset of $X$ with $\codim_XW\ge d$, and let $U$ be a dense open subset of $W^0(d)$. Then there is a canonical isomorphism
\[
\sigma_d:(E^{(d/d+1)})^W(X)\xrightarrow{\sim} (\Omega_T^dE)^{(0/1)}(U).
\]
in $\SH$.
 \end{cor}

\begin{proof}  Let $X^0=X\setminus(W\setminus U)$. By  the localization property for 
$E^{(d/d+1)}$, the restriction
\[
(E^{(d/d+1)})^W(X)\to (E^{(d/d+1)})^{U}(X^0)
\]
is a weak equivalence, so we reduce to the case $X=X^0$, $W=U$. 

By considering the deformation to the normal bundle, we have a canonical isomorphism 
\[
(E^{(d/d+1)})^W(X)\cong (E^{(d/d+1)})^W(N)
\]
in $\SH$, where $N$ is the normal bundle of $W$ in $X$ and $W$ is included in $N$ by the zero section $i_0:W\to N$.  

Let $N^0:=N\setminus\{i_0(W)\}$ with projection $q:N^0\to W$. Using  Corollary~\ref{cor:LayerFunct} and the localization property for 
$E^{(d/d+1)}$ again, the pull-back by $q$ induces weak equivalences
\begin{align*}
&(E^{(d/d+1)})^W(N)\xrightarrow{q^*} (E^{(d/d+1)})^{N^0}(q^*N)\\
&(\Omega^d_TE)^{(0/1)}(W)\xrightarrow{q^*} (\Omega_T^dE)^{(0/1)}(N^0).
\end{align*}
Using the diagonal section $\delta:N^0\to q^*N^0\subset q^*N$ as 1, we have a canonical isomorphism
\[
\phi:q^*N\cong N^0\times \A^1.
\]
This in turn gives a canonical trivialization of the normal bundle of $i_0(N^0)$ in $q^*N$, hence a canonical isomorphism in $\SH$
\[
 (E^{(d/d+1)})^{N^0}(q^*N)\cong  (\Omega_T^dE)^{(0/1)}(N^0).
\]
This completes the construction.
\end{proof}

\section{Generalized cycles}  We use the results of the previous sections to give an interpretation of the layers in the homotopy coniveau tower.

\subsection{The semi-local $\Delta^*$} We recall that $\Delta^n$ has the vertices
$v_0,\ldots, v_n$, where $v_i$ is the closed subscheme defined by $t_j=0$ $j\neq i$.
For a scheme $X$, we let $\Delta^n_0(X)$ be the intersection of all open subschemes $U\subset X\times\Delta^n$ with $X\times v_i\subset U$ for all $i$.

\begin{rem} If $X$ is a semi-local scheme with closed points $x_1,\ldots, x_m$, then
$\Delta^n_0(X)$ is just the semi-local scheme $\Spec \sO_{X\times\Delta^n,S}$, where
$S$ is the closed subset $\{ x_i\times v_j\ |\ i=1,\ldots, m, \ j=0,\ldots, n\}$. In
particular $\Delta^n_0(X)$ is an affine scheme if $X$ is semi-local.
\end{rem}

For a scheme $T$, we let $\Delta^*_0(T)$ denote the cosimplicial ind-scheme 
$n\mapsto\Delta^n_0(T)$; if $T$ is semi-local, then $\Delta^*_0(T)$ is a cosimplicial
semi-local scheme.  For  $F$ a field, we write $\Delta^*_{0,F}$ for 
$\Delta^*_0(\Spec F)$ 

\subsection{Some vanishing theorems}  We fix a presheaf $E:\Sm/k^\op\to\Spt$. For this section, we assume that $E$ is homotopy invariant and  satisfies Nisnevic excision; if $k$ is finite, we also suppose that $E$ satisfies the axiom A3. We note that these hypotheses pass to $E^{(p)}$ and $E^{(p/p+r)}$ for all $p,r\ge0$. Finally, we assume that $k$ is perfect.

\begin{lem} \label{lem:Vanishing}   Let $F=E^{(p)}:\Sm/k^\op\to \Spt$ with $p>0$. Then for $X$ in $\Sm/k$, $F^{(0/1)}(X,-)$ is weakly contractible.
\end{lem}

\begin{proof}  Noting that $F^{(0/1)}(X,-)\cong F^{(0/1)}(\Spec k(X),-)$ in $\SH$ (Corollary~\ref{cor:LayerFunct}), we reduce to the case of a field $K$. In this case, we have $F^{(0/1)}(K,-)=E^{(p)}(\Delta^*_{0,K})$. Since each $\Delta^n_{0,K}$ is semi-local, and hence affine, it follows from the construction of the functor $E^{(p)}$ in \cite[\S 7]{LevineChowMov} that we have the natural weak equivalences of presheaves on $\Ord$ (i.e. simplicial spectra)
\[
[n\mapsto E^{(p)}(\Delta^n_{0,K}]\cong
[n\mapsto E^{(p)}(\Delta^n_{K,0},-)_{f_n}],
\]
where $f_n:\amalg\Delta^m_{0,K}\to\Delta^n_{0,K}$ is the disjoint union of the inclusions of faces $\Delta^m_{0,K}\to \Delta^n_{0,K}$.  

Thus, $F^{(0/1)}(K,-)$ is weakly equivalent to  the total space of the bisimplicial spectrum 
\[
(n,m)\mapsto E^{(p)}(n,m), 
\]
where $E^{(p)}(n,m)$ is the limit of the spectra with support 
\[
E^W(\Delta^n_{0,K}\times_K \Delta^m_K), 
\]
as $W$ runs over all closed subsets of $\Delta^n_{0,K}\times_K \Delta^m_K$ satisfying
\[
\codim_{F\times F'}(W\cap F\times F')\ge p
\]
for all faces $F'\subset\Delta^m_K$, $F\subset\Delta^n_{0,K}$.

For each $m$, we have the restriction to a face (say the face $t_{m+1}=0$) 
\[
\delta^*:E^{(p)}(-,m+1)\to E^{(p)}(-,m),
\]
with right inverse given by pull-back by the corresponding codegeneracy map
\[
\sigma^*:E^{(p)}(-,m)\to E^{(p)}(-,m+1). 
\]
By the same argument as for the homotopy property for $E^{(p)}(X,-)$ (see \cite[Theorem 3.3.5]{LevineChowMov}), one shows that 
$\sigma^*\circ\delta^*$ is homotopic to the identity, hence $\delta^*$ is a homotopy equivalence.

Thus, the inclusion $E^{(p)}(-,0)\to E^{(p)}(-,-)$ is a weak equivalence. However, if $W$ is an irreducible closed subset of $\Delta^n_{0,K}$ which intersects all faces in codimension $\ge p>0$, then in particular, $W$ contains no vertex of $\Delta^n_{0,K}$. Since $\Delta^n_{0,K}$ is semi-local with closed points the vertices, this implies that $W$ is empty, that is, $E^{(p)}(-,0)$ is weakly contractible, proving the result.
\end{proof}

\begin{prop} \label{prop:Vanishing}  Let $F=E^{(p)}:\Sm/k^\op\to\Spt$ with $p>0$.   Then $F^{(q/q+1)}$ is weakly contractible for all $0\le q<p$. Similarly, $(E^{(p/p+1)})^{(q/q+1)}$ is weakly contractible for $0\le q<p$.
\end{prop}

\begin{proof}  Since the operation $F\mapsto F^{(q/q+1)}$ is compatible with taking homotopy cofibers, the second assertion follows from the first. 

 For the first assertion,  the case $p=1$ is handled by Lemma~\ref{lem:Vanishing}; we prove the general case by induction on $p$.

Note that  Theorem~\ref{thm:Funct}(3) gives us the weak equivalence
\[
\psi^d:\Omega_T^dF\to (\Omega_T^dE)^{(p-d)}
\]
 in $\Spt(\Sm\ds k)$. Thus, by our inductive assumption, $(\Omega_T^dF)^{(q/q+1)}$ is weakly contractible
for $0\le q<p-d$. By the Gersten spectral sequence (Corollary~\ref{cor:GerstenSS}), this implies that the restriction map
\[
F^{(q/q+1)}(\Delta_K^n)\to F^{(q/q+1)}(\Delta_{0,K}^n)
\]
is a weak equivalence for $0\le q<p$, for all fields $K$ finitely generated over $k$ and for all $n$. Thus, we have the isomorphisms in $\SH$
\[
F^{(q/q+1)}(K)\cong (F^{(q/q+1)})^{(0)}(K,-)\cong (F^{(q/q+1)})^{(0/1)}(K,-)
\]
for $0\le q<p$. Applying  Lemma~\ref{lem:Vanishing} to $F$, $(F^{(q/q+1)})^{(0/1)}(K,-)$ is weakly contractible for $q>0$, which completes the proof.
\end{proof}

\begin{prop} \label{prop:Vanishing2}  Let $F=E^{(p/p+1)}:\Sm/k^\op\to\Spt$ with $p>0$. Then $F^{(p+r)}$ is weakly contractible for all $r>0$.
\end{prop}

\begin{proof} $F^{(p+r)}(X)$ is isomorphic in $\SH$ to the total space of the simplicial spectrum $F^{(p+r)}(X,-)$. $F^{(p+r)}(X,n)$ in turn is the limit of the spectra with support $F^W(X\times\Delta^n)$, where $W$ is  a closed subset of $X\times\Delta^n$ which, among other properties, has codimension $\ge p+r>p$. By Corollary~\ref{cor:Gysin}, it follows that  $F^{(p+r)}(X,n)$ is weakly contractible, whence the result.
\end{proof}

\subsection{The main result}  

\begin{thm}\label{thm:GenCyc} Let $k$ be a perfect field.   Let $E\in\Spt(k)$ be a homotopy invariant presheaf satisfying Nisnevic excision; if $k$ is finite, we also assume that $E$ satisfies the axiom A3.  Take integers $0\le p\le q$. Then
\begin{enumerate}
\item Applying the functor ${-}^{(q)}$ to the canonical map $E^{(p)}\to E$ induces a weak equivalence
\[
(E^{(p)})^{(q)}\to E^{(q)}.
\]
in  $\Spt(k)$.
\item Applying the natural transformation ${-}^{(p)}\to \id$ to   $E^{(q)}$ induces a weak equivalence
\[
(E^{(q)})^{(p)}\to E^{(q)}
\]
in  $\Spt(k)$.
\item We have a natural isomorphism in $\sH\Spt(k)$ 
\[
(E^{(q/q+1)})^{(p/p+1)}\cong \begin{cases} 0&\quad\text{for }q\neq p,\\
E^{(p/p+1)}&\quad\text{for }q=p.
\end{cases}
\]
\end{enumerate}
\end{thm}

\begin{proof} For (1), we apply ${-}^{(q)}$ to the tower
\[
E^{(p)}\to E^{(p-1)}\to \ldots\to E
\]
giving the tower
\[
(E^{(p)})^{(q)}\to (E^{(p-1)})^{q)}\to \ldots\to E^{(q)}
\]
with layers $(E^{(r/r+1)})^{(q)}$, $r=0,\ldots, p-1$, which all vanish by   Proposition~\ref{prop:Vanishing2}. For (2), we use the same argument, applying the tower of functors
\[
{-}^{(p)}\to {-}^{(p-1)}\to \ldots\to \id
\]
to $E^{(q)}$ and using Proposition~\ref{prop:Vanishing}.

For (3), the case $p> q$ follows from   Proposition~\ref{prop:Vanishing2}. The same argument as for (2), replacing $E^{(q)}$ with $E^{(q/q+1)}$ handles the case $p<q$ and shows that the map
\[
(E^{(q/q+1)})^{(q)}\to E^{(q/q+1)}
\]
induced by applying the natural transformation ${-}^{(q)}\to   \id$ to $E^{(q/q+1)}$ is an isomorphism. Since $(E^{(q/q+1)})^{(q+1)}$ is weakly contractible, the natural map $(E^{(q/q+1)})^{(q)}\to (E^{(q/q+1)})^{(q/q+1)}$ is also an isomorphism, completing the proof.
\end{proof}

\begin{cor}\label{cor:GenCyc}  Let $E:\Sm/k^\op\to\Spt$ be a presheaf satisfying the same hypotheses as in Theorem~\ref{thm:GenCyc} .  Let $X$ be in $\Sm/k$. Then $E^{(p/p+1)}(X)$ is naturally isomorphic in $\SH$ to the total spectrum of a simplicial spectrum $E^{(p/p+1)}_\loc(X,-)$, with
\[
E^{(p/p+1)}_\loc(X,n)\cong \coprod_{x\in X^{(p)}(n)}(\Omega_T^pE)^{(0/1)}(k(x))
\]
in $\SH$.
\end{cor}

\begin{proof} By Theorem~\ref{thm:GenCyc}, $E^{(p/p+1)}(X)$ is isomorphic (in $\SH$) to the total spectrum of the simplicial spectrum $n\mapsto (E^{(p/p+1)})^{(p/p+1)}(X,n)$. By  Corollary~\ref{cor:Gysin}, we have the isomorphism
\[
(E^{(p/p+1)})^{(p/p+1)}(X,n)\cong \coprod_{x\in X^{(p)}(n)}(\Omega_T^pE)^{(0/1)}(k(x)),
\]
in $\SH$, as desired.
\end{proof}

\section{Computations}
In this section, we consider a special type of theory $E:\Sm/k^\op\to\Spt$ for which the ``cells" $(\Omega_{\P^1}^pE)^{(0/1)}$ are particularly simple, namely, that  for a field $F$, $(\Omega_{\P^1}^pE)^{(0/1)}(F)$ is a $K(\pi,0)$ with $\pi=\pi_0((\Omega_{\P^1}^pE)(F))$. For such theories, we can define an associated cycle theory $\CH^p(-;E,n)$ which generalizes the higher Chow groups of Bloch. We show that $K$-theory is of this type, and thus recover the Bloch-Lichtenbaum/ Friedlander-Suslin spectral sequence \cite{BL, FriedSus} as our homotopy coniveau spectral sequence. This gives a new proof that this spectral sequence has the expected $E_2$-terms consisting of motivic cohomology. Motivic cohomology itself is also of this form, and, being the associated cycle theory of another theory, has a particularly simple spectral sequence.

\subsection{Well-connected theories}

\begin{Def}\label{Def:WellConn} Let $E:\Sm/k^\op\to\Spt$ be a functor satisfying the hypotheses of Thereom~\ref{thm:GenCyc}.  We call $E$ {\em well-connected} if 
\begin{enumerate}
\item for $X\in\Sm/k$ and $W\subset X$ a closed subset, $E^W(X)$ is -1-connected.
\item for $F$ a field finitely generated over $k$, $\pi_n((\Omega_T^dE)^{(0/1)}(F))=0$ for $n\neq 0$ and all $d\ge0$.
\end{enumerate}
\end{Def}

\begin{rem} Suppose $E$ satisfies part (1) of Definition~\ref{Def:WellConn}. Since   
$\Omega_TE(X)=E^{X\times0}(X\times\A^1)$, it follows that $\Omega_T^dE$ also satisfies (1) for all $d\ge0$, so $E$ is well-connected if and only if $\Omega_T^dE$ is well-connected for all $d\ge0$.
\end{rem}

\begin{lem}\label{lem:pi0WellConn} Suppose $E$ is well-connected. Let $F$ be a field finitely generated over $k$. Then the natural map
\[
\pi_0((\Omega_T^pE)(F))\to \pi_0((\Omega_T^pE)^{(0/1)}(F))
\]
is an isomorphism for all $p\ge0$.
\end{lem}

\begin{proof} We give the proof for $p=0$ to simplify the notation. Since $E(\Delta^n_{0,F})$ is -1-connected, we have the exact sequence
\[
\pi_0(E(\Delta^1_{0,F}))\xrightarrow{\delta_1^*-\delta^*_0}
\pi_0(E(\Delta^0_{0,F}))\to \pi_0(E^{(0/1)}(F))\to0.
\]
Similarly, we have the  surjections $\pi_0(E(\Delta^n_{F}))\to \pi_0(E(\Delta^n_{0,F}))$. Using the homotopy property, we find that the natural map
\[
p^*:\pi_0(E(F))\to  \pi_0(E(\Delta^n_{0,F}))
\]
is an isomorphism for all $n$, so the above exact sequence becomes
\[
\pi_0(E(F))\xrightarrow{0}
\pi_0(E(F))\cong \pi_0(E^{(0/1)}(F)).
\]
\end{proof}

\subsection{Cycles} Let $E:\Sm/k^\op\to\Spt$ be a well-connected theory. For $X\in\Sm/k$ and $W\subset X$ a closed subset, set
\[
z^p_W(X;E):=\oplus_{\substack{x\in X^{(p)}\\\bar{x}\subset W}} \pi_0(\Omega_T^pE(k(x))).
\]
We write $z^p(X;E)$ for $z^p_X(X;E)$.

Let $f:Y\to X$ be a morphism in $\Sm/k$, and let $W\subset X$ be a codimension $p$ closed subset such that $f^{-1}(W)$ has codimension $p$ on $Y$. Take $x\in X^{(p)}$ with $\bar{x}\subset W$, and let $y\in Y^{(p)}$ be a point in $f^{-1}(\bar{x})$.
We have the pull-back homomorphism 
\[
f^*_{y/x}: \pi_0(\Omega_T^pE)(k(x))\to  \pi_0(\Omega_T^pE)(k(y))
\]
defined as the sequence
\begin{multline*}
\pi_0(\Omega_T^pE)(k(x))\cong  \pi_0(\Omega_T^pE)^{(0/1)}(k(x))\cong  \pi_0(E^{(p/p+1)})^{\bar{x}}(X)\\
\xrightarrow{f^*} \pi_0(E^{(p/p+1)})^{f^{-1}\bar{x}}(Y)\xrightarrow{\rm res}
\pi_0(E^{(p/p+1)})^{y}(\Spec\sO_{Y,y})\\
\cong  \pi_0(\Omega_T^pE)^{(0/1)}(k(y))\cong \pi_0(\Omega_T^pE)(k(y)).
\end{multline*}

Taking the sum of the $f^*_{y/x}$ defines the pull-back
\[
f^*:z^p_W(X;E)\to z^p_{f^{-1}W}(Y;E)
\]
which is easily seen to be functorial. 

For $X\in\Sm/k$, we define the {\em higher cycles with $E$-coefficients}  as
\[
z^p(X;E,n):=\lim_{\substack{\to\\W\in\sS_X^{(p)}(n)}}z^p_W(X\times\Delta^n;E),
\]
forming the simplicial abelian group $n\mapsto z^p(X;E,n)$ and the associated complex $z^p(X;E,*)$.

\begin{Def} Let $X$ be in $\Sm/k$. The higher Chow groups of $X$ with $E$-coefficients are the groups
\[
\CH^p(X;E,n):= H_n(z^p(X;E,*)).
\]
\end{Def}

\begin{rem} The higher Chow groups of Bloch, $\CH^p(X,n)$, are defined without reference to an underlying cohomology theory $E$. Instead, one uses the usual cycle groups 
\[
z^p(X):=\oplus_{x\in X^{(p)}}\Z
\]
as the building blocks for the cycle complex $z^p(X,*)$, where the pull-back map $f^*$ is defined via Serre's intersection multiplicity formula.

The properties we have established for the spectra $E^{(p/p+1)}(X,-)$, namely: homotopy invariance, localization, extension to a functor,are  all based on the analogous properties for the complexes $z^p(X,-)$ ({\it cf.} \cite{AlgCyc, Loc, MixMot}). In the sequel, we will often identify $z^p(X,-)$ with the associated simplicial Eilenberg-Maclane spectrum, so as to enable a comparision with other simplicial spectra.
\end{rem}

\subsection{Well-connectedness}  We give an alternative description of this property.

Recall the category $\Ord$ with objects the finite ordered sets $[n]:=\{0<1<\ldots<n\}$ and morphisms the order-preserving maps of sets.  For each $n$, we have the $n+1$-cube diagram with objects the injective maps $[m]\to[n]$ in $\Ord$, with $m\le n$ (plus the inclusion of the empty set $\0\to[n]$) and maps being commutative triangles of injective maps. Denote this $n+1$-cube by $([n],\del)$. Similarly, we have the $n$-cube diagram of injective maps $[m]\to[n]$ which send 0 to 0; denote this by $([n],\del^+)$. We have the evident inclusion $i_{n+}:([n],\del^+)\to ([n],\del)$. Sending $[m]$ with $m\le n-1$ to $[m+1]$  by $i\mapsto i+1$ defines the functor $j_n:([n-1],\del)\to
([n],\del^+)$. We have as well the projections $\pi_n: ([n],\del)\to \Ord$, $\pi_n^+: ([n],\del^+)\to \Ord$ sending $[m]\to[n]$ to $[m]$;  the maps $i\mapsto i+1$ from $[m]$ to $[m+1]$ gives us the natural transformation $f_n:\pi_{n-1}\to \pi_{n+}\circ j_n$. 

Let $E:\Ord^\op\to\Spt$ be a simplicial spectrum. Define $E([n],\del)$ to be the iterated homotopy fiber of $E\circ\pi_n$ over the $n+1$-cube $([n],\del)$ and define $E([n],\del^+)$ similarly. The inclusion $i_n$ and natural transformation $E(f_n):E\circ \pi_{n+}\circ j_n\to E\circ \pi_{n-1}$ give rise to the homotopy fiber sequence
\[
E([n],\del)\xrightarrow{i_n^*}E([n],\del^+)\xrightarrow{f_n^*}
E([n-1],\del).
\]

Similarly, for a functor $E:\Ord^\op\to \Ab$, the diagram $([n],\del)$ defines an $n+1$-dimensional complex $E\circ \pi_n$ and we let $E([n],\del)$ be the associated (homological) total complex (with $E([n])$ in degree 0). We define $E([n],\del^+)$ similarly. We extend this construction to simplicial complexes $E:\Ord^\op\to C(\Ab)$ in the obvious way.  Explicitly,  for  $E:\Ord^\op\to \Ab$, the the complex $E([n],\del)$ is:
\[
E([n],\del)_{-m}:=\oplus_{g:[n-m]\to [n]}E([n-m]),
\]
and differential $d_{-m}:E([n],\del)_{-m+1}\to E([n],\del)_{-m}$ the signed sum of the maps
\begin{align*}
&E(f)_g:=E(f):(E([n-m+1]),g)\to (E([n-m]),g\circ f); \\
&f:[n-m]\to [n-m+1]\in\Ord^{inj},
\end{align*}
$E([n],\del^+)_{-m}$ is the quotient of $E([n],\del)_{-m}$ by the subgroup
 \[
\oplus_{g,\  g(0)\neq0}E([n-m]),
\]

Finally, for $E:\Ord^\op\to C(\Ab)$, let $E_*$ be the complex associated to the simplicial abelian complex $E$, i.e., $E_m=E([m])$ and $d_m:E_{m+1}\to E_m$ is the usual sum of the maps $(-1)^iE(\delta_i):E([m+1])\to E([m])$, where $\delta_i:[m]\to[m+1]\in\Ord^{inj}$ is the map which omits $i$.

\begin{lem}\label{lem:DoldKan} Let $E:\Ord^\op\to C(\Ab)$ be a functor with $E([n])$ -1-connected for all $n$. Then there is an exact sequence
\[
H_0(E([n+1],\del^+))\xrightarrow{E(\delta_0)}H_0(E([n],\del))\to
H_n(E_*)\to0
\]
for all $n\ge0$. If $m\mapsto H_0(E([m]))$ is the constant functor, then 
\[
H_0(E([n],\del))\to H_n(E_*)
\] 
is an isomorphism.
\end{lem}

\begin{proof} The second assertion follows from the first. Indeed, the degeneracy maps $\sigma_i^n:[n]\to[n-1]$ (where $\sigma^n_i$  is the unique surjective order-preserving map $f:[n]\to[n-1]$ with $f(i)=f(i-1)$) define a splitting of the complexes $E([m],\del^+)$, giving the exact sequence
\[
0\to H_m(E([n],\del^+))\to H_m(E([n]))\xrightarrow{\sum_iE(\delta_i )}
\oplus_{i=1}^nH_m(E([n-1])).
\]
Thus, the assumption that $m\mapsto H_0(E([m]))$ is constant
 implies that $H_0(E([n],\del^+))=0$.

The first assertion is an easy consequence of the Dold-Kan correspondence, and is proved in, e.g., \cite[Lemma 2.6]{Loc}.
\end{proof}

\begin{lem} \label{lem:DoldKan2} Let $E:\Ord^\op\to \Spt$ be a simplicial spectrum. Suppose that $E([n])$ is 0-connected for all $n$. Then the total spectrum $|E|$ is weakly contractible if and only if $\pi_0(E([n],\del))=0$ for all $n\ge0$.
\end{lem}

\begin{proof} The strongly convergent spectral sequence
\[
E^{a,b}_1=\pi_{-a}E([-b])\Longrightarrow \pi_{-a-b}|E|
\]
shows us that $|E|$ is 0-connected.  We may therefore replace $E$ with $\Z E:\Ord^\op\to C(\Ab)$, and we need to show that $H_0((\Z E)([n],\del))=0$ for all $n$ if and only if $H_N(E_*)=0$. 

As in the proof of Lemma~\ref{lem:DoldKan}, we have inclusions $H_m(E([n],\del^+))\to H_m(E([n]))$. In particular, $H_0(E([n],\del^+))=0$. The result now follows directly from Lemma~\ref{lem:DoldKan}.
\end{proof}

\begin{rem} Given $E\in\Spt(k)$ and a cosimplicial $k$-scheme $Y:\Ord\to\Sm/k$, we can consider the simplicial spectrum $E\circ Y$ and construct the spectra $(E\circ Y)([n],\del)$, $(E\circ Y)([n],\del^+)$. We denote these by $E(Y([n]),\del)$ and $E(Y([n]),\del^+)$, respectively. 

For example, taking $Y=\Delta^*_{F,0}$, we have the spectra $E(\Delta^n_{F,0},\del)$ and 
$E(\Delta^n_{F,0},\del^+)$.

We use a similar notation for a presheaf of complexes $E:\Sm/k^\op\to C(\Ab)$.
\end{rem}

\begin{prop}\label{prop:Ax5Rev} Let $E\in\Spt(k)$ be a homotopy invariant presheaf satisfying Nisnevic excision; if $k$ is finite assume in addition that $E$ satisfies the axiom A3.  Suppose that $E$ satisfies Definition~\ref{Def:WellConn}(1). Then $E$ is well-connected  if and only if
\[
\pi_0((\Omega_T^dE)(\Delta^n_{0,F},\del))=0
\]
for all $n\ge1$, all $d\ge0$ and all finitely generated field extensions $F$ of $k$.
\end{prop}

\begin{proof}  We have the spectral sequence
\[
E^1_{p,q}=\pi_{p+q}(\Omega_T^dE)(\Delta^p_{0,F})\Longrightarrow \pi_{p+q}(\Omega_T^dE)^{(0/1)}(F,-).
\]
Since $\pi_m((\Omega_T^dE)(\Delta^p_{0,F}))=0$ for $m<0$, we have the exact sequence
\[
\pi_0((\Omega_T^dE)(\Delta^1_{0,F}))\xrightarrow{d_1^{1,0}}\pi_0((\Omega_T^dE)(F))\to \pi_0((\Omega_T^dE)^{(0/1)}(F,-))\to0.
\]
As in the proof of Lemma~\ref{lem:pi0WellConn}, $\pi_0((\Omega_T^dE)(\Delta^1_{0,F}))=\pi_0((\Omega_T^dE)(F))$ and $d_1^{1,0}$ is the zero map, hence the map 
\[
\pi_0((\Omega_T^dE)(F))\to \pi_0((\Omega_T^dE)^{(0/1)}(F,-))
\]
 is an isomorphism.

Thus, to prove the proposition, we must show that 
$\pi_n(E^{(0/1)}(F,-))=0$ for all $n\ge1$ if and only if $\pi_0(E(\Delta^n_{0,F},\del))=0$ for all $n\ge1$.

By assumption, the presheaf $E$ is pointwise -1-connected. Let 
\[
\tau_{\ge1}E\to E 
\]
be the 0-connected cover of $E$ (defined pointwise), giving us the fiber sequence
\[
\tau_{\ge1}E\to E\to \pi_0E
\]
We have already seen that $\pi_0E(F)\to \pi_0E(\Delta^n_{0,F})$ is an isomorphism for all $n$, and that
\[
\pi_0(E(\Delta^*_{0,F}))\cong (\pi_0E)(\Delta^n_{0,F})\cong \pi_0(E(F)).
\]
Similarly, $\pi_m((\pi_0E)(\Delta^*_{0,F})=0$ for $m\neq0$. Thus the weak homotopy fiber sequence
\[
(\tau_{\ge1}E)(\Delta^*_{0,F})\to E(\Delta^*_{0,F})\to (\pi_0E)((\Delta^*_{0,F})
\]
shows us that $\pi_0(\tau_{\ge1}E)(\Delta^*_{0,F})=0$ and
\[
\pi_m(\tau_{\ge1}E)(\Delta^*_{0,F})\to \pi_m E(\Delta^*_{0,F})
\]
is an isomorphism for all $m\neq0$. In particular, $\pi_n(E^{(0/1)}(F,-))=0$ for all $n\ge1$ if and only if
$(\tau_{\ge1}E)(\Delta^*_{0,F})$ is weakly contractible.

Similarly, $\pi_m(\tau_{\ge1}E)(\Delta^n_{F,0},\del)\to \pi_mE(\Delta^n_{F,0},\del)$ is an isomorphism for all $m\neq -n$. Thus, we are reduced to showing that $(\tau_{\ge1}E)(\Delta^*_{0,F})$ is weakly contractible if and only if $\pi_0(\tau_{\ge1}E)(\Delta^n_{F,0},\del)=0$ for all $n\ge0$. This follows from Lemma~\ref{lem:DoldKan2}.
\end{proof}

\subsection{The case of $K$-theory} 

\begin{thm} The algebraic $K$-theory functor $K:\Sm/k\to\Spt$ is homotopy invariant, satisfies Nisnevic excision and  satisfies axiom A3. In addition, $K$ is well-connected.
\end{thm}

\begin{proof} We use the basic results on $K$-theory and $G$-theory proved in \cite {Quillen}. By Quillen's resolution theorem, we have a weak equivalence of $K$-theory and $G$-theory for regular schemes. Thus, the homotopy invariance for $G$-theory  gives the homotopy invariance property for $K$-theory on $\Sm/k$. Quillen's localization theorem yields the weak equivalence
\[
K^W(X)\sim G(W)
\]
for $W\subset X$ a closed subset, $X\in\Sm/k$, hence $K$-theory  satisfies Nisnevic excision.

For axiom A3, let $k\subset k'$ be a finite degree $n$ field extension. We have the pushforward map $\pi_*:K(X\times_kk')\to K(X)$ and pull-back $\pi^*:K(X)\to K(X\times_kk')$ with $\pi_*\pi^*=\times n$ on the $K$-groups $K_p(X)$. This readily implies A3.

 For the well-connectedness, property (1) of Definition~\ref{Def:WellConn}  follows by the localization theorem, since $G$-theory is -1-connected (for $U\subset W$ open, the map $G_0(W)\to G_0(U)$ is surjective).

For part (2), we use Weibel's homotopy $K$-theory, $KH$ \cite{Weibel}.  By Vorst \cite{Vorst}, the normal crossing divisor $\del\Delta^n_{0,F}\subset\Delta^n_{0,F}$ is $K_1$-regular, hence we have the isomorphism
\[
K_n(\del\Delta^n_{0,F})\cong KH_n(\del\Delta^n_{0,F}).
\]
for $n\le1$. Since $KH$ satisfies Mayer-Vietoris for unions of closed subschemes, we have the weak equivalence of $KH(\Delta^n_{0,F},\del)$ with the homotopy fiber of the restriction
\[
KH(\Delta^n_{0,F})\to KH(\del\Delta^n_{0,F})
\]
By the $K_1$-regularity, we thus have the exact sequence
\[
K_1(\Delta_{0,F}^n)\to K_1(\del\Delta_{0,F}^n)\to K_0(\Delta_{0,F}^n,\del)\to
K_0(\Delta_{0,F}^n)\to K_0(\del\Delta_{0,F}^n).
\]
 $\Delta_{0,F}^n$ is semi-local and affine, so let $R$ be the ring of functions on $\Delta_{0,F}^n$ and $I$ is the ideal defining $\del\Delta_{0,F}^n$. Since $R$ is semi-local, the restriction map
$\GL(R)\to \GL(R/I)$ is surjective. Since $\Delta_{0,F}^n$ is affine, we have surjections
\[
\GL(R)\to K_1(\Delta_{0,F}^n);\quad
\GL(R/I)\to K_1(\del\Delta_{0,F}^n).
\]
Also, $K_0(R)=\Z=K_0(R/I)$, so 
$K_0(\Delta_{0,F}^n,\del)=0$.
Using Proposition~\ref{prop:Ax5Rev}, we see that $K$-theory is well-connected.
\end{proof}

\begin{thm}\label{thm:KThySS} There is a natural isomorphism in $\SH$
\[
z^p(X,-)\cong K^{(p/p+1)}(X,-).
\]
\end{thm}

\begin{proof} By Theorem~\ref{thm:GenCyc}, we have the isomorphism in $\SH$
\[
K^{(p/p+1)}(X,-)\cong (K^{(p/p+1)})^{(p/p+1)}(X,-).
\]
By Corollary~\ref{cor:GenCyc}, we have the isomorphism in $\SH$
\[
(K^{(p/p+1)})^{(p/p+1)}(X,n)\cong \coprod_{x\in X^{(p)}(n)} (\Omega^p_TK)^{(0/1)}(k(x)).
\]
Since $K$ is well-connected and $K_0(k(x))=\Z$, $ (\Omega^p_TK)^{(0/1)}(k(x))$ is the Eilen\-berg-Maclane spectrum $K(\Z,0)$. Thus, we have the weak equivalence
\begin{equation}\label{eqn:WeakEquiCyc}
(K^{(p/p+1)})^{(p/p+1)}(X,n)\xrightarrow{\sim} K(z^p(X,n),0).
\end{equation}
It remains to see that the two sides agree as simplicial spectra. 

The map \eqref{eqn:WeakEquiCyc} is just the weak equivalence of $(K^{(p/p+1)})^{(p/p+1)}(X,n)$ with its 0th Postnikov layer.   Thus, we need only see that $z^p(X,-)$ and $\pi_0(K^{(p/p+1)})^{(p/p+1)}(X,-)$ agree as simplicial abelian groups.

For this, take $x\in X^{(p)}(n)$. We have the natural map
\[
G(\bar{x})\sim K^{\bar{x}}(X\times\Delta^n)\to (K^{(p/p+1)})^{(p/p+1)}(X,n)
\]
and, for each face map $g:\Delta^m\to\Delta^n$, the commutative diagram
\[
\xymatrix{
\pi_0G(\bar{x})\ar[r]\ar[d]_{g^*}&\pi_0(K^{(p/p+1)})^{(p/p+1)}(X,n)\ar[d]^{g^*}\\
\pi_0G(g^{-1}(\bar{x}))\ar[r]&\pi_0(K^{(p/p+1)})^{(p/p+1)}(X,m)
}
\]
Similarly, we have the surjection $\pi_0G(\bar{x})\to  \pi_0G(k(x))$, $\pi_0G(g^{-1}(\bar{x}))\to
\pi_0G(k(g^{-1}(\bar{x})))$ and the identifications
\begin{align*}
&\pi_0G(k(x))=z^p_{\bar{x}}(X\times\Delta^n)\\
&\pi_0G(k(g^{-1}(\bar{x})))=z^p_{g^{-1}(\bar{x})}(X\times\Delta^m)
\end{align*}
Since the pull-back on cycles is defined via Serre's intersection multiplicity formula and Serre's vanishing theorem, we have the commutative diagram
\[
\xymatrix{
\pi_0G(\bar{x})\ar[r]\ar[d]_{g^*}&z^p_{\bar{x}}(X\times\Delta^n)\ar[d]^{g^*}\\
\pi_0G(g^{-1}(\bar{x}))\ar[r]&z^p_{g^{-1}(\bar{x})}(X\times\Delta^m)
}
\]
with surjective rows.

Putting these two commutative diagrams together with the weak equivalence from Corollary~\ref{cor:GenCyc} gives the functoriality of the weak equivalence \eqref{eqn:WeakEquiCyc} with respect to the simplicial structure.
\end{proof}

\subsection{Bloch motivic cohomology}\label{sec: BlochMotCoh} As in Theorem~\ref{thm:Funct}, one can make the cycle complexes $z^p(X,*)$ functorial in $X\in\Sm/k$ (see \cite[\S9]{LevineChowMov}). Specifically, there are fibrant complexes of Nisnevic sheaves on $\Sm/k$, $\sZ^p$, whose image in the derived category of Nisnevic sheaves on $\Sm\ds k$ is isomorphic to the functor 
\[
X\mapsto z^p(X,*).
\]
Thus, for each $X\in \Sm/k$, we have the complex $\sZ^p(X)_*$, with a natural isomorphism $\sZ^p(X)_*\cong z^p(X,*)$ in $\bD^-_\Nis(X)$. The (shifted) homology of $\sZ^p(X)$ is the {\em Bloch motivic cohomology} of $X$:
\[
H^n(X,\Z(p)):=H_{2p-n}(\sZ^p(X))=H_{2p-n}(z^p(X,*)).
\]

We consider $X\mapsto \sZ^p(X)$ as a functor $\sZ^p:\Sm/k^\op\to \Spt$ by taking the associated Eilenberg-Maclane spectrum.  As in  Theorem~\ref{thm:Funct}, the localization theorem for the complexes $z^p(X,*)$ yields the natural weak equivalences
\[
\sZ^p\to \Omega_T(\sZ^{p+1}).
\]

\begin{thm}\label{thm:MotCohUnstable} For each $p\ge0$ and each $q\ge0$ we have the isomorphism in $\sH\Spt(k)$
\[
(\sZ^p)^{(q/q+1)}\cong\begin{cases} 0&\quad\text{for }q\neq p\\
\sZ^p&\quad\text{for }q= p.
\end{cases}
\]
\end{thm}

\begin{proof}  From Theorem~\ref{thm:KThySS}, we have the isomorphism in $\sH\Spt(k)$
\[
\sZ^p\cong  K^{(p/p+1)},
\]
giving the  isomorphism in $\sH\Spt(k)$
\[
(\sZ^p)^{(q/q+1)}\cong  (K^{(p/p+1)})^{(q/q+1)}.
\]
The result then follows from Theorem~\ref{thm:GenCyc}.
\end{proof}

\section{The slice tower for $S^1$-spectra}\label{sec:S1Slice}

Voevodsky has defined the  {\em slice tower} for both $\P^1$-spectra and $S^1$-spectra. In this section, we look at the version for $S^1$-spectra, and show that the slice tower agrees with the homotopy coniveau tower. In this section $k$ will be a perfect infinite field.

\subsection{The slice tower} We consider $\P^1$ as a pointed space over $k$ using $\infty$ as the base-point. For a pointed space $Z$ over $k$, we write $\Sigma_{\P^1}Z$ for $\P^1\wedge Z$ and extend this notation to spectra over $k$ in the evident manner.

 In analogy with a categorical construction of the classical Postnikov tower, Voevodsky considers the localizing subcategory $\Sigma^d_{\P^1}\SH_{S^1}(k)$ of the category $\SH_{S^1}(k)$ of $S^1$-spectra over $k$ generated by the $\P^1$-sus\-pensions $\Sigma^d_{\P^1}E$ for $E\in \SH_{S^1}(k)$. This forms the tower of localizing subcategories
\[
\ldots \subset \Sigma^{d+1}_{\P^1}\SH_{S^1}(k)\subset \Sigma^d_{\P^1}\SH_{S^1}(k)\subset\ldots\subset
\Sigma^0_{\P^1}\SH_{S^1}(k)=\SH_{S^1}(k).
\]
The inclusion $i_d: \Sigma^d_{\P^1}\SH_{S^1}(k)\to   \SH_{S^1}(k)$ admits a right adjoint $r_d: \SH_{S^1}(k)\to  \Sigma^d_{\P^1}\SH_{S^1}(k)$, and Voevodsky defines  $f_d:=i_dr_d$. For $E\in \SH_{S^1}(k)$ this yields the natural tower
\begin{equation}\label{eqn:S1SliceTower}
\ldots\to f_{d+1}E\to f_dE\to\ldots\to E
\end{equation}
called the {\em $S^1$-slice tower}. Clearly the map $f_dE\to E$ is universal for maps $F\to E$ with $F\in\Sigma^d_{\P^1}\SH_{S^1}(k)$.

Our main theorem is
\begin{thm}\label{thm:S1Slice} Suppose that $k$ is a perfect infinite field. For $E\in \SH_{S^1}(k)$ and $d\ge0$ an integer, $E^{(d)}$ is in
$\Sigma^d_{\P^1}\SH_{S^1}(k)$, and the map $\phi_d:E^{(d)}\to f_dE$ adjoint to the canonical map $E^{(d)}\to E$ is an isomorphism.
\end{thm}

\subsection{The splitting} To prove Theorem~\ref{thm:S1Slice}, we first need  to construct a splitting to the canonical map $\kappa_d:E^{(d)}\to E$, in case $E=\Sigma^d_{\P^1}F$ for some $F\in \SH_{S^1}(k)$. The construction relies on the de-looping theorem \ref{thm:Delooping}, as extended in Theorem~\ref{thm:Funct}(3).

The canonical map $\kappa_d:E^{(d)}\to E$ induces the map 
\[
\Omega(\kappa_d):\Omega^d_{\P^1}(E^{(d)})\to \Omega^d_{\P^1}E, 
\]
which is an isomorphism by Corollary~\ref{cor:Delooping}.

Recalling that $E=\Sigma^d_{\P^1}F$, we have the canonical map
\[
\iota:F\to \Omega^d_{\P^1}\Sigma^d_{\P^1}F=\Omega^d_{\P^1}E,
\]
adjoint to the identity on $\Sigma^d_{\P^1}F$. Taking adjoints again, the map $[\Omega(\kappa_d)]^{-1}\circ\iota$ induces the map in $\SH_{S^1}(k)$
\[
\iota_d:E\to E^{(d)}
\]
with $\kappa_d\circ\iota_d=\id$, giving us the desired splitting.

\subsection{The proof of Theorem~\ref{thm:S1Slice}} We proceed by a series of lemmas.

\begin{lem}\label{lem:S1Slice1} For $E\in \Sigma_{\P^1}^d\SH_{S^1}(k)$, the canonical map $\kappa_q:E^{(q)}\to E$ is an isomorphism for all $q\ge d$.
\end{lem}

\begin{proof} As $E\mapsto E^{(q)}(X,-)$ is compatible with filtered colimits in $E$, we may assume that $E=\Sigma_{\P^1}^dF$ for some $F\in\SH_{S^1}(k)$. The splitting $\iota_d$ defined above gives us the commutative diagram
\[
\xymatrix{
E\ar[rr]^\id\ar[dr]_{\iota_q}&&E\\
&E^{(q)}\ar[ur]_{\kappa_q}
}
\]
Apply the functor ${(-)}^{(p/p+1)}$. By Theorem~\ref{thm:GenCyc}
\[
(E^{(q)})^{(p/p+1)}\cong0
\]
for $0\le p<q$, hence $E^{(p/p+1)}\cong0$ for $0\le p<q$. Thus the layers in the coniveau tower
\[
E^{(q)}\to E^{(q-1)}\to\ldots\to E^{(0)}=E
\]
are all zero, hence $E^{(q)}\to E$ is an isomorphism.
\end{proof}

\begin{lem} \label{lem:S1Slice2} Let $W\subset Y$ be a closed subset of some $Y\in\Sm/k$. Suppose that $\codim_YW\ge d$, and take $E\in \SH_{S^1}(k)$. Then the map $(f_dE)^W(Y)\to E^W(Y)$ induced by the canonical map $f_dE\to E$ is an isomorphism.
\end{lem}

\begin{proof} We proceed by descending induction  on the codimension, starting with codimension $\dim Y+1$. Since $k$ is perfect, $W$ admits a stratification with smooth strata; this reduces us to the case of $W\subset Y$ a smooth closed subscheme of codimension $d$. Similarly, we may assume that $W$ has trivial normal bundle in $Y$.

We have the canonical isomorphism
\[
\pi_n(E^W(Y))\cong \Hom_{\SH_{S^1}(k)}(\Sigma^\infty_s(Y/Y\setminus W),\Sigma^{-n}E),
\]
and a similar description of $\pi_n((f_dE)^W(Y))$. By the Morel-Voevodsky purity theorem \cite[Theorem 2.23]{MorelVoev} and the triviality of the normal bundle $\nu_{W/Y}$, we have the isomorphisms in $\SH_{S^1}(k)$:
\[
\Sigma^\infty_s(Y/Y\setminus W)\cong \Sigma^\infty_s\Th(\nu_{W/Y}\cong
\Sigma^\infty_s(\Sigma^d_{\P^1}W_+)\cong 
\Sigma^d_{\P^1}\Sigma^\infty_sW_+.
\]

The universal property of $f_dE\to E$ shows that the induced map
\begin{multline*}
 \Hom_{\SH_{S^1}(k)}(\Sigma^d_{\P^1}\Sigma^\infty_sW_+,\Sigma^{-n}f_dE)\\
 \to
  \Hom_{\SH_{S^1}(k)}(\Sigma^d_{\P^1}\Sigma^\infty_sW_+,\Sigma^{-n}E)
\end{multline*}
 is an isomorphism, whence the result.
\end{proof}

\begin{lem} \label{lem:S1Slice3} For $E\in \SH_{S^1}(k)$,  and $p\ge d$, the map $(f_dE)^{(p)}\to E^{(p)}$ induced by applying ${-}^{(p)}$ to the canonical map $f_dE\to E$ is an isomorphism.
\end{lem}

\begin{proof} Fix an $X\in \Sm/k$. Then $(f_dE)^{(p)}(X)\to E^{(p)}(X)$ is isomorphic to the map on the total spectra of the simplicial spectra
\[
(f_dE)^{(p)}(X,-)\to E^{(p)}(X,-).
\]
Since
\[
E^{(p)}(X,n)=\hocolim_{W\in \sS_X^{(p)}(n)}E^W(X\times\Delta^r),
\]
and as all $W$ in the limit have codimension $\ge p\ge d$,   Lemma~\ref{lem:S1Slice2} shows that 
$(f_dE)^{(p)}(X,n)\to E^{(p)}(X,n)$ is a weak equivalence for all $n$, and hence 
$(f_dE)^{(p)}(X,-)\to E^{(p)}(X,-)$ is a weak equivalence.
\end{proof}

\begin{lem} \label{lem:S1Slice4} Take $E\in \SH_{S^1}(k)$. Then $E^{(d)}$ is in $\Sigma^d_{\P^1}\SH_{S^1}(k)$.
\end{lem}

\begin{proof} Clearly $f_dE$ is in $\Sigma^d_{\P^1}\SH_{S^1}(k)$. By Lemma~\ref{lem:S1Slice1}, $(f_dE)^{(d)}$ is in $\Sigma^d_{\P^1}\SH_{S^1}(k)$ as well. By Lemma~\ref{lem:S1Slice3}, the map
$(f_dE)^{(d)}\to E^{(d)}$ is an isomorphism, hence $E^{(d)}$ is in $\Sigma^d_{\P^1}\SH_{S^1}(k)$.
\end{proof}

\begin{proof}[Conclusion of the proof] Take $E\in \SH_{S^1}(k)$. By Lemma~\ref{lem:S1Slice4}, $E^{(d)}$ is in $\Sigma^d_{\P^1}\SH_{S^1}(k)$, giving us the map $\phi:E^{(d)}\to f_dE$ adjoint to $\kappa:E^{(d)}\to E$. Apply the natural transformation $\tilde\kappa:{-}^{(d)}\to\id$ to the commutative diagram
\[
\xymatrix{
E^{(d)}\ar[d]_{\phi}\ar[r]^\kappa&E\\
f_dE\ar[ru]_\tau
}
\]
giving us the map of commutative diagrams
\[
\xymatrix{
(E^{(d)})^{(d)}\ar[d]_{\phi^{(d)}}\ar[r]^{\kappa^{(d)}}&E^{(d)}\\
(f_dE)^{(d)}\ar[ru]_{\tau^{(d)}}
}
\quad\hbox{\lower25pt\hbox{$\xrightarrow{\tilde{\kappa}}$}}\quad
\xymatrix{
E^{(d)}\ar[d]_{\phi}\ar[r]^\kappa&E\\
f_dE.\ar[ru]_\tau
}
\]

The maps $\kappa^{(d)}:(E^{(d)})^{(d)}\to E^{(d)}$ and
$\tilde{\kappa}:(E^{(d)})^{(d)}\to E^{(d)}$ are isomorphisms by Theorem~\ref{thm:GenCyc}.  By Lemma~\ref{lem:S1Slice3},
$\tau^{(d)}$ is an isomorphism, hence $\phi^{(d)}$ is also an isomorphism. By Lemma~\ref{lem:S1Slice1},  $\tilde{\kappa}:(f_dE)^{(d)}\to f_dE$ is an isomorphism, so finally 
$\phi$ is an isomorphism.
\end{proof}

\subsection{An $S^1$-connectedness result}

A well-known result in topology states that if a pointed space $X$ is $n$-connected, then so is $\Omega\Sigma X$. The analogous statement for $\SH_{S^1}(k)$ with respect to the $\P^1$-suspension and 
-loops functors was conjectured by Voevodsky \cite{VoevSlice1bis}:

\begin{conj} \label{conj:S1Connectedness} Let $E$ be in $SH_{S^1}(k)$. If $E$ is  in $\Sigma^d_{\P^1}\SH_{S^1}(k)$, $\Omega_{\P^1}\Sigma_{\P^1}E$ is also in 
$\Sigma^d_{\P^1}\SH_{S^1}(k)$.
\end{conj}

This conjecture is a consequence of

\begin{thm} \label{thm:S1slice2} There is a natural isomorphism of endo-functors on $\SH_{S^1}(k)$:
\[
\Omega_{\P^1}\circ f_{d+1}\cong f_d\circ \Omega_{\P^1}
\]
\end{thm}

Indeed, to prove the Conjecture, it suffices to show that  $\Omega_{\P^1}$ maps $\Sigma^{d+1}_{\P^1}\SH_{S^1}(k)$ to $\Sigma^d_{\P^1}\SH_{S^1}(k)$ for all $d\ge0$. This follows directly from 
Theorem~\ref{thm:S1slice2} and the fact that $f_{d+1}$ is the identity on $\Sigma^{d+1}_{\P^1}\SH_{S^1}(k)$ and $f_d$ maps $\SH_{S^1}(k)$ to $\Sigma^d_{\P^1}\SH_{S^1}(k)$.

The theorem is in turn an immediate consequence of the existence of the natural isomorphisms $\phi_p:{(-)}^{(p)}\to f_p$ from Theorem~\ref{thm:S1Slice} and the natural isomorphism of  Theorem~\ref{thm:Funct}(3),
\[
\Omega_{\P^1}\circ {(-)}^{(p)}\cong {(-)}^{(p-1)}\circ \Omega_{\P^1}.
\]

\subsection{The connectivity conjecture} The $S^1$-connectedness result has a consequence for $\P^1$-spectra as well, which we state here even though we will not recall the  definition of the homotopy category of $\P^1$-spectra $\SH(k)$ until \S\ref{sec:P1Spectra}.

There are adjoint functors
\begin{align*}
&\Sigma^\infty_{\P^1}:\SH_{S^1}(k)\to \SH(k)\\
&\Omega^\infty_{\P^1}:\SH(k)\to \SH_{S^1}(k).
\end{align*}
The connectivity conjecture is
\begin{conj}[\cite{VoevSlice1bis}]\label{conj:P1Connectedness} Let $E$ be in $\Sigma_{\P^1}^n\SH_{S^1}(k)$. Then $\Omega_{\P^1}^\infty\Sigma^\infty_{\P^1}E$ is in $\Sigma_{\P^1}^n\SH_{S^1}(k)$.
\end{conj}

This conjecture is a direct consequence of the $S^1$-connectedness conjecture 
\ref{conj:S1Connectedness}, since $\Omega_{\P^1}^\infty\Sigma^\infty_{\P^1}E$ is represented by the homotopy colimit (in $\Spt(k)$) $\hocolim_{n\to\infty}\Omega^n_{\P^1}\Sigma^n_{\P^1}E$. Thus, Theorem~\ref{thm:S1slice2}  yields a proof of Conjecture~\ref{conj:P1Connectedness}.

\section{The $\P^1$-stable theory}\label{sec:P1Spectra} We now pass to the setting of $\P^1$-spectra. In this section, we recall the definition of the $\P^1$-spectra and the construction of the slice tower. With the help of the de-looping isomorphism of Theorem~\ref{thm:Funct}(3), we extend in the following section the homotopy coniveau filtration of $S^1$-spectra to a tower of $\P^1$-spectra. The identification of the homotopy coniveau tower with Voevodsky's slice tower for $S^1$-spectra easily extends to a similar identification for $\P^1$-spectra. We conclude by showing that the 0th slice of the sphere spectrum is motivic cohomology, and describing consequences of these results for the Atiyah-Hirzebruch spectral sequence.

\subsection{$\P^1$-spectra} We give three definitions to allow for some flexibility:

\begin{Def}\label{Def:P1OmegaSpec} A $\P^1$-$\Omega$-spectrum $\sE$ over $X$ is given by
\begin{enumerate}
\item A sequence $(E_0, E_1,\ldots)$, where
$E_j\in\Spt(X)$ is a homotopy invariant presheaf satisfying Nisnevic excision.
\item Weak equivalences in $\Spt(X)$, $\epsilon_n:E_n\to \Omega^1_{\P^1}E_{n+1}$, $n=0,1\ldots$
\end{enumerate}
Maps are maps of sequences respecting the maps in (2).  We denote the category of $\P^1$-$\Omega$-spectra over $X$ by $\Spt_{\P^1}^\Omega(X)$.
\end{Def}

For the next definition, we use the category $\Spc(X)$ of presheaves of spaces over $X$.

\begin{Def}\label{Def:P1Spec} A $\P^1$-spectrum $\sE$ over $X$ is given by
\begin{enumerate}
\item A sequence $(E_0, E_1,\ldots)$, where each
$E_j$ is a pointed space over $X$.
\item  Maps of spaces over $k$, $\epsilon_n:\Sigma_{\P^1}E_n\to E_{n+1}$, $n=0,1\ldots$
\end{enumerate}
Maps are maps of sequences respecting maps in (2).   We denote the category of $\P^1$-spectra over $X$ by $\Spt_{\P^1}(X)$.
\end{Def}

Finally, we have the category of {\em $(s,p)$-spectra over $X$}.
\begin{Def}\label{Def:spSpec}
The category $\Spt_{(s,p)}(X)$ has  objects sequences $\sE:=(E_0, E_1,\ldots )$ of presheaves $E_n\in \Spt(X)$, together with bonding morphisms $\epsilon_n: \Sigma_{\P^1} E_n\to E_{n+1}$. Maps are sequences of maps in $\Spt(k)$ respecting the bonding morphisms.  
\end{Def}
Clearly  $\P^1$-$\Omega$-spectra over $X$ form a full subcategory of $\Spt_{(s,p)}(X)$.

If $\sE=(E_0, E_1,\ldots)$ is a $\P^1$-spectrum, a $\P^1$-$\Omega$-spectrum or an $(s,p)$-spectrum, we have the suspensions 
\begin{align*}
&\Sigma_{\P^1}\sE:=(E_1, E_2,\ldots) 
&\Sigma_{\P^1}^{-1}\sE:=(\Omega^1_{\P^1}E_0, E_0,E_1,\ldots). 
\end{align*}

\subsection{Model structure and homotopy categories} We recall the category $\SH(k)$ and its relation to the three categories of spectra defined above. For details, we refer the reader to \cite{MorelLec, Morel}.

For an $(s,p)$-spectrum $(\sE,\epsilon_n)$ over $k$,  the $\epsilon_n$ induce, for each $Y\in \Sm/k$, the map
\[
\epsilon_n(Y): E_n(Y)\to E_{n+1}(\Sigma_{\P^1}Y_+).
\]
For $Y$ in $\Sm/k$, we have the bi-graded stable homotopy groups
\[
\pi^s_{a,b}(\sE(Y))=\lim_{n\to\infty}\pi^s_{a+2n}(E_n(\Sigma_{\P^1}^{n+b}Y_+)),
\]
using the maps $\epsilon_n(-)$ for the transition maps in the inductive system of homotopy groups. The $\pi^s_{a,b}(\sE(Y))$ form a presheaf of abelian groups on $\Sm/k$; we let $\pi^s_{a,b}(\sE)$ denote the associated Nisnevic sheaf.

A map $f:\sE\to \sF$ of $(s,p)$-spectra is called a weak equivalence if $f$ induces an isomorphism $f_*:\pi^s_{*,*}(\sE)\to \pi^s_{*,*}(\sF)$ on the homotopy sheaves. $f$ is a cofibration if  $f_0$ is a cofibration in $\Spt_\Nis(k)$, and for each $n\ge0$, the map
\[
E_{n+1}\coprod_{\P^1\wedge E_n}\P^1\wedge F_n\to F_{n+1}
\]
is a cofibration in $\Spt_\Nis(k)$. The fibrations are characterized by having the RLP with respect to trivial cofibrations. This gives us the model category of  $(s,p)$-spectra on $\Sm/k$. The homotopy category is denoted $\SH(k)$.

If $\sE=(E_0, E_1,\ldots)$ is a $\P^1$-spectrum, we can form the associated $(s,p)$-spectrum by 
taking the term-wise (simplicial) suspension spectra $(\Sigma^\infty_sE_0, \Sigma^\infty_sE_1,\ldots)$. If $\sE=(E_0, E_1,\ldots)$ is a $\P^1$-$\Omega$-spectrum, the maps $\epsilon_n:E_n\to\Omega_{\P^1}E_{n+1}$ induce by adjointness the maps $\epsilon'_n:\P^1\wedge E_n\to E_{n+1}$, forming the associated $(s,p)$-spectrum.  Similarly, if $\sE=((E_0, E_1,\ldots),\epsilon_n)$ is an $(s,p)$-spectrum, we have by adjointness the maps $\epsilon_n':E_n\to\Omega_{\P^1}E_{n+1}$; if $\sE$ is fibrant, this forms a $\P^1$-$\Omega$-spectrum.   Thus, we may pass from $\P^1$-spectra to $\P^1$-$\Omega$-spectra by first forming the associated $(s,p)$-spectrum, taking a fibrant model, and then using adjointness. We denote this functor by $\sE\mapsto\Omega^*_{\P^1}\sE$.

We can also take the simplicial 0-spaces of an $(s,p)$-spectrum or a $\P^1$-$\Omega$-spectrum, forming a $\P^1$-spectrum: 
\[
\sE=(E_0, E_1,\ldots)\mapsto  (\lim_m\Omega^m_sE_{0m}, \lim_m\Omega^m_sE_{1m},\ldots).
\]
Here $\Omega^m_s$ is the loop-space functor with respect to the simplicial structure. 

Via these functors, the model structure on $\Spt_{(s,p)}(k)$ induces model structures on the categories of $\P^1$-spectra and on $\P^1$-$\Omega$-spectra, and gives a Quillen equivalence of these three model categories. In particular, we can consider  $\P^1$-spectra and  $\P^1$-$\Omega$-spectra as objects in $\Spt_{(s,p)}(k)$ or in $\SH(k)$.  

\begin{rems} (1) The above definition of $\SH(k)$ is slightly different than the one given by  \cite{Morel}. First of all, Morel uses the suspension functor with respect to $(\A^1\setminus\{0\},1)$,  rather than $(\P^1,\infty)$. Secondly, the individual spaces occurring in the spectra $E_n$ are required to be sheaves of pointed simplicial sets rather than presheaves. Since $(\P^1,\infty)$ is homotopy equivalent to $S^1\wedge(\A^1\setminus\{0\},1)$, the two different choices of suspensions lead to Quillen equivalent model categories, and as the cofibrations and weak equivalences in the presheaf  category are defined stalk-wise,  using presheaves or sheaves also yield Quillen equivalent model categories. Thus, we may use the same notation $\SH(k)$ for the homotopy category.

\medskip
\noindent
(2) We have the functor $\Sigma^\infty_{\P^1}:\Spt_{S^1}(k)\to \Spt_{(s,p)}(k)$, defined by sending $E$ to the sequence $\Sigma^\infty_{\P^1}E:=(E, \Sigma^1_{\P^1}E, \Sigma^2_{\P^1}E\ldots)$, with the evident bonding maps. $\Sigma^\infty_{\P^1}$ is a left Quillen functor, with right adjoint the 0-space functor $\sE\mapsto \Omega^\infty_{\P^1}\sE$, where
\[
\Omega^\infty_{\P^1}\sE:=\lim_{n\to \infty}\Omega_{\P^1}^nE_n
\]
if $\sE=(E_0, E_1,\ldots)$. In particular, this shows that a weak equivalence $\sE\to \sF$ between fibrant objects  in  $\Spt_{(s,p)}(\Sm/k)$ induces a point-wise weak equivalence $f_n:E_n\to F_n$ on the various $S^1$-spectra. 
\end{rems}

\begin{exs}\label{exs:Spectra} (1) Each $X\in\Sm/k$ determines the {\em $\P^1$-suspension spectrum} 
\[
\Sigma^\infty_{\P^1} X_+:=
(X_+,\Sigma^1_{\P^1}X_+, \Sigma^2_{\P^1}X_+,\ldots)
\]
and the corresponding $\P^1$-$\Omega$-spectrum $\Omega^*_{\P^1}\Sigma^\infty_{\P^1} X_+$. For $X=\Spec k$, we write $S^0_k$ for $\Spec k_+$. We have the {\em $\P^1$-sphere spectrum} 
\[
\Sigma^\infty_{\P^1}S^0_k=:(S^0_k, \P^1_k,\ldots,\Sigma_{\P^1}^dS^0_k,\ldots)
\]
 and the associated $\P^1$-$\Omega$-spectrum $\1:=\Omega_{\P^1}^*\Sigma^\infty_{\P^1}S^0_k$.

\medskip
(3) For $X\in\Sm/k$, let $\Cyc^d(X)$ denote the set of effective cycles $W=\sum_in_iW_i$, $n_i>0$, in  
$(\P^1)^d\times X$, with each irreducible component $W_i$ finite over $X$, and dominating some component of $X$. This defines the Nisnevic sheaf  $X\mapsto \Cyc^d(X)$. The reduced version is the quotient 
\[
\RCyc^d(X):=\Cyc^d(X)/\sum_{i=1}^d p_{i*}(\Cyc^{d-1}(X)),
\]
where $p_{i*}$ is the map induced by the inclusion $p_i:(\P^1)^{d-1}\to(\P^1)^d$ defined by inserting $\infty$ in the $i$th spot.

We have the map $\P^1\to\Cyc^1$ defined by taking the graph of a map $f:X\to \P^1$. Taking the product over $X$ gives the map
\[
\RCyc^d(X)\wedge\RCyc^{d'}(X)\to \RCyc^{d+d'}(X),
\]
which thus gives us the bonding morphisms
\[
\P^1\wedge\RCyc^d(X)\to \RCyc^{d+1}(X)
\]
This structure defines the $\P^1$-spectrum $\HZ$; in characteristic zero, $\HZ$ is represented by the symmetric powers of $(\P^1)^{\wedge d}$ in the evident way.

The associated $\P^1$-$\Omega$-spectrum $\Omega_{\P^1}^*\HZ$ is equivalent to the sequence of Bloch motivic cohomology presheaves
\[
(\sZ^0, \sZ^1,\ldots, \sZ^d)
\]
with connecting maps the localization weak equivalences
\[
 \sZ^p\to \Omega_{\P^1}\sZ^{p+1}.
\]
(see \S\ref{sec: BlochMotCoh}).

A direct map relating the two constructions is given as follows: Send $W\in \Cyc^d(X\times\Delta^n)$ to the cycle $W\in z^d((\P^1)^d\times X,n)$, which we then restrict to $W^0\in z^d(\A^d\times X,n)$ (with $\A^1=\P^1\setminus\{\infty\}$). This gives a natural transformation 
\[
\RCyc^d(X\times\Delta^*)\to z^d(\A^d\times X,*)\sim z^d(X,*),
\]
which gives the direct relation. That this map gives a Zariski-local weak equivalence is proved by Voevodsky-Suslin \cite{FSV} assuming resolution of singularities, and in general by Voevodsky \cite{VoevCyc}.
\end{exs}

\subsection{The stable homotopy coniveau tower} For a $\P^1$-$\Omega$-spec\-trum $\sE=((E_0, E_1,\ldots),\epsilon_*)$, and integer $p$, set
$\fil_p\sE:=(E_0^{(p)}, E_1^{(p+1)},\ldots)$, where the maps $\epsilon_d$ are given by
the delooping weak equivalences of Theorem~\ref{thm:Funct}(3) 
\[
(E_d)^{(d+p)}\xrightarrow{\epsilon_d^{(d+p)}}(\Omega_{\P^1}E_{d+1})^{(d+p)}
\xrightarrow{\psi_{d+p+1}}
\Omega_{\P^1}(E_{d+1}^{(d+p+1)}).
\]
Here we have chosen a lifting of $\psi_p$ to a map in $\Spt(k)$, which we can do because $\Omega_{\P^1}(E_{d+1}^{(d+p+1)})$ is fibrant and $(\Omega_{\P^1}E_{d+1})^{(d+p)}$ is cofibrant.

\begin{rem} By replacing $E^{(p+d)}_d$ with the appropriate limit using the natural transformations $\theta$ of Thereom~\ref{thm:Funct}(3), we can use the natural transformations $\tau$ 
of Thereom~\ref{thm:Funct}(3) to give a canonical lifting of the $\psi_p$. We will assume that we have done this, so that the maps $\psi_p$ are now lifted to $\Spt(k)$, satisfying the compatibilities listed in Theorem~\ref{thm:Funct}(3) in $\Spt(k)$ rather than in the homotopy category $\sH\Spt(k)$.
\end{rem}

The natural maps $E_j^{(p+j)}\to E_j$ define the map of $\P^1$-$\Omega$-spectra
\[
\fil_p\sE\to \sE.
\]
Recall that $E^{(n)}=E^{(0)}$ for $n<0$.
 
We thus have the tower of $\P^1$-$\Omega$-spectra
\begin{equation}\label{eqn:StHCTower}
\ldots\to \fil_{p+1} \sE\to  \fil_p\sE\to \ldots\to  \fil_0\sE\to \fil_{-1}\sE\to \ldots\to \sE.
\end{equation}
We write $\fil_{p/p+r}\sE$ for the cofiber $\fil_{p+r}\sE\to \fil_p\sE$ and $\barfil_p\sE$ for $\fil_{p/p+1}\sE$.

\begin{rem}\label{rem:Suspension} For $q\ge 0$ and $p\ge0$, we have the identity 
\[
\Sigma^q_{\P^1}(\fil_p\sE)=\fil_{p+q}\Sigma^q_{\P^1}\sE;
\]
for $q\ge0$ and $p<0$, we have this identity in ``sufficiently large degree"; in any case, a weak equivalence. Similarly, for $q<0$, the localization weak equivalence
\[
\Omega_{\P^1}^{-q}(E_n^{(m)})\sim E_{n+q}^{(m+q)}
\]
gives a natural isomorphism
\[
\Sigma^q_{\P^1}(\fil_p\sE)\sim \fil_{p+q}\Sigma^q_{\P^1}\sE
\]
in $\SH(k)$.
\end{rem}

\section{The slice tower in $\SH(k)$}

Voevodsky \cite{VoevSlice1} defines the slice filtration in   $\SH(k)$ just as it is defined in $\SH_{S^1}(k)$, the main difference being that the filtration is infinite in both the postive and negative directions.

Let $\SH^\eff(k)$ be the smallest localizing subcategory of $\SH(k)$   containing all suspension spectra $\Sigma_{\P^1}^\infty X_+$ with $X\in \Sm/k$; this is the same as the smallest localizing subcategory containing all the $\P^1$-suspension spectra $\Sigma_{\P^1}^\infty E$ for $E\in \SH_{S^1}(k)$. For each integer $p$, let $\Sigma^p_{\P^1}\SH^\eff(k)$ denote the smallest localizing subcategory of $\SH(k)$ containing the $\P^1$-spectra $\Sigma^p_{\P^1}\sE$ for $\sE\in \SH^\eff(k)$. Voevodsky remarks that the inclusion $i_p:\Sigma^p_{\P^1}\SH^\eff(k)\to \SH^\eff(k)$ admits the right adjoint $r_p:\SH^\eff(k)\to \Sigma^p_{\P^1}\SH^\eff(k)$; setting $f_p:=i_p\circ r_p$, one has for each $\sE\in \SH(k)$ the functorial {\em slice tower}
\[
\ldots\to f_{d+1}\sE\to f_d\sE\to\ldots\to f_0\sE\to f_{-1}\sE\to\ldots \to \sE.
\]
As for the slice tower in $\SH_{S^1}(k)$, the map $f_p\sE\to \sE$ is universal for maps 
$\sF\to \sE$, $\sF\in \Sigma^p_{\P^1}\SH^\eff(k)$.  The cofiber of $f_{d+1}\sE\to f_d\sE$ is denoted $s_d\sE$.

Our main result is the identification of the slice tower with the stable homotopy coniveau tower \eqref{eqn:StHCTower}.

\begin{thm} \label{thm:Slice} Let $k$ be a perfect field. For $\sE\in \SH(k)$, $\fil_p\sE$ is in  $\Sigma^p_{\P^1}\SH^\eff(k)$, and the map $\fil_p\sE\to f_p\sE$ adjoint to $\fil_p\sE\to \sE$ is an isomorphism.
\end{thm}

The proof follows the same line as that of Theorem~\ref{thm:S1Slice}, and relies on the lemmas used in that proof. We fix a perfect field $k$.

\begin{rem} In discussing the $S^1$ slice tower, we required $k$ to be an {\em infinite} perfect field. The reason that we required $k$ to be infinite was to have the functor $E\mapsto E^{(p)}$ defined for all fibrant $E$ in $\Spt_{S^1}(k)$. For a fibrant $(s,p)$-spectrum $\sE=(E_0, E_1, \ldots)$ the presheaves $E_n$ are all 0-spectra of a fibrant $(s,p)$-spectrum. In particular, the $E_n$ satisfy axiom A3 in case $k$ is a finite field, and hence the operation $E_n\mapsto E_n^{(p)}$ is well-defined. Thus all the results of \S\ref{sec:S1Slice} can be applied in our setting without requiring that $k$ be infinite.
\end{rem}

\subsection{The proof of Theorem~\ref{thm:Slice}}

\begin{lem}\label{lem:Slice1} For $\sE\in \Sigma_{\P^1}^d\SH(k)$, the canonical map $\kappa_d:\fil_d\sE\to \sE$ is an isomorphism.
\end{lem}

\begin{proof} We may assume that $\sE=\Sigma_{\P^1}^d\Sigma^\infty_{\P^1}E$ for some $E\in \SH_{S^1}(k)$, i.e., writing $\sE=(\sE_0,  \sE_1,\ldots)$,
we have
 $\sE_n=\Sigma_{\P^1}^{n+d}E$,
where $\Sigma_{\P^1}^rE=\Omega_{\P^1}^{-r}E$ for $r<0$. Thus 
\[
\fil_d\sE=((\fil_d\sE)_0,(\fil_d\sE)_1,\ldots)
\]
with $(\fil_d\sE)_n=(\Sigma_{\P^1}^{n+d}E)^{(n+d)}$. The canonical map
$(\fil_d\sE)_n\to \sE_n$ is by Lemma~\ref{lem:S1Slice1}  an isomorphism in $\SH_{S^1}(k)$ for all $n$, whence the result.
\end{proof}

\begin{lem} \label{lem:Slice2} Take a fibrant $\sE$ in $\Spt_{(s,p)}(k)$, and write
\[
\sE=((E_0, E_1,\ldots), \epsilon_n:E_n\to\Omega_{\P^1}E_{n+1})
\]
with each $E_n\in  \Spt(k)$. For an integer $n$, write 
\[
f_n\sE=((f_n\sE)_0, (f_n\sE)_1,\ldots),
\]
where $f_n\sE$ is assumed to be fibrant in $\Spt_{(s,p)}(k)$.
 Let $X$ be in $\Sm/k$ and $W\subset X$ a closed subset with $\codim_XW\ge n+m$ for some integer $m\ge0$. Then the map
\[
(f_n\sE)_m^W(X)\to E_m^W(X)
\]
is a weak equivalence.
\end{lem}

\begin{proof} As in the proof of Lemma~\ref{lem:S1Slice2}, we may assume that $W$ is a smooth codimension $m+n$ closed subscheme of $X$, with trivial normal bundle $\nu$.  We note that, if $\sF=(F_0, F_1,\ldots)\in \Spt_{(s,p)}(k)$ is fibrant, then
\[
\Hom_{\SH(k)}(\Sigma^\infty_{\P^1}A\to \sF)\cong \Hom_{\SH_{S^1}(k)}(\Sigma^p_{\P^1}A, F_p)
\]
for all $A$ in $\Spt(k)$ and all $p\ge0$. Indeed,  $F_0$ is the 0-spectrum $\Omega^\infty_{\P^1}\sF$. Since $\Sigma_{\P^1}^\infty$ and $\Omega^\infty_{\P^1}$ are adjoint,  our claim is verified for $p=0$. For general $p$, we have
\[
\Hom_{\SH_{S^1}(k)}(\Sigma^p_{\P^1}A, F_p)\cong \Hom_{\SH_{S^1}(k)}(A, \Omega^p_{\P^1}F_p)
\]
and $F_0\cong \Omega^p_{\P^1}F_p$ since $\sF$ is fibrant.

Taking $A=\Sigma^n_{\P^1}W_+$, and using the universal property of $f_n\sE\to \sE$, we find that
\[
 \Hom_{\SH_{S^1}(k)}(\Sigma^{n+m}_{\P^1}W_+, (f_n\sE)_m)\to
  \Hom_{\SH_{S^1}(k)}(\Sigma^{n+m}_{\P^1}W_+, \sE_m)
\]
is an isomorphism for all $m\ge0$. But, as in the proof of Lemma~\ref{lem:S1Slice2},  a choice of a trivialization of $\nu$ gives a natural isomorphism
\[
 \Hom_{\SH_{S^1}(k)}(\Sigma^{n+m}_{\P^1}W_+, E)\cong E^W(X)
 \]
 for all fibrant $E\in \Spt(k)$, proving the result.
 \end{proof}

\begin{lem} \label{lem:Slice3} For $\sE\in \SH(k)$,  and $p\ge d$, the map $\fil_p(f_d\sE)\to \fil_p\sE$ induced by applying $\fil_p$ to the canonical map $f_d\sE\to \sE$ is an isomorphism.
\end{lem}

\begin{proof} If we write $\sE=(E_0,E_1,\ldots)$, $f_d\sE=((f_d\sE)_0, (f_d\sE)_1,\ldots)$, then
\begin{align*}
&\fil_p\sE=(E_0^{(p)},E_1^{(p+1)},\ldots)\\
&\fil_pf_d\sE=((f_d\sE)_0^{(p)}, (f_d\sE)_1^{(p+1)},\ldots).
\end{align*}
Using Lemma~\ref{lem:Slice2} and arguing as in the proof of Lemma~\ref{lem:S1Slice2}, we see that
\[
 (f_d\sE)_m^{(p+m)}(X,n)\to E_m^{(p+m )}(X,n)
 \]
 is a weak equivalence for all $X\in\Sm/k$ and all $m\ge0$. This yields the desired result.
 \end{proof}

\begin{lem} \label{lem:Slice4} Take $\sE\in \SH(k)$. Then $\fil_d\sE$ is in $\Sigma^d_{\P^1}\SH(k)$.
\end{lem}

\begin{proof} The proof is exactly the same as the proof of Lemma~\ref{lem:S1Slice4}, using Lemma~\ref{lem:Slice3} instead of Lemma~\ref{lem:S1Slice3}.
\end{proof}

Theorem~\ref{thm:Slice} is now proved exactly as was Theorem~\ref{thm:S1Slice}, using Lemmas \ref{lem:Slice1}-\ref{lem:Slice4} in place of Lemmas  \ref{lem:S1Slice1}-\ref{lem:S1Slice4}.

\section{The sphere spectrum and the $\HZ$-module structure}

In this section, we analyze the layer $\barfil_0\1$, and show that this spectrum is isomorphic to the motivic cohomology spectrum $\HZ$. By Proposition~\ref{prop:module}, this gives the  $\barfil_p\sE$ an $\HZ$-module structure, and shows that the $E_1$-terms in the homotopy coniveau  spectral sequence (Proposition~\ref{SpecSeqProp2}) may be interpreted as generalized motivic cohomology. Throughout this section, we assume that the base-field $k$ is perfect.

\subsection{The fundamental class of a system} We consider the following situation: For each $n\ge0$ we are given $X_n\in \Sm/k$ and  a closed subscheme $D_n\subset X_n$.  We have as well morphisms $i_n:X_n\to X_{n+1}$, $n\ge0$, and an integer $d\ge0$. We assume
\begin{enumerate}
\item $D_n$ is smooth over $k$ of pure codimension $d$ in $X_n$.
\item The diagram
\[
\xymatrix{
D_n\ar[r]\ar[d]_{i_n}&X_n\ar[d]^{i_n}\\
D_{n+1}\ar[r]&X_{n+1}
}
\]
is cartesian.
\end{enumerate}
This data gives us the following inverse system of spectra:
\[
\ldots\to (\1^{(d/d+1)})^{D_{n+1}}(X_{n+1})\xrightarrow{i_n^*}
 (\1^{(d/d+1)})^{D_n}(X_n)\to\ldots
 \]
 We will construct an element $[D_*]\in\pi_0\holim_n (\1^{(d/d+1)}_d)^{D_n}(X_n)$, which is natural in the system $(X_n,D_n)$, that is,  given a map of systems
 \[
  f_n:(Y_n,E_n)\to (X_n,D_n)
  \]
  such that  the diagram
  \[
\xymatrix{
D_n\ar[r]\ar[d]_{f_n}&X_n\ar[d]^{f_n}\\
E_n\ar[r]&Y_n
}
\]
is cartesian, then we have $f_*^*([D_*])=[E_*]$.

 To construct $[D_*]$, we start with the system in $\Sm/k$
 \[
 D_0\xrightarrow{i_0} D_1\xrightarrow{i_1}\ldots
 \]
 Letting $p_n:D_n\to \Spec k$ be the structure morphism, we have the system of maps of Nisnevic sheaves of pointed sets
 \[
 p_n:D_{n+}\to S^0_k=\Spec k_+
 \]
 with $p_{n+1}\circ i_n=p_n$. Taking $\Sigma^\infty_s$ and composing with the canonical maps in $\Spt(k)$
 \[
 \Sigma^\infty_sS^0_k\to \1_0\to \1_0^{(0/1)}
 \]
 yields the map in $\Spt(k)$
 \[
 p_{D_*}:\hocolim_n\Sigma^\infty_sD_{n+}\to \1_0^{(0/1)}.
 \]

Let  $\sHom_{\Spt(k)}(-,-)$ denote the $\Spt$-valued Hom-functor. We have the weak equivalences in $\SH$:
 \begin{multline*}
 \sHom_{\Spt(k)}(\hocolim_n\Sigma^\infty_sD_{n+},\1_0^{(0/1)})\\
 \to
 \holim_n\ \sHom_{\Spt(k)}(\Sigma^\infty_sD_{n+},\1_0^{(0/1)})
 \\\to
 \holim_n \1_0^{(0/1)}(D_n).
 \end{multline*}
 Applying the extended purity theorem Corollary\ref{cor:Delooping}, we have the isomorphisms
 \[
 \1_0^{(0/1)}(D_n)\to (\1_d^{(d/d+1)})^{D_n}(X_n),
 \]
in $\SH$, hence theisomorphism in $\SH$
 \[
 \holim_n \1_0^{(0/1)}(D_n)\to \holim_n  (\1_d^{(d/d+1)})^{D_n}(X_n).
 \]
Putting these all together  gives us the isomorphisms
 \begin{multline*}
  \pi_0(\holim_n(\1_d^{(d/d+1)})^{D_n}(X_n))\\
  \cong \pi_0(\sHom_{\Spt(k)}(\hocolim_n\Sigma^\infty_sD_{n+},\1_0^{(0/1)}))\\
 \cong
  \Hom_{\Spt(k)}(\hocolim_n\Sigma^\infty_sD_{n+},\1_0^{(0/1)}).
  \end{multline*}
  Thus, the map $p_{D_*}$ yields the element
  \[
  [\bar{D}_*]\in  \pi_0(\holim_n(\1_d^{(d/d+1)})^{D_n}(X_n))
  \]
  as desired. 
  
  The naturality in the system $X_n, D_n$ follows from the naturality of the de-looping weak equivalences.

\subsection{The reverse cycle map}  In this section, we show how to map $\HZ$ back to the layer $\barfil_0\1$. 

Let $\sE=(E_0, E_1,\ldots)$, $\epsilon_n:E_n\to\Omega_{\P^1}E_{n+1}$, be a $\P^1$-$\Omega$-spectrum. Define the $\P^1$-$\Omega$-spectrum $\sE^{(\A^*)}$ as follows: 
\[
\sE^{(\A^*)}=(E_0^{(\A^0)}, E_1^{(\A^1)}, \ldots, E_n^{(\A^n)},\ldots), 
\]
where $E_n^{(\A^n)}(Y):=E_n(Y\times\A^n)$. We have the $\A^1$-weak equivalences $p^*_n:E_{n+1}^{(\A^n)}\to E_{n+1}^{(\A^{n+1})}$ induced by the projection $p_n:\A^{n+1}\to \A^n$ on the first $n$ factors. The bonding maps $E_n^{(\A^n)}\to E_{n+1}^{(\A^{n+1})}$ are the composition
\[
E_n^{(\A^n)}\xrightarrow{\epsilon_n^{(\A^n)}}\Omega_{\P^1}E_{n+1}^{(\A^n)}\xrightarrow{p_n^*}
\Omega_{\P^1}E_{n+1}^{(\A^{n+1})}
\]
The maps $\pi_n^*:E_n\to E_n^{(\A^n)}$ induced by the projections $\pi_n:Y\times\A^n\to Y$ clearly define an $\A^1$-weak equivalence $\pi^*:\sE\to\sE^{(\A^*)}$.

We construct a map $\rcyc: \HZ\to \barfil_0\1^{(\A^*)}$ by first constructing maps
\[
\rcyc_d:\Sigma^\infty_s\HZ_d\to (\barfil_0\1)_d^{(\A^d)}
\]
in $\Spt(k)$, $d\ge1$, which we then patch together to yield the map $\rcyc$.

We first consider the case $d=1$. Identify $\Sym^m\P^1$ with $\P^m$ by noting that the symmetric homogeneous functions on $(\P^1)^m$ are the same as the homogeneous functions in variables $X_0,\ldots, X_m$. $\Cyc^1$ is represented by the union  
\[
\coprod_{m=1}^\infty\Sym^m\P^1=\coprod_{m=1}^\infty\P^m
\]
 via the incidence subvariety $D'_m\subset \P^1\times \P^m$. $D'_m$ is defined by the bihomogeneous polynomial $\sum_{i=0}^mX_iT_0^{m-i}T_1^m$; evidently, $D'_m$ is smooth over $k$. We identify $\A^1$ with $\P^1\setminus\infty$ as usual.

$\RCyc^1$ is represented by the system of pointed (by $\infty$) schemes
\begin{gather*}
\infty \xrightarrow{q_0} \Sym^1\P^1\xrightarrow{q_1} \Sym^2\P^1\xrightarrow{q_2}\ldots\xrightarrow{q_{m-1}} \Sym^m\P^m\xrightarrow{q_m}\ldots,\\
q_m(\hbox{$\sum_i$}x_i)=\hbox{$\sum_i$}x_i+\infty,
\end{gather*}
via the system of  incidence subvarieties $D_m\subset \A^1\times \P^m$,
\[
D_m:=D'_m\cap\A^1\times \Sym^m\P^1.
\]
Taking $X_m= \A^1\times \P^m$, we have the system $(X_m,D_m)$ satisfying the conditions of the previous section (for $d=1$), and thus the element
\[
[D_*]\in\pi_0\holim_m(\1_1^{(1/2)})^{D_m}(\A^1\times \Sym^m\P^1).
\]

We have the  ``forget supports" map
\[
(\1_1^{(1/2)})^{D_m}(X_m)=(\1_1^{(1/2)})^{D_m}(\A^1\times\Sym^m\P^1)
\xrightarrow{\alpha}\1_1^{(1/2)}(\A^1\times\Sym^m\P^1) 
\]
Thus  $[D_*]$ gives us the element
\[
[\Rcyc^1]\in \pi_0\holim_m(\1_1^{(1/2)})^{(\A^1)}(\Sym^m\P^1).
\]

Since the maps $\Sym^m\P^1\to \Sym^{m+1}\P^1$ are closed embeddings, the induced maps
in $\SH_{S^1}(k)$, $\Sigma^\infty_s\Sym^m\P^1\to \Sigma^\infty_s\Sym^{m+1}\P^1$, are cofibrations. Thus the canonical map
\[
\hocolim_m\Sigma^\infty_s\Sym^m\P^1\to \colim_m\Sigma^\infty_s\Sym^m\P^1
\]
is a weak equvalence. 
We therefore have the isomorphisms
\begin{multline*}
 \pi_0\holim_m(\1_1^{(1/2)})(\Sym^m\P^1)\\
 \cong \Hom_{\SH_{S^1}(k)}(\hocolim_m\Sigma^\infty_s\Sym^m\P^1, (\1_1^{(1/2)})^{(\A^1)})\\
 \cong \Hom_{\SH_{S^1}(k)}(\colim_m\Sigma^\infty_s\Sym^m\P^1, (\1_1^{(1/2)})^{(\A^1)}).
 \end{multline*}
Since the presheaf $\Sigma^\infty_s\RCyc^1$  is represented by  $\colim_m\Sigma^\infty_s\Sym^m\P^1$, the element $[\Rcyc^1]\in \pi_0\holim_m(\1_1^{(1/2)})^{(\A^1)}(\Sym^m\P^1)$ thus determines (up to homotopy) the map
\[
\Rcyc^1:\Sigma^\infty_s\RCyc^1\to  (\1_1^{(1/2)})^{(\A^1)}
\]
in $\Spt(k)$. Since $\RCyc^1=\HZ_1$ and $(\1_1^{(1/2)})^{(\A^1)}= (\barfil_0\1)_1^{(\A^1)}$, we have the map we wanted:
\[
\rcyc^1:\Sigma^\infty_s\HZ_1\to  (\barfil_0\1)_1^{(\A^1)}.
\]

For $d>1$,  let $W$ be in $\Cyc^d(X)$. We first consider the case of semi-local $X$, with a finite set of chosen points $x_1,\ldots, x_s$, to explain the idea of the construction. Then $P:=\cup_i W\cap (\P^1)^d\times x_i$ is a finite subset of $(\P^1)^d\times X$. Thus, there is a $k$-point $*$ in $\P^1(k)$ with $|W|\subset
(\P^1\setminus \{*\})^d\times X$, i.e., $W$ is a finite cycle on $\A^d\times X$.

Choosing a general linear projection $\pi:\A^d\to\A^1$, we map $W$ birationally to $\pi_*(W)$, and the set $P$ isomorphically to $\pi(P)$. We have the push-forward weak-equivalence
\[
\pi_*: \1_0^{(0/1)}(|W|)\to  \1_0^{(0/1)}(|\pi_*W|),
\]
defined as follows: Choosing suitable coordinates on $\A^d$ gives an isomorphism $\A^d\cong \A^1\times\A^{d-1}$ for which $\pi$ becomes identified with the projection on $\A^1$. We may therefore embed $\A^d$ as an open subset of $\A^1\times(\P^1)^{d-1}$. Let $\bar{\pi}:\A^1\times(\P^1)^{d-1}\to \A^1$ be the projection. We thus have the weak equivalences
\begin{multline*}
 \1_0^{(0/1)}(|W|):= (\1_0^{(d/d+1)})^{|W|}(\A^d\times X)\sim
(\1_0^{(d/d+1)})^{|W|}(\A^1\times X\times(\P^1)^{d-1})\\
\to (\1_0^{(d/d+1)})^{\bar{\pi}^{-1}(|\pi_*W|)}(\A^1\times X\times(\P^1)^{d-1})\sim
 (\1_0^{(1/2)})^{|\pi_*W|}(\A^1\times X)\\
=:\1_0^{(0/1)}(|\pi_*W|),
\end{multline*}
giving the definition of $\pi_*$.

Thus, the class $\rcyc^1(\pi_*W)$ gives the class $\rcyc^d(W)\in \1_0^{(0/1)}(|W|)$. This class is functorial with respect to restriction to the points $x_1$ $,\ldots,$ $x_s$.

To make this  canonical, let $K=k(\P^1)^d$, let $*\in(\P^1)^d(K)$ be the canonical point and let $(\P^1)^d_\infty=(\P^1)^d\setminus (\P^1\setminus\{\infty\})^d$. There is a unique isomorphism $\psi_1:\A^1_K\to \P^1_K\setminus *$ with $\psi_1(0)=(1:0)$, $\psi_1(1)=(1:1)$. We let $\psi_d=(\psi_1)^d:\A^d_K\to (\P^1)^d_K$
 be the resulting open immersion. Let $L\supset K$ be the field $K(x_1,\ldots,x_d)$ and let  $\pi_L:\A^d_L\to \A^1_L$ be the linear projection $\pi_L(y_1,\ldots, y_d)=\sum_ix_iy_i$. We identify $ (\P^1\setminus\{*\})^d_L$ with $\A^d_L$ via $\psi$ without further reference to $\psi$. We note that $L$ is a purely transcendental extension of $k$.

Take $X\in \Sm/k$ and let $W\in\Cyc^d(X)$ be an effective finite cycle, $q: (\P^1)^d\times X\to X$ the projection, $W_L^0=W_L\cap \A^d_L\times X$.  
For each point $x\in X$, $W_L^0\cap  \A^d_L\times x$ is dense in $W_L$ and  $\pi_L$ gives an isomorphism (of reduced schemes) from $|W_L^0|\cap\A^d_L\times x$ to its image $\pi(|W^0_L|\cap \A^d_L\times x)$ in 
$\A^1_L\times x$.

We apply the results of Corollary~\ref{cor:LayerFunct} and Corollary~\ref{cor:Gysin}, giving us the sequence of isomorphisms in $\SH$
\begin{align}
(\1_d^{(d/d+1)})^{|W|}(( \P^1)^d\times X)&\cong (\1_d^{(d/d+1)})^{|W_L|}(( \P^1)^d\times X)
\label{align:CycPushforward}\\
&\cong(\1_d^{(d/d+1)})^{|W_L^0|}(\A^d_L\times X)\notag\\
&\xrightarrow{\pi_*} (\1_1^{(1/2)})^{\pi(|W^0_L|)}(\A^1_L\times X)\notag
\end{align}

Let $\Cyc^1_{L.d}(X)$ be the set of finite cycles on $\P^1\times X_L$ of the form $\overline{\pi_*(W^0_L)}$ for some $W\in\Cyc^d(X)$, and let $\RCyc^1_{L.d}(X)$ be the quotient of $\Cyc^1_{L.d}(X)$ by the subset of cycles of the form $\overline{\pi_*(W^0_L)}$ for some $W\in\Cyc^d(X)$ which is supported in $(\P^1)^d_\infty$. This defines the presheaf  $\Cyc^1_{L.d}$ and quotient presheaf $\RCyc^1_{L.d}$  on $\Sm/k$. Letting $\Cyc^1_L$ be the presheaf $X\mapsto \Cyc^1(X_L)$, we have  the natural inclusion $\iota_d:\Cyc^1_{L.d}\to \Cyc^1_L$.  We note that sending $W\subset ( \P^1)^d\times X$ to  $\overline{\pi_*(W^0_L)}$ defines isomorphisms of presheaves
\begin{align*}
\bar{\pi}_*:\Cyc^d\to \Cyc^1_{L.d}\\
\pi_*:\RCyc^d\to \RCyc^1_{L.d}
\end{align*}

Define $(\1_d^{(d/d+1)})_\fin(X, ( \P^1)^d)$ to be the limit
\[
(\1_d^{(d/d+1)})_\fin(X, ( \P^1)^d)=\hocolim_W(\1_d^{(d/d+1)})^W(( \P^1)^d\times X)
\]
as $W$ runs over codimension $d$ closed subsets of $(\P^1)^d\times X$, finite over $X$, such that each irreducible component of $W$ dominates some component of $X$.  Sending $X$ to $(\1_d^{(d/d+1)})_\fin(X, ( \P^1)^d)$ defines the presheaf
\[
(\1_d^{(d/d+1)})_\fin(?, ( \P^1)^d):\Sm/k^\op\to\Spt.
\]
Define $(\1_d^{(d/d+1)})(X, \A^d)$  similarly, replacing $( \P^1)^d\times X$ with $\A^1\times X$ and $W$ with $W\cap \A^1\times X$; $\A^1=\P^1\setminus\{\infty\}$.

Let $(\1_1^{(1/2)})_\fin( X,\A^1)_{L,d}$ be the limit
\[
(\1_1^{(1/2)})_\fin( X,\A^1)_{L,d}:=\hocolim_{\substack{\to\\D}}(\1_1^{(1/2)})^D( \A^1\times X_L),
\]
where now $D$ runs over the closed subsets of the form $\overline{\pi(|W^0_L|)}$, where $W$ runs over all codimension $d$ closed subsets of $(\P^1)^d\times X$, with the finiteness and dominance conditions as above.  Let $(\1_1^{(1/2)})( X,\A^1)_{L,d}$ be the limit 
\[
(\1_1^{(1/2)})_\fin( X,\A^1)_{L,d}:=\hocolim_{\substack{\to\\D\supset D'}}(\1_1^{(1/2)})^{D\setminus D'}( \A^1\times X_L\setminus D'),
\]
where  $D$ is as above, and $D'$  runs over the closed subsets of the form $\overline{\pi(|W^0_L|)}$, where $W\subset (\P^1)_\infty$ has codimension $d-1$, with the finiteness and dominance conditions as above.

Sending $X$ to $(\1_1^{(1/2)})_\fin( X,\A^1)_{L,d}$ defines the presheaf
\[
(\1_1^{(1/2)})_\fin( ?,\P^1)_{L,d}:\Sm/k^\op\to\Spt.
\]
and the sequence of maps \eqref{align:CycPushforward} gives us the isomorphism
\[
\hat{\pi}_*^{(d)}:(\1_d^{(d/d+1)})_\fin(?, ( \P^1)^d)\to (\1_1^{(1/2)})_\fin(?,\A^1)_{L,d}
\]
in $\SH_{S^1}(k)$. Similarly, we have the presheaves $(\1_1^{(1/2)})( ?,\A^1)_{L,d}$ and
$(\1_d^{(d/d+1)})(?, \A^d)$ and  the weak equivalence
\[
\pi_*^{(d)}:(\1_d^{(d/d+1)})(?,  \A^d)\to (\1_1^{(1/2)})(?,\A^1)_{L,d}
\]
in $\SH_{S^1}(k)$.

View the element  $[\Rcyc^1]\in \pi_0\holim_m(\1_1^{(1/2)})^{(\A^1)}(\Sym^m\P^1)$ as a map
\[
[\Rcyc^1]: \colim_m \Sigma^\infty_s\A^1\times \Sym^m\P^1/(\A^1\times \Sym^m\P^1\setminus D_m)\to
\1_1^{(1/2)}.
\]
Each element $W$ of $\Cyc^1_{L.d}(X)$ gives us the element $\iota W\in \Cyc^1(X_L)$, and hence a morphism $\phi_W:X_L\to  \Sym^m\P^1_L$ for some $m$, with
\[
W=(\phi_W\times\id)^*(D_m).
\]
Thus we have 
\[
\id\times \phi:\A^1\times X/(\A^1\times X\setminus |W|)\to \A^1\times \Sym^m\P^1/(\A^1\times \Sym^m\P^1\setminus D_m)
\]
and hence we may compose $\phi_W\times\id$ with $[\Rcyc^1]$, giving the map of spectra
\[
[\Rcyc^1]\circ(\phi_W\times\id):\Sigma^\infty_s\A^1\times X/(\A^1\times X\setminus |W|)\to
\1_1^{(1/2)}
\]
This then defines the natural transformation
\begin{multline*}
\Sigma^\infty_s\Cyc^1_{L.d}(X)
\xrightarrow{\hat{\cyc}_{L.d}(X)} \\\Hom_{\Spt}(\holim_{W\in \Cyc^1_{L.d}(X)}\Sigma^\infty_s\A^1\times X/(\A^1\times X\setminus |W|), \1_1^{(1/2)}).
\end{multline*}
Composing with the natural map
\begin{multline*}
\Hom_{\Spt}(\holim_{W\in \Cyc^1_{L.d}(X)}\Sigma^\infty_s\A^1\times X/(\A^1\times X\setminus |W|), \1_1^{(1/2)})\\\to
(\1_1^{(1/2)})_\fin( X,\A^1)_{L,d},
\end{multline*}
we have thus defined the natural transformation
\[
\hat{\cyc}_{L.d}:\Sigma^\infty_s\Cyc^1_{L.d}\to (\1_1^{(1/2)})_\fin( ?,\A^1)_{L,d}.
\]
Composing with the isomorphism $\bar{\pi}_*:\Cyc^d\to \Cyc^1_{L.d}$ gives the natural transformation
\[
\cyc_{L.d}:\Sigma^\infty_s\Cyc^d \to (\1_1^{(1/2)})_\fin( ?,\A^1)_{L,d}.
\]

Similarly, the natural transformation $\cyc_{L.d}$ descends to the natural transformation
\[
\Rcyc_{L.d}:\Sigma^\infty_s\RCyc^d\to (\1_1^{(1/2)})( ?,\A^1)_{L,d}.
\]
Using the isomorphism $\pi_*^{(d)}$ defines the natural transformation
\[
\Rcyc^{(d)}_d:\Sigma^\infty_s\RCyc^1_{L.d}\to (\1_d^{(d/d+1)})(?, \A^d);
\]
forgetting the supports and changing notation gives us the map we wanted:
\[
\rcyc_d:\Sigma^\infty_s\HZ_d\to (\barfil_0\1)_d^{(\A^d)}.
\]

\subsection{The extension to a map in $\SH(k)$}

There is a minor lack of compatibility among the maps $\rcyc_d$, which we need to correct to yield a map of $(s,p)$-spectra. For this, we use a homotopy colimit construction.

Suppose we have $(s,p)$-spectra 
\begin{gather*}
\sE=(E_0, E_1,\ldots);\ \epsilon_n:\Sigma_{\P^1}E_n\to E_{n+1}\\
\sE'=(E_0', E_1',\ldots);\ \epsilon_n':\Sigma_{\P^1}E_n'\to E_{n+1}'
\end{gather*}
and maps
\[
\phi_n:E_n\to E_n'
\]
in $\Spt(k)$.  We consider the sequence of closed embeddings
\[
\Delta^0\xrightarrow{i_0}\Delta^1\xrightarrow{i_1}\ldots\xrightarrow{i_{n-1}}\Delta^n\xrightarrow{i_n}\ldots
\]
with $i_n(v_i^n)=v_i^{n+1}$, $i=0,\ldots, n$, where $v_i^n$ is the vertex $t_j=0$, $j\neq i$, of $\Delta^n$. Let 
\[
\tilde{E}_n:=\hocolim_{0\le j\le n} \Sigma^{n-j}_{\P^1}E_j\wedge \Delta^j,
\]
where the maps
\[
\Sigma^{n-j}_{\P^1}E_j\wedge \Delta^j\to  \Sigma^{n-j-1}_{\P^1}E_{j+1}\wedge\Delta^{j+1}
 \]
 in this diagram  are $\id_{(\P^1)^{n-j-1}}\wedge \epsilon_j\wedge i_j$. As the diagram $\P^1\wedge\Sigma^{n-j}_{\P^1}E_j\wedge\Delta^j$, $0\le j\le n$, is evidently a subdiagram of 
 $\Sigma^{n+1-j}_{\P^1}E_j\wedge \Delta^j$, $0\le j\le n+1$, we have the canonical map
 \[
 \tilde{\epsilon}_n:\P^1\wedge\tilde{E}_n\to \tilde{E}_{n+1},
 \]
 forming the $(s,p)$-spectrum $\tilde{\sE}$.
 We have as well the canonical map $p_n:\tilde{E}_n\to E_n$, giving the map of $(s,p)$-spectra
$p:\tilde{\sE}\to\sE$, which is evidently an $\A^1$-weak equivalence. Thus, in order to define a map $\phi:\sE\to \sE'$ in $\SH(k)$, which assembles the given maps $\phi_n:E_n\to E_n'$, it suffices to extend the collection of maps $\phi_n$ to a strictly compatible family of maps
\[
\tilde{\phi}_{n,j}:\Sigma^{n-j}_{\P^1}E_j\wedge \Delta^j\to \Sigma^{n-j}_{\P^1}E_j'.
\]

With this in mind, we proceed to the case of interest.

\begin{prop}\label{prop:CycleMap} The maps $\rcyc^d$ give rise to a  map  $\rcyc:\HZ\to
(\barfil_0\1)^{(\A^*)}$ in $\SH(k)$.
\end{prop}

\begin{proof}  The connecting maps for $(\barfil_0\1)^{(\A^*)}$  
\[
(\barfil_0\1)_d^{(\A^d)}\to \Omega_{\P^1}(\barfil_0\1)_{d+1}^{\A^{d+1}}
\]
are adjoint to  maps
\[
\P^1\wedge (\barfil_0\1)_d^{(\A^d)} \to (\barfil_0\1)_{d+1}^{\A^{d+1}}
\]
which in turn are induced by the natural maps
\[
\rho_d:\P^1(X)\times (\barfil_0\1)_d(X\times\A^d)\to (\barfil_0\1)_{d+1}(X\times\A^{d+1})
\]
defined by the following: Let $f:X\to \P^1$ be a morphism. The graph of $f$ gives the inclusions $X\times\A^d\to X\times\P^1\times\A^d$, which then gives the map
\[
f_*: (\1_d)^{(d/d+1)}(X\times\A^d)\to   (\1_{d+1})^{(d+1/d+2)}(X\times\P^1\times\A^d).
\]
Composing with the restriction 
\[
(\1_{d+1})^{(d+1/d+2)}(X\times\P^1\times\A^d)\to(\1_{d+1})^{(d+1/d+2)}(X\times\A^1\times\A^d)
\]
and using the canonical weak equivalences
\begin{align*}
&(\1_{d+1})^{(d+1/d+2)}(X\times\P^1\times\A^d)\sim  (\barfil_0\1)_d(X\times\A^d)\\
&(\1_{d+1})^{(d+1/d+2)}(X\times\A^1\times\A^d)\sim
(\barfil_0\1)_{d+1}(X\times\A^{d+1})
\end{align*}
 completes the definition of $\rho_d$.

The connecting maps for $\HZ$ are induced by maps
\[
\P^1(X)\times\Cyc^d(X)\to \Cyc^{d+1}(X)
\]
which are defined similarly, by taking the fiber product of a cycle $W\in\Cyc^d(X)$ with the graph of a map $f:X\to \P^1$ to define the resulting cycle $\Gamma_f\times_XW\in\Cyc^{d+1}(X)$.

We need to see that the maps $\rcyc^d$ are compatible with these connecting maps, up to a compatible family of $\A^1$ homotopies, as explained above.

For this, note that the projection $p_d$ on the last $d$ factors gives an isomorphism of $|\Gamma_f\times_XW|$ with $|W|$. Thus, the only difference between $\rho_d(f\times \rcyc^d(W))$ and  $\rcyc^{d+1}(\Gamma_f\times_XW)$ arises in the use of two different projections: $\pi_d\circ p_d$ and $\pi_{d+1}$, where $\pi_d$ and $\pi_{d+1}$ are the generic projections used in the definition of $\rcyc^d$ and $\rcyc^{d+1}$.

Let $L_{d+1}=K(x_1,\ldots, x_{d+1})$ be the field extension used to define $\pi_{d+1}$.  Clearly, we can define a family of linear projections
 \[
\pi_{d,1}:\A^{d+1}\times\Delta^1_{0,L_{d+1}}\to \A^1\times\Delta^1_{0,L_{d+1}}
\]
which agrees with $\pi_d\circ p_d$ at $(1,0)$ and $\pi_{d+1}$ at $(0,1)$, by making  the linear interpolation:
\[
\pi_{d,1}(y_1,\ldots, y_{d+1};t_0,t_1)=(t_0+t_1)\sum_{i=1}^dx_iy_i+t_1x_{d+1}y_{d+1}.
\]

 Using this family,   the rational invariance of $\1_0^{(0/1)}$ gives a homotopy between $\psi_d(f\times \rcyc^d(W))$ and $\rcyc^{d+1}(\Gamma_f\times_XW)$. Indeed, recall the canonical map
\[
[\Rcyc^1]: \colim_m \Sigma^\infty_s\A^1\times \Sym^m\P^1/(\A^1\times \Sym^m\P^1\setminus D_m)\to
\1_1^{(1/2)}),
\]
and the semi-local $n$-simplicies $\Delta^n_0$. 
The linear interpolation $\pi_{d,1}$ gives us the effective divisor
$W^1:=(\pi_{d,1},\id)_*(W\times\Delta^1_0)$ on $\A^1\times X\times\Delta^1_0$, which is classified by a morphism 
\[
\phi^1_W:X\times\Delta^1_0\to \lim_m\Sym^m\P^1,
\]
and hence the morphism
\begin{multline*}
\A^1\times X\times\Delta^1_0/(\A^1\times X\times\Delta^1_0\setminus |W^1|)\\
\xrightarrow{\id\times\phi^1_W}
 \A^1\times \Sym^m\P^1/(\A^1\times \Sym^m\P^1\setminus D_m).
\end{multline*}
Composing with $[\Rcyc^1]$ yields the morphism 
 \[
[\Rcyc^1]\circ(\phi^1_W\times\id):\Sigma^\infty_s\A^1\times X\times\Delta^1_0/(\A^1\times X\times\Delta^1_0/\setminus |W^1|)\to
\1_1^{(1/2)}
\]

But by extended purity (Corollary~\ref{cor:Gysin}), $[\Rcyc^1]\circ(\phi^1_W\times\id)$ extends canonically to a morphism
\[
\overline{[\Rcyc^1]\circ(\phi^1_W\times\id)}:\Sigma^\infty_s\A^1\times X\times\Delta^1/(\A^1\times X\times\Delta^1/\setminus\overline{|W^1|})\to
\1_1^{(1/2)}
\]
Threading this extension through the same process as we used to define $\rcyc_d$ yields the $\A^1$-homotopy for the two maps  $\P^1\wedge\Sigma^\infty_s\HZ_{d}\to  (\barfil_0\1)_{d+1}^{\A^{d+1}}$
in the diagram
\[
\xymatrix{
\Sigma^\infty_s\HZ_{d+1}\ar[r]^{\rcyc_{d+1}}& (\barfil_0\1)_{d+1}^{\A^{d+1}}\\
\P^1\wedge\Sigma^\infty_s\HZ_{d}\ar[r]_{\rcyc_{d}}\ar[u]^{\epsilon_n}& (\barfil_0\1)_{d}^{\A^{d}}\ar[u]_{\epsilon_n'}
}
\]

Similarly, for each $n$, we construct by linear interpolation a family of linear projections 
 \[
\pi_{d,n}:\A^{d+n}\times\Delta^n_{0,L_{d+n}}\to \A^1\times\Delta^n_{0,L_{d+n}}
\]
such that 
\[
\pi_{d,n}(v_i)=  \pi_{d+i}\circ p_{d+i}\circ p_{d+i+1}\circ\ldots\circ p_{d+n-1},
\]
where $v_i$ is the vertex $t_i=1$, $t_j=0$, $j\neq i$. By the same process as above, these give the  higher  $\A^1$-homotopies required to define the desired map of $\P^1$-spectra.
\end{proof}

\subsection{The cycle class map}\label{subsec:class}  We denote the Bloch motivic cohomology spectrum $(\sZ^0, \sZ^1,\ldots)$ by
$\sZ$.

\begin{lem}\label{lem:HZComp}
\[
\barfil_q\Sigma^p_{\P^1}\sZ\sim\begin{cases}
0&\quad\text{for } q\neq p\\
\Sigma^p_{\P^1}\sZ&\quad\text{for } q=p
\end{cases}
\]
\end{lem}

\begin{proof} By Remark~\ref{rem:Suspension}, it suffices to prove the case $p=0$. This follows directly from the identification
\[
\Omega^*_{\P^1}\HZ\cong(\sZ^0,\sZ^1,\ldots)
\]
and Theorem~\ref{thm:MotCohUnstable}.
\end{proof}

The canonical map $\1\to\sZ$ thus induces the map
\[
\cl:\barfil_0\1\to \barfil_0\sZ\cong
\sZ.
\]
On the zero-spaces, this is a natural transformation
\[
\cl^0:\1^{(0/1)}_0\to \sZ^0.
\]
Note that $\sZ^0$ is the constant sheaf $\Z$ (for the Zariski topology) on $\Sm/k$.

By the naturality of $\cl$, and the explicit description of the $d$th space $\1^{(d/d+1)}_d$ of $\barfil_0\1$ given by  Corollary~\ref{cor:GenCyc}, we find that, for $X\in\Sm/k$,  
\[
\cl^d:\1^{(d/d+1)}_d(X,-)\to z^d(X,-)
\]
is induced by the map on $n$-simplices
\[
\coprod_{x\in X^{(d)}(n)}\1^{(0/1)}_0(k(x))\xrightarrow{\coprod_x\cl^0}
\coprod_{x\in X^{(d)}(n)}\Z.
\]

Replacing $\sZ^d(X)$ with $\sZ^d(X\times\A^d)$, we have the modified spectrum $\Sigma\sZ^{(\A^*)}$, and the map $\cl:\barfil_0^{(\A^*)}\1\to  \sZ^{(\A^*)}$.

\subsection{$\barfil_0\1$ and $\HZ$}  The map $\rcyc:\HZ\to \barfil^{(\A^*)}_0\1$ extends canonically to the map on the associated $\P^1$-$\Omega$-spectrum $\rcyc:\Omega^*_{\P^1}\HZ\to \barfil^{(\A^*)}_0\1$.
  
\begin{thm}\label{thm:Comp} The maps $\cl:\barfil_0^{(\A^*)}\1\to  \sZ^{(\A^*)}$ and 
 $\rcyc:\Omega^*_{\P^1}\HZ\to \barfil^{(\A^*)}_0\1$ are  isomorphisms in $\SH(k)$.
\end{thm}

\begin{proof}  We first consider the composition $\psi$:
\[
\Sigma^\infty_{\P^1}S^0_k  \xrightarrow{\cl} \HZ \xrightarrow{\rcyc} \barfil^{(\A^*)}_0\1.
\]
Looking at the $d$th spaces gives the map
\[
\psi_d: \Sigma^d_{\P^1}S^0_k  \to (\barfil_0^{(\A^*)}\1)_d,
\]
 i.e.,  an element
\[
p_d \in \pi_0(\1_d^{(d/d+1)}(\Sigma^d_{\P^1}\A^d_+)).
\]

After correcting by the $\A^1$-homotopies defined in Proposition~\ref{prop:CycleMap}, we can compute $\psi_d$ as follows: Let $\delta\subset 
(\P^1)^d\times\A^d$ be the transpose of the graph of the standard inclusion $\A^d\to(\P^1)^d$. We have the de-looping isomorphism (in $\SH$)
\[
(\1_d^{(d/d+1)})^\delta((\P^1)^d\times\A^d) \cong \1_0^{(0/1)}(\delta).
\]
Also, since $\delta\cap (\P^1)^d_\infty=\0$, the canonical map
\[
(\1_d^{(d/d+1)})^\delta(\Sigma^d_{\P^1}\A^d_+)\to
(\1_d^{(d/d+1)})^\delta((\P^1)^d\times\A^d)
\]
is an isomorphism. We have as well the canonical map $\delta\to\Spec k$, which gives the canonical element $[\delta]\in \pi_0(S^0_k(\delta))$. The canonical map $\Sigma_s^\infty S^0_k\to
 \1_0^{(0/1)}$ with the above isomorphisms, the element $[\delta]$ yields the element $[\delta]^{(d/d+1)}\in \pi_0((\1_d^{(d/d+1)})^\delta(\Sigma^d_{\P^1}\A^d_+)))$. Forgetting supports, the element 
  $[\delta]^{(d/d+1)}$ maps to $p_d$.
  
We can go through the same procedure, replacing $\delta$ with the embedding $0^d\times  \A^d\to 
(\P^1)^d\times\A^d$, giving us the element $[0^d\times\A^d]^{(d/d+1)}]$ in  $\pi_0((\1_d^{(d/d+1)})^\delta(\Sigma^d_{\P^1}\A^d_+)))$, and, after forgetting supports, the element $p_d'\in 
 \pi_0(\1_d^{(d/d+1)}(\Sigma^d_{\P^1}\A^d_+))$. The graph of the map
 \begin{align*}
& \iota:\A^d\times\A^1\to (\P^1)^d
 &\iota(x_1,\ldots, x_d,t)=((1:tx_1),\ldots,(1tx_d))
 \end{align*}
 similarly gives an element $[\iota]\in  \pi_0(\1_d^{(d/d+1)}(\Sigma^d_{\P^1}\A^d\times\A^1_+))$ which defines an $\A^1$-homotopy between $p_d$ and $p_d'$, i.e. $p_d=p_d'$.
 
On the other hand, we have the canonical map $\1\to \sigma_0\1$, which on the $d$th spectrum level is the canonical map $\kappa_d:\Sigma^d_{\P^1}S^0_k\to \1_d^{(d/d+1)}$, which is the adjoint to the canonical element $1\in \pi_0(\1_0^{(0/1)}(k))$. It is clear from this description that $p_d'$ is the element of 
$\pi_0(\1_d^{(d/d+1)}(\Sigma^d_{\P^1}S^0_k)$ corresponding to $\kappa_d$, after pull-back by $\Sigma^d_{\P^1}\A^d_+\to \Sigma^d_{\P^1}S^0_k$. Therefore, after identifying
$ (\barfil_0^{(\A^*)}\1)_d$ with $ (\barfil_0\1)_d$  by this pull-back, 
$\psi_d=\kappa_d$. Thus the canonical extension of $\psi$ to $\tilde\psi:
\barfil_0\1\to \barfil^{(\A^*)}_0\1$ is the weak equivalence  given by projecting $X\times\A^*\to X$.

We now consider the composition $\phi(X)$:
\[
\HZ(X)\xrightarrow{\rcyc} \barfil^{(\A^*)}_0\1(X)\xrightarrow{\cl}
\sZ^{(\A^*)}(X).
\]
From the explicit description of $\cl$ given in \S\ref{subsec:class}, we  see that $\phi(X)$ is given by the map (on the $d$th space) which associates to a  cycle $W$ on 
$X\times\Delta^n\times(\P^1)^d$ the restriction to $X\times\Delta^n\times\A^d$. This map is the weak equivalence described above in Example~\ref{exs:Spectra}(2). Thus the extension of $\phi$ to $\phi:\Omega^*_{\P^1}\HZ\to \sZ^{(\A^*)}$ is an isomorphism in $\SH(k)$.
\end{proof} 

\begin{cor}\label{cor:SModule} Assume Conjecture~\ref{conj:ring}. For each  $\sE\in\SH(k)$  there is a natural $\HZ$-module structure on $s_n\sE$.
\end{cor}

\begin{proof}  The unit $\1\to \HZ$ lifts canonically to the unit $\1^\Sigma\to\HZ$ in $\Spt^\Sigma_{\P^1}(k)$.  By Proposition~\ref{prop:module}, $s_n\sE$ has the structure of a $s_0\1^\Sigma$-module, while by  Theorem~\ref{thm:Slice} and Theorem~\ref{thm:Comp}  the unit  induces an equivalence of the category of $s_0\1^\Sigma$-modules with $\HZ$-modules.
\end{proof}

\section{The motivic Atiyah-Hirzebruch spectral sequence}

We collect our results on the homotopy coniveau spectral sequence. For the results on $\DM(k)$ and $\SH(k)$ we use in this section, we refer the reader to the lectures of Morel \cite{MorelLec, Morel} and Voevodsky \cite{VoevSlice1bis}, as well as the papers of Ostvar-Roendigs \cite{OstvarRoendigs, RoendigsBigMot} and Spitzweck  \cite{Spitzweck}. Throughout this section, we are assuming the validity of Conjecture~\ref{conj:ring} described below.

\subsection{Products}   Let $\sE$  and $\sE'$ be $\P^1$-$\Omega$-spectra. The canonical maps $f_n\sE\to \sE$ and $f_m\sE'\to \sE'$ induce the map 
$\mu:f_n\sE\wedge f_m\sE'\to \sE\wedge \sE'$. $f_n\sE$ is in $\Sigma^n_{\P^1}\SH^\eff(k)$ and $f_m\sE$ is in $\Sigma^m_{\P^1}\SH^\eff(k)$, hence $f_n\sE\wedge f_m\sE'$ is in  $\Sigma^{n+m}_{\P^1}\SH^\eff(k)$. 

Applying $f_{n+m}$ to $\mu$ and using the universal property of $f_{n+m}$, we have the diagram
\[
\xymatrixcolsep{40pt}
\xymatrix{
f_{n+m}(f_n\sE\wedge f_m\sE')\ar[r]^-{ f_{n+m}(\mu)}\ar[d]_{i_{n+m}}&
f_{n+m}(\sE\wedge \sE')\\
f_n\sE\wedge f_m\sE'
}
\]
with $i_{n+m}$ a weak equivalence. Thus, we have the multiplication
\[
\mu_{n,m}: f_n\sE\wedge f_m\sE'\to f_{n+m}(\sE\wedge \sE')
\]
One checks that the $\mu_{*,*}$ are associative in $\SH(k)$ and are compatible with respect to increasing $n$ and $m$. Passing to the layers, we have the functor $\oplus_n$  on 
$\SH(k)$ and associative, graded natural transformation
\[
\mu_{*,*}:\oplus_ns_n\wedge \oplus_ns_n\to \oplus_ns_n\circ \wedge.
\]

Jardine \cite{Jardine2} has defined a model category of {\em symmetric spectra} over $k$  (with the stable model structure),  $\Spt^\Sigma_{\P^1}(k)$. The objects are similar to those in the category of $\P^1$-spectra, i.e., sequences of presheaves of spectra $(E_0, E_1,\ldots)$ on $\Sm/k$, with bonding maps $\P^1\wedge E_n\to E_{n+1}$, with the extra data of a symmetric group action  $\Sigma_n\times E_n\to E_n$ such that the composite bonding maps $(\P^1)^{\wedge p}\wedge E_n\to E_{n+p}$ are $\Sigma_p\times\Sigma_n$-equivariant. Forgetting the symmetric group actions defines the functor $\Spt^\Sigma_{\P^1}(k)\to \Spt_{\P^1}(k)$ which induces an equivalence on the homotopy categories $\SH^\Sigma(k)\to \SH(k)$. 

The advantage of the category $\Spt^\Sigma_{\P^1}(k)$ is that it admits a well-defined smash product, which gives  $\Spt^\Sigma_{\P^1}(k)$ the structure of a tensor category with unit the symmetric sphere spectrum $\1^\Sigma$. In particular, one has for each ring-object $\sE\in \Spt^\Sigma_{\P^1}(k)$ the catgory of $\sE$-modules in $\Spt^\Sigma_{\P^1}(k)$. Jardine also shows that the motivic cohomology spectrum $\HZ$ has a canonical lifting to a ring-object $\HZ\in \Spt^\Sigma_{\P^1}(k)$, giving the category of $\HZ$-modules. Recently, Ostvar-Roendigs \cite{OstvarRoendigs} have shown that the homotopy category of  the category of $\HZ$-modules is equivalent to the so-called ``big" category of motives over $k$, $\DM(k)$.

In this regard, we make the following working assumption or conjecture:

\begin{conj}\label{conj:ring} The functor $\oplus_n s_n$ and natural transformation $\mu_{*,*})$ on $\SH(k)$ lifts to a functor $\oplus_n s_n$ on $\Spt^\Sigma_{\P^1}(k)$ to graded objects in $\Spt^\Sigma_{\P^1}(k)$ and a natural transformation 
\[
\mu_{*,*}:\oplus_ns_n\wedge \oplus_ns_n\to \oplus_ns_n\circ \wedge.
\]
of bi-graded functors $\Spt^\Sigma_{\P^1}(k)\otimes \Spt^\Sigma_{\P^1}(k)\to \Spt^\Sigma_{\P^1}(k)$.
\end{conj}

Assuming this conjecture, we have

\begin{prop} \label{prop:module} For $\sE\in\SH(k)$,  $s_n\sE$ has a canonical structure of an  $s_0\1$-module.
\end{prop}

\begin{proof} We may represent $\sE$ by $\sE\in \Spt^\Sigma_{\P^1}(k)$.  The result then follows from the above discussion applied to the canonical $\1^\Sigma$-module structure on $\sE$.
\end{proof}

\subsection{$\HZ$-modules and $\DM$}

Following \cite{OstvarRoendigs}, we have the big triangulated tensor category of motives $\DM(k)$, containing the triangulated category of effective motives $\DM^\eff_-(k)$. There is an Eilenberg-Maclane functor
\[
\sH:\DM(k)\to \SH(k)
\]
sending $\DM^\eff_-(k)$ to $\SH^\eff(k)$, and a ``Suslin homology" functor
\[
h_S:  \SH(k)\to \DM(k),
\]
which is left adjoint to $\sH$ and sends $\SH^\eff(k)$ to $\DM^\eff_-(k)$. We denote the unit object of $\DM(k)$ by $\Z$. There are canonical isomorphisms $\sH(\Z)\cong \HZ$, $h_S(\1)\cong\Z$. 

\begin{thm}[Ostvar-Roendigs]\label{thm:OR} Via the equivalance of $\SH^\Sigma(k)$ with $\SH(k)$ given by the forgetful functor, the Eilenberg-Maclane functor identifies $\DM(k)$ with the homotopy category of $\HZ$-modules in $\Spt^\Sigma(k)$.
\end{thm}

This,  together with Theorem~\ref{thm:Slice} and Proposition~\ref{prop:module}, yields:

\begin{cor} For each $n$ there is a functor
\[
\sM(\barfil_n):\SH(k)\to \DM(k)
\]
and a natural isomorphism $\sH\circ\sM(\barfil_n)\cong \barfil_n$.
\end{cor}

\subsection{The spectral sequence}
\begin{Def} Let $\sE$ be a $\P^1$-$\Omega$-spectrum. Define the object $\pi^\mu_p\sE$ of $\DM(k)$ by 
\[
\pi^\mu_p\sE:=\sM(\barfil_p\sE)[-p].
\]
\end{Def}

We have the functor
\begin{align*}
m:\Sm/k\to\DM(k),\\
m(X)=h_S(\Omega_{\P^1}^\infty\Sigma^\infty_{\P^1}X_+).
\end{align*}
 
For an object $M$ of $\DM(k)$ and $X\in \Sm/k$, we have the {\em motivic cohomology}
\[
\H^n(X,M):=\Hom_{\DM}(m(X), M[n])
\]
and the natural isomorphism
\[
\H^n(X,M)\cong \Hom_{\SH(k)}(\Omega_{\P^1}^\infty\Sigma^\infty_{\P^1}X_+,\Sigma^{-n}_s\sH(M)).
\]

\begin{thm} Let $\sE$ be a $\P^1$-$\Omega$-spectrum,  let $E:\Sm/k^\op\to \Spt$ be the 0th spectrum of $\sE$ and let $X$ be in $\Sm/k$. Then (after re-indexing) the homotopy coniveau spectral sequence for $E(X)$ is 
\[
E_2^{p,q}=\H^p(X,\pi^\mu_{-q}\sE)\Longrightarrow
\pi_{-p-q}\hat{E}(X).
\]
\end{thm}

\begin{proof}  The $E_1$-term is given by
\[
E_1^{p,q}=\pi_{-p-q}(E^{(p/p+1)}(X)).
\]
Also $E^{(p/p+1)}$ is the $0$th spectrum in $\barfil_p(\sE)$, so we have
\begin{align*}
\pi_{-p-q}(E^{(p/p+1)}(X)) 
&=\Hom_{\SH(k)}(\Sigma^{-p-q}_s\Sigma^\infty_{\P^1}X_+,
\barfil_p(\sE))\\
&=\Hom_{\DM(k)}(m(X)[-p-q], \pi^\mu_p(\sE)[p])\\
&=\H^{2p+q}(X,\pi^\mu_p(\sE)).
\end{align*}

As the transformation $(p,q)\mapsto (p+r,q-r+1)$ sends $(2p+q,p)$ to $(2p+q+r+1, p+r)$, we can reindex to form an $E_2$-spectral sequence by replacing $\H^{2p+q}$ with $\H^p$ and $\pi^\mu_p$ with $\pi^\mu_{-q}$:
\[
E_2^{p,q}:=\H^p(X,\pi^\mu_{-q}(\sE))\Longrightarrow \pi_{-p-q}\hat{E}(X).
\]
\end{proof}

To aid in concrete computations, we use:

\begin{lem}\label{lem:Shift} We have an isomorphism in $\DM(k)$:
\[
\pi^\mu_{p+q}\sE\cong \pi^\mu_q(\Sigma_{\P^1}^{-p}\sE)\otimes\Z(p)[p].
\]
\end{lem}

\begin{proof} The functor $\sH$ is a tensor functor and there is a canonical weak equivalence
\[
\sH(\Z(q)[2q]\otimes M)\sim \Sigma^q_{\P^1}\sH(M).
\]
where $\Z(q)$ is the Tate object in $\DM(k)$.
Thus, if $\sF$ is an $\HZ$-module, we have the canonical isomorphism
\[
\sM(\Sigma^q_{\P^1}\sF)\cong \Z(q)[2q]\otimes \sM(\sF).
\]

By Remark~\ref{rem:Suspension}, we have the canonical isomorphisms in $\DM(k)$
\begin{align*}
\pi^\mu_{p+q}\sE=\sM(\barfil_{p+q}\sE)[-p-q]&\cong \sM(\barfil_{p+q}(\Sigma^p_{\P^1}\Sigma^{-p}_{\P^1}\sE))[-p-q]\\
&\cong \sM(\Sigma^p_{\P^1}\barfil_q(\Sigma^{-p}_{\P^1}\sE))[-p-q]\\
&\cong \sM(\barfil_q(\Sigma^{-p}_{\P^1}\sE))\otimes\Z(p)[p-q]\\
&\cong \pi^\mu_q(\Sigma_{\P^1}^{-p}\sE)\otimes\Z(p)[p].
\end{align*}
\end{proof}

So, for example, if $\sE$ is the $K$-theory $\P^1$-$\Omega$-spectrum $\sK:=(K,K,\ldots)$, then 
$\Sigma^q_{\P^1}\sK=\sK$, and $\barfil_0\sK=\HZ$. Thus $\pi^\mu_0(\sK)=\Z$ and
\[
\pi^\mu_p(\sK)=\pi^\mu_0(\sK)\otimes\Z(p)[p]=\Z(p)[p].
\]

\begin{rem} This identity does not rely on Conjecture~\ref{conj:ring}, rather, we have the isomorphism
\[
\sH(\pi^\mu_p(\sK)[p])\cong \barfil_p\sK
\]
by direct computation.
\end{rem}

Thus, our $E_2$-spectral sequence is the Bloch-Lichtenbaum, Fried\-lander-Suslin spectral sequence  \cite{BL, FriedSus} 
\[
E_2^{p,q}=H^{p-q}(X,\Z(-q))\Longrightarrow
K_{-p-q}(X).
\]

\begin{rem}
The $\P^1$-$\Omega$-spectrum $\fil_0\sK$:
\[
\fil_0\sK=(K,K^{(1)},K^{(2)},\ldots)
\]
gives an explicit model for  $\P^1$-connected algebraic $K$-theory.
\end{rem}
 
As a second example, we recall that, for $\sE\in\SH(k)$, $X\in\Sm/k$, we have the bi-graded homotopy groups
\[
\sE_{a,b}(X):=\Hom_{\SH(k)}(\Sigma^{a-2b}\Sigma_{\P^1}^b\Sigma^\infty X_+,\sE).
\]
Letting $E(b)$ be the $0$th spectrum of   $\Sigma_{\P^1}^{b}\sE$, we thus have the identity
\[
\sE_{a,-b}(X)=\pi_{a+2b}(E(b)(X)).
\]
Thus we have the spectral sequence
\[
E_2^{p,q}=\H^p(X,\pi^\mu_{-q}(\Sigma_{\P^1}^{b}\sE))\Longrightarrow
\pi_{-p-q}(\hat{E}(b)(X))=\hat{\sE}_{-p-q-2b,-b}(X).
\]

Via Lemma~\ref{lem:Shift}, we have
\begin{align*}
\H^p(X,\pi^\mu_{-q}(\Sigma_{\P^1}^{b}\sE))&=
\H^p(X,\pi^\mu_{-b-q}(\sE)\otimes\Z(b)[b])\\&=
\H^{p+b}(X,\pi^\mu_{-b-q}(\sE)\otimes\Z(b)).
\end{align*}
Thus, making the translation $(p,q)\mapsto(p-b,q-b)$, we have the spectral sequence
\[
E_2^{p,q}=\H^p(X,\pi^\mu_{-q}(\sE)\otimes\Z(b))\Longrightarrow
\hat{\sE}^{p+q,b}(X).
\]

\end{document}